\title{TABLE OF  FAMILIES OF ALTERNATING KNOTS WITH THEIR CONWAY'S FUNCTION}
\author{E. PI\~NA \\
Departamento de F\'\i sica\\ Universidad Aut\'onoma Metropolitana - Iztapalapa, \\
P. O. Box 55 534 \ Mexico, D. F., 09340 Mexico \\
e-mail: pge@xanum.uam.mx}

\documentclass[12pt]{article}
\usepackage{graphicx}
\usepackage{psfrag}

\begin{document}
\date {Dans le centenaire d'Henri Poincar\'e 1854-1912}
\maketitle

\abstract{A table of the families of alternating knots formed by conways is presented. The Conway's function is shown with thw use of linear algebra
in terms of natural numbers, called conways, that represent the number of crossings along a direction, as it was used by J. Conway for the classification
of knots.

Colored figures and tangles show the parts of the knots or tangles with a definite handedness: all the colored parts of the knot family are associated
to a particular orientation. For example all the colored conways have a right hand screw thread, and all the white conways have the opposite handedness.
Figures for six conways were colored with two different colors for forty two families in order to show the dissection of the knot in two tangles
corresponding to a particular factorization of the Conway's function. The Conway's function of each family is expressed as the internal product of
two vectors corresponding to each of two colored family 2-tangles, and with a full factorization which is not unique.
}

\pagebreak

\section{INTRODUCTION}
This document presents a table of the 65 families of prime alternating knots (see links also) formed by conways, from one to six conways. The conway
is represented in these figures by an "eye" with a $a_j$ inside it. Such a conway represent a tangle formed by the twisting of two strands which in
the projection to form the figures has $a_j$ crossings. The $a_j$'s represent a natural number labeled with $j$. Each figure of whatever family in this
table has as a figure caption the Conway's function of the family. This terminology gives credit to the seminal paper presented in 1967 and published
in 1970 by J. Conway \cite{co}.

The Conway's function of a family of knots is a polynomial with coefficients one, formed by monomials which are products of conways. It correspond to
the sum of the absolute values of the coefficients of the Alexander polynomial which Conway highlights for each knot or link in his tabulation of the
simplest knots \cite{co}. The number of monomials of the Conway's function is the Conway's number of the seed, which is the knot corresponding to the
family that has $a_j =1$ for all the conways of the family.

This table reproduces the same figures in \cite{p1} where all the terminology not fully defined in the present introduction is explained and discussed.
The difference with this document and \cite{p1} is the emphasis because here we present only the figures of the different prime families of knots.
A detailed account is found in \cite{p1}. Another difference is that the Conway's functions in the figure caption of the figure of each different family
appears here in terms of an explicit matrix product of vectors and matrices associated to 2-tangles or 3-tangles, highlighting the stelar presence of
the matrix metric in the 2-tangle case,
\begin{equation}
M \equiv \left( \begin{array}{cc}
0 & 1 \\
1 & 0
\end{array} \right)
\end{equation}
and the metric matrix in the 3-tangle case
\begin{equation}
\left( \begin{array}{ccccc}
0 & 0 & 0 & 0 & 1 \\
0 & 0 & 1 & 1 & 0 \\
0 & 1 & 0 & 1 & 0 \\
0 & 1 & 1 & 0 & 0 \\
1 & 0 & 0 & 0 & 0
\end{array} \right) \, .
\end{equation}

The full explanation of the table requires reference to \cite {p1}, however I believe it is very useful to have in the web this document independently
of the multiple definitions and properties of the algebra in \cite{p1}. The factorizations presented in this table show in an evident form many
(but not all) of the possible representations of the Conway's function of a particular family.

The Conway's function may be presented in all the cases as the interior product of two vectors corresponding to 2-tangles. For 35 families formed by
6 conways we give the Conway function as an interior product of two 2-tangles formed each by 3 conways recognized by the colors red and green.
A 2-dimensional vector of Conway's functions is associated to each 2-tangle; the interior product is defined through the metric matrix $M$.

By using the colors it is interesting to note that the last families of the table are factorized in terms of 3-tangles because it is not possible to
separate those families in 2-tangles formed, at least one of them, by two conways. The families with three colors do not allow to dissect in two
2-tangles of three conways as the other families of six conways; these families corresponding to one particular seed could be separated either
with two colors (2-tangles formed by 2 and 4 conways) or with two colors (two 3-tangles formed by three conways).

Note that I have not developed very far the algebra of 3-tangles. The 5-dimensional vectors that appear in the factorizations of 3-tangles in this
table may be identified as only four different vectors.

Many factorizations in this table could be modified in multiple forms by using the fact that the matrices of the form
\begin{equation}
\left( \begin{array}{cc}
A_1 & B_1 \\
B_1 & 0
\end{array} \right)\, , \quad \quad \left( \begin{array}{cc}
A_2 & B_2 \\
B_2 & 0
\end{array} \right)
\end{equation}
commute when multiplied with the metric $M$ as it occurs in many terms of this table namely
\begin{equation}
\left( \begin{array}{cc}
A_1 & B_1 \\
B_1 & 0
\end{array} \right) M \left( \begin{array}{cc}
A_2 & B_2 \\
B_2 & 0
\end{array} \right) = \left( \begin{array}{cc}
A_2 & B_2 \\
B_2 & 0
\end{array} \right) M \left( \begin{array}{cc}
A_1 & B_1 \\
B_1 & 0
\end{array} \right)\, .
\end{equation}

The transformation of this equation by the $M$ matrix gives the equivalent expression which could be used to obtain other different factorizations
\begin{equation}
\left( \begin{array}{cc}
0 & B_1 \\
B_1 & A_1
\end{array} \right) M \left( \begin{array}{cc}
0 & B_2 \\
B_2 & A_2
\end{array} \right) = \left( \begin{array}{cc}
0 & B_2 \\
B_2 & A_2
\end{array} \right) M \left( \begin{array}{cc}
0 & B_1 \\
B_1 & A_1
\end{array} \right)\, .
\end{equation}

Notice that these properties are related to the existence of abelian groups
\begin{equation}
\left( \begin{array}{cc}
A_1 & B_1 \\
B_1 & 0
\end{array} \right) M \left( \begin{array}{cc}
A_2 & B_2 \\
B_2 & 0
\end{array} \right) = \left( \begin{array}{cc}
A_1 B_2 + A_2 B_1 & B_1 B_2 \\
B_1 B_2 & 0
\end{array} \right) \, ,
\end{equation}
and
\begin{equation}
\left( \begin{array}{cc}
0 & B_1 \\
B_1 & A_1
\end{array} \right) M \left( \begin{array}{cc}
0 & B_2 \\
B_2 & A_2
\end{array} \right) = \left( \begin{array}{cc}
0 & B_1 B_2 \\
B_1 B_2 & A_1 B_2 + A_2 B_1
\end{array} \right) \, .
\end{equation}

Particular cases of this property at the beginning and ending of the factorization are remarked. We list a few examples at the beginning, but transposing these matrix equations you find similar expressions at the end of a product of factors.
$$
\left( \begin{array}{cc}
A_1 & B_1
\end{array} \right) M \left( \begin{array}{cc}
0 & A_2 \\
A_2 & B_2
\end{array} \right) = \left( \begin{array}{cc}
0 & 1
\end{array} \right) \left( \begin{array}{cc}
0 & A_1 \\
A_1 & B_1
\end{array} \right) M \left( \begin{array}{cc}
0 & A_2 \\
A_2 & B_2
\end{array} \right) =
$$
\begin{equation}
\left( \begin{array}{cc}
A_2 & B_2
\end{array} \right) M \left( \begin{array}{cc}
0 & A_1 \\
A_1 & B_1
\end{array} \right)\, .
\end{equation}
and
$$
\left( \begin{array}{cc}
A_1 & B_1
\end{array} \right) M \left( \begin{array}{cc}
A_2 & B_2 \\
B_2 & 0
\end{array} \right) = \left( \begin{array}{cc}
1 & 0
\end{array} \right) \left( \begin{array}{cc}
A_1 & B_1 \\
B_1 & 0
\end{array} \right) M \left( \begin{array}{cc}
A_2 & B_2 \\
B_2 & 0
\end{array} \right) =
$$
\begin{equation}
\left( \begin{array}{cc}
A_2 & B_2
\end{array} \right) M \left( \begin{array}{cc}
A_1 & B_1 \\
B_1 & 0
\end{array} \right)\, .
\end{equation}

\begin{figure}[b]

\section{THE FAMILIES OF KNOTS OF ONE AND TWO CONWAYS}

\psfrag{p}{\LARGE{$a_1$}}
\psfrag{q}{\LARGE{$a_2$}}
\hfil\scalebox{0.5}{\includegraphics{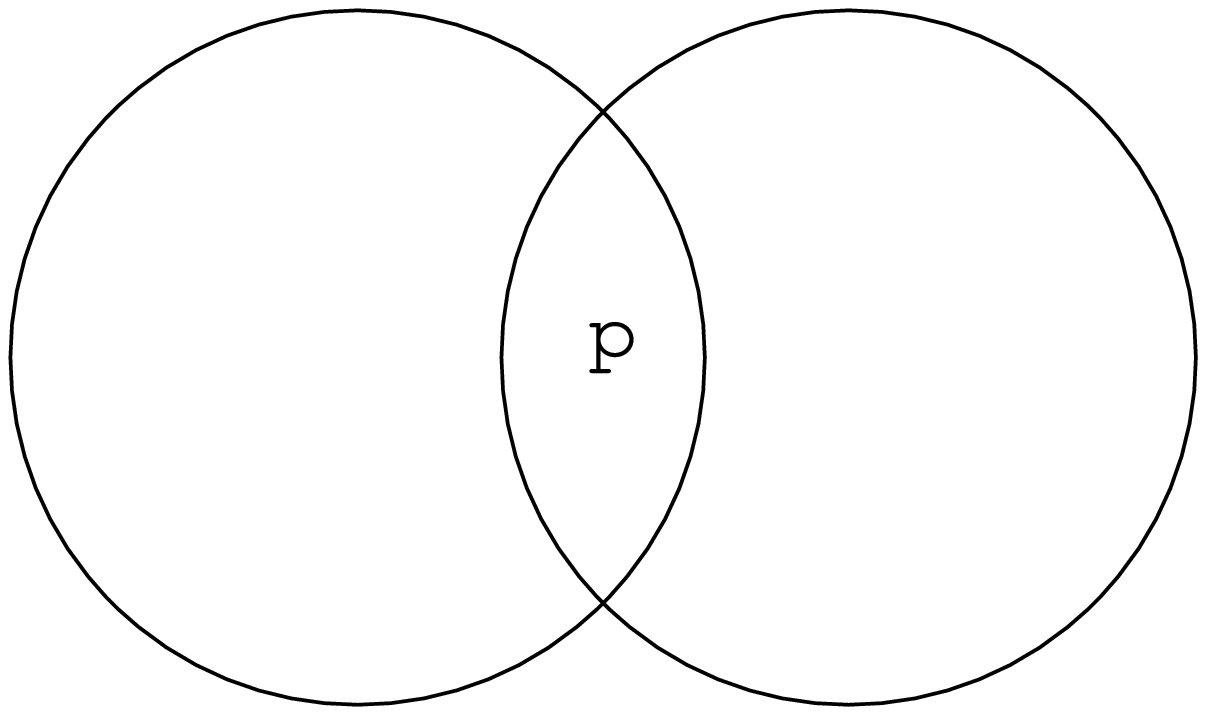}}\hfil\scalebox{0.5}{\includegraphics{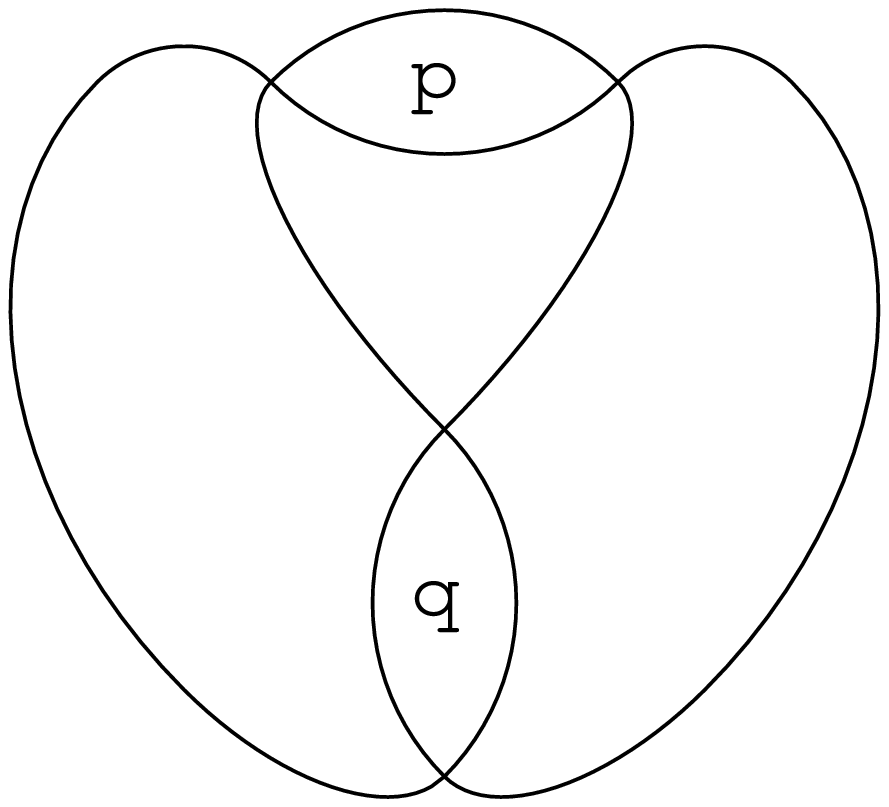}}\hfil

\

\caption{ \ }
Rational knots, its Conway's functions are $a_1$ and $1 + a_1 a_2$.

\

\end{figure}

\newpage

\begin{figure}

\section*{THE FAMILIES OF KNOTS OF THREE CONWAYS WITH SEED THE TREFOIL TORUS KNOT $3_1$ (2 CASES)}

\centering

\psfrag{a}{\LARGE{$a_1$}}
\psfrag{b}{\LARGE{$a_2$}}
\psfrag{c}{\LARGE{$a_3$}}
\psfrag{d}{\LARGE{$a_4$}}
\scalebox{0.50}{\includegraphics{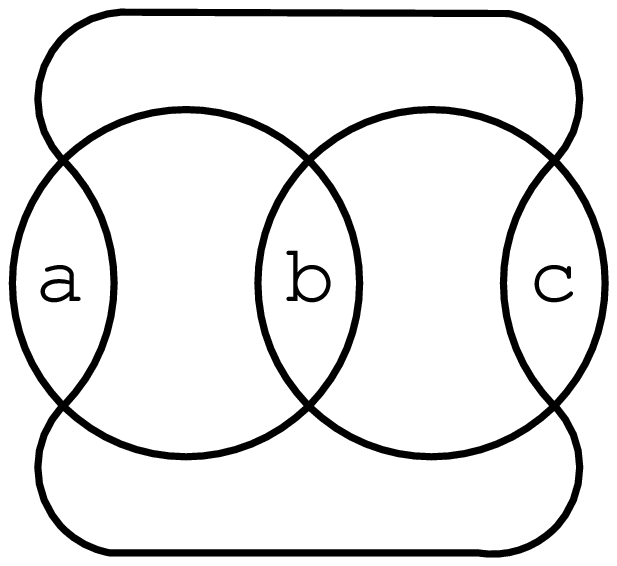}}

\caption{ \ }
$$
a_1 a_2 + a_2 a_3 + a_3 a_1 = (a_1, 1) M \left( \begin{array}{cc}
0 & a_2 \\
a_2 & 1
\end{array} \right) M \left( \begin{array}{c}
a_3 \\
1
\end{array} \right)
$$
\end{figure}

\begin{figure}

\centering

\psfrag{a}{\LARGE{$a_1$}}
\psfrag{b}{\LARGE{$a_2$}}
\psfrag{c}{\LARGE{$a_3$}}
\psfrag{d}{\LARGE{$a_4$}}
\scalebox{0.50}{\includegraphics{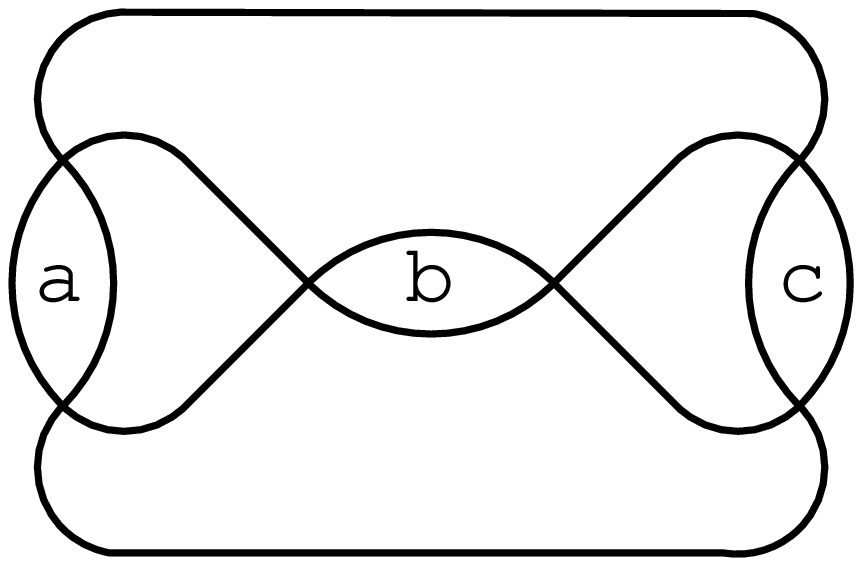}}

\caption{ \ }
$$
a_1 a_2 a_3 + a_1 + a_3 = (a_1, 1) M \left( \begin{array}{cc}
0 & 1 \\
1 & a_2
\end{array} \right) M \left( \begin{array}{c}
a_3 \\
1
\end{array} \right)
$$
\end{figure}

\begin{figure}

\section{THE FAMILIES OF KNOTS OF FOUR CONWAYS (5 CASES)}

\subsection{The families of knots with seed the Solomon link (two cases)}

\centering

\psfrag{a}{\LARGE{$a_1$}}
\psfrag{b}{\LARGE{$a_2$}}
\psfrag{c}{\LARGE{$a_3$}}
\psfrag{d}{\LARGE{$a_4$}}
\scalebox{0.50}{\includegraphics{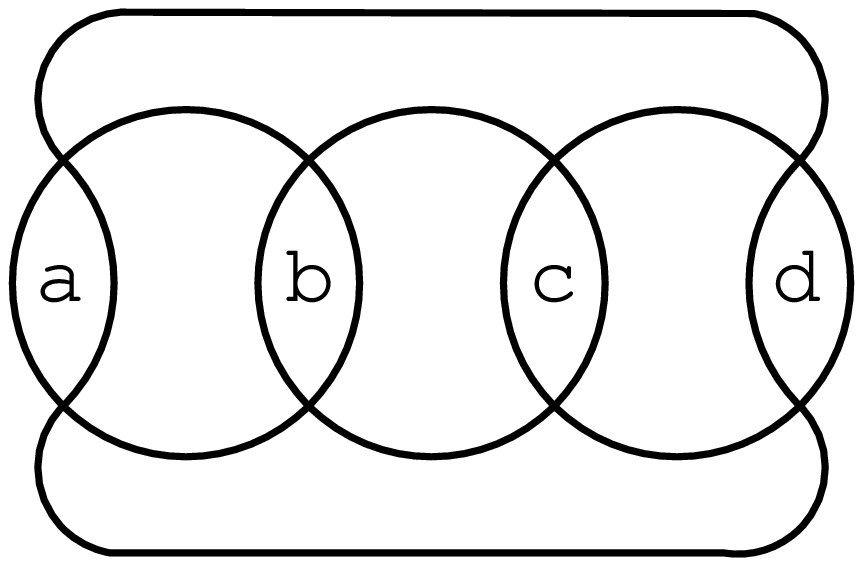}}

\caption{ \ }
$$
(a_1 a_2,  a_1 + a_2) M \left( \begin{array}{c}
a_3 a_4 \\
a_3 + a_4
\end{array} \right) =
$$
$$
(a_1, 1) M \left( \begin{array}{cc}
0 & a_2 \\
a_2 & 1
\end{array} \right) M \left( \begin{array}{cc}
0 & a_3 \\
a_3 & 1
\end{array} \right) M \left( \begin{array}{c}
a_4 \\
1
\end{array} \right)
$$
\end{figure}

\begin{figure}

\centering

\psfrag{a}{\LARGE{$a_1$}}
\psfrag{b}{\LARGE{$a_2$}}
\psfrag{c}{\LARGE{$a_3$}}
\psfrag{d}{\LARGE{$a_4$}}
\scalebox{0.50}{\includegraphics{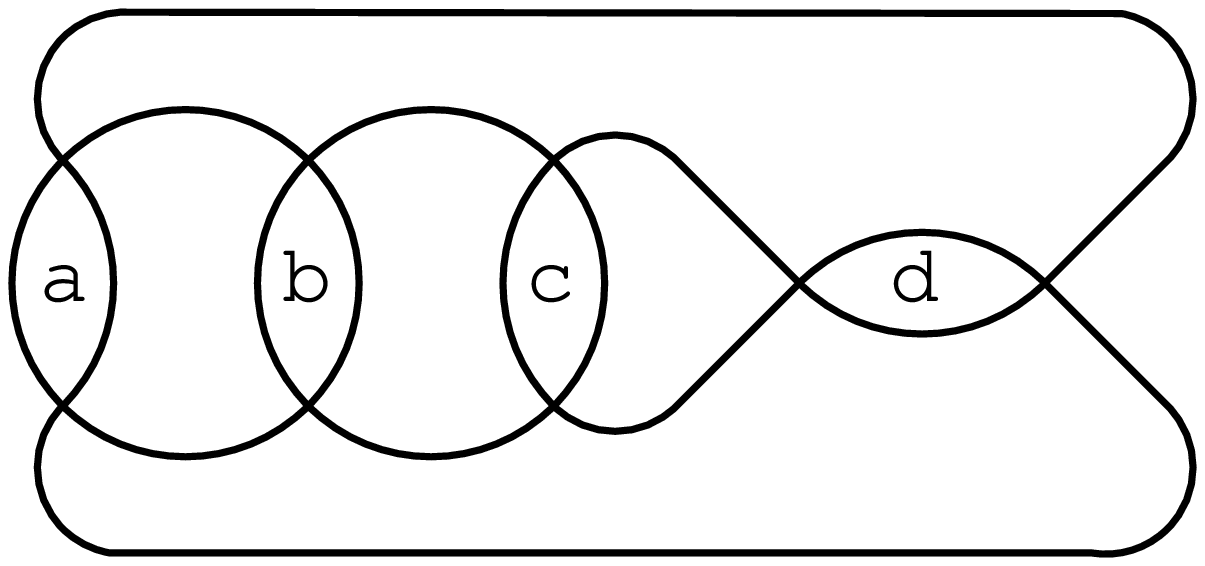}}

\caption{ \ }
$$
(a_1 a_2,  a_1 + a_2) M \left( \begin{array}{c}
a_3 \\
a_3 a_4 + 1
\end{array} \right) =
$$
$$
(a_1, 1) M \left( \begin{array}{cc}
0 & a_2 \\
a_2 & 1
\end{array} \right) M \left( \begin{array}{cc}
0 & a_3 \\
a_3 & 1
\end{array} \right) M \left( \begin{array}{c}
1 \\
a_4
\end{array} \right)
$$

\end{figure}

\begin{figure}

\subsection{The families of knots with seed the knot $4_1$ (three cases)}

\centering

\psfrag{a}{\LARGE{$a_1$}}
\psfrag{b}{\LARGE{$a_2$}}
\psfrag{c}{\LARGE{$a_3$}}
\psfrag{d}{\LARGE{$a_4$}}
\scalebox{0.50}{\includegraphics{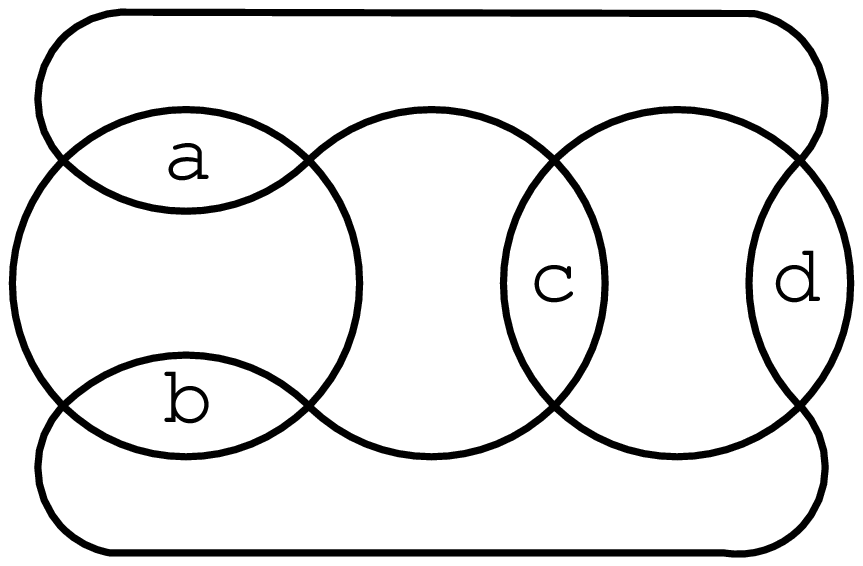}}

\caption{ \ }
$$
(a_1 + a_2,  a_1 a_2) M \left( \begin{array}{c}
a_3 a_4\\
a_3 + a_4
\end{array} \right) =
$$
$$
(1, a_1) M \left( \begin{array}{cc}
1 & a_2 \\
a_2 & 0
\end{array} \right) M \left( \begin{array}{cc}
0 & a_3 \\
a_3 & 1
\end{array} \right) M \left( \begin{array}{c}
a_4 \\
1
\end{array} \right)
$$
\end{figure}

\begin{figure}

\centering

\psfrag{a}{\LARGE{$a_1$}}
\psfrag{b}{\LARGE{$a_2$}}
\psfrag{c}{\LARGE{$a_3$}}
\psfrag{d}{\LARGE{$a_4$}}
\scalebox{0.50}{\includegraphics{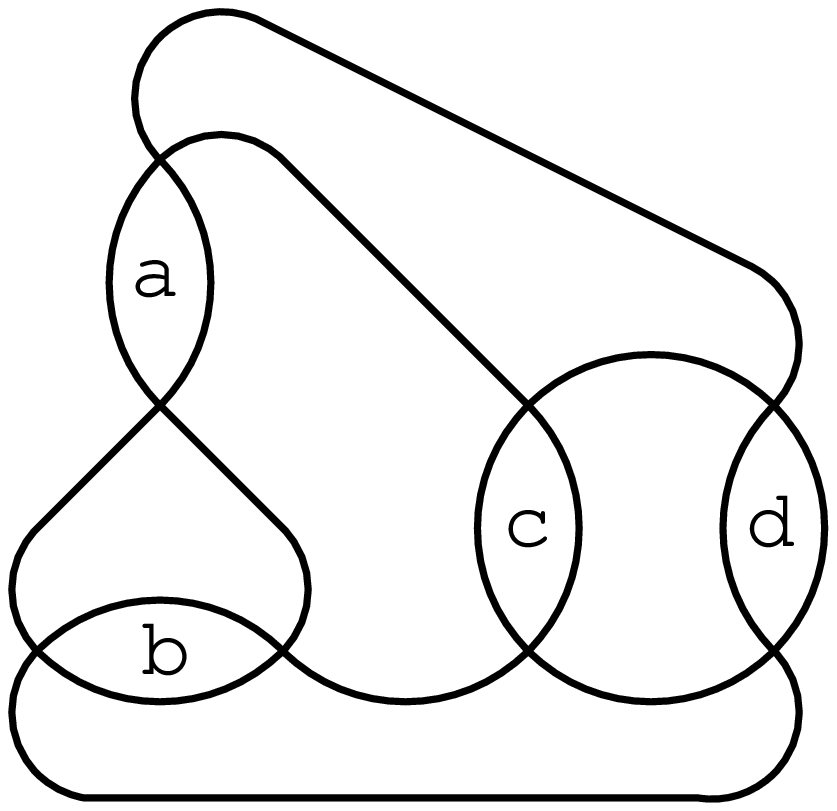}}

\caption{ \ }
$$
(a_1 a_2 + 1, a_2) M \left( \begin{array}{c}
a_3 a_4\\
a_3 + a_4
\end{array} \right) =
$$
$$
(a_1, 1) M \left( \begin{array}{cc}
1 & a_2 \\
a_2 & 0
\end{array} \right) M \left( \begin{array}{cc}
0 & a_3 \\
a_3 & 1
\end{array} \right) M \left( \begin{array}{c}
a_4 \\
1
\end{array} \right)
$$
\end{figure}

\begin{figure}

\centering

\psfrag{a}{\LARGE{$a_1$}}
\psfrag{b}{\LARGE{$a_2$}}
\psfrag{c}{\LARGE{$a_3$}}
\psfrag{d}{\LARGE{$a_4$}}
\scalebox{0.50}{\includegraphics{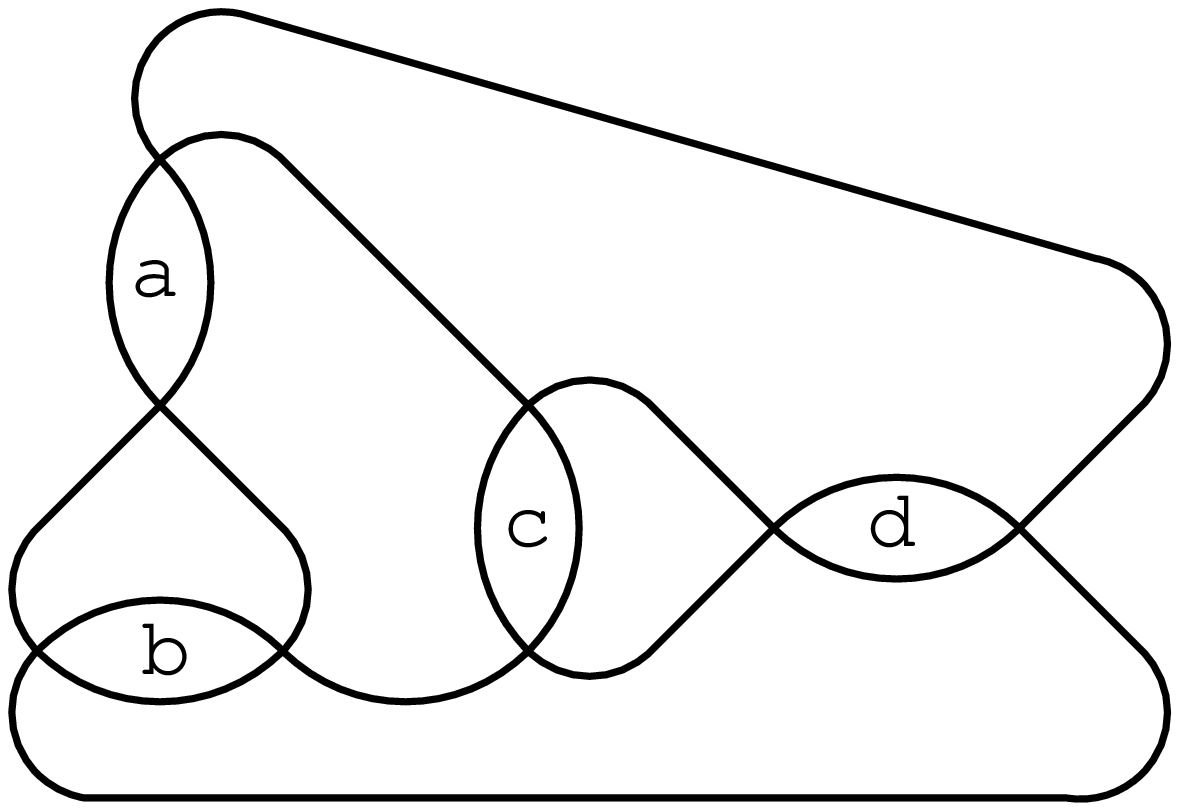}}

\caption{ \ }
$$
(a_1 a_2 + 1, a_2) M \left( \begin{array}{c}
a_3 \\
a_3 a_4 + 1
\end{array} \right) =
$$
$$
(a_1, 1) M \left( \begin{array}{cc}
1 & a_2 \\
a_2 & 0
\end{array} \right) M \left( \begin{array}{cc}
0 & a_3 \\
a_3 & 1
\end{array} \right) M \left( \begin{array}{c}
1 \\
a_4
\end{array} \right)
$$
\end{figure}

\begin{figure}

\section{THE FAMILIES OF KNOTS OF FIVE CONWAYS (12 CASES)}

\subsection{The families of knots with seed the torus knot $5_1$ (two cases)}

\centering

\psfrag{a}{\LARGE{$a_1$}}
\psfrag{b}{\LARGE{$a_2$}}
\psfrag{c}{\LARGE{$a_3$}}
\psfrag{d}{\LARGE{$a_4$}}
\psfrag{e}{\LARGE{$a_5$}}
\scalebox{0.50}{\includegraphics{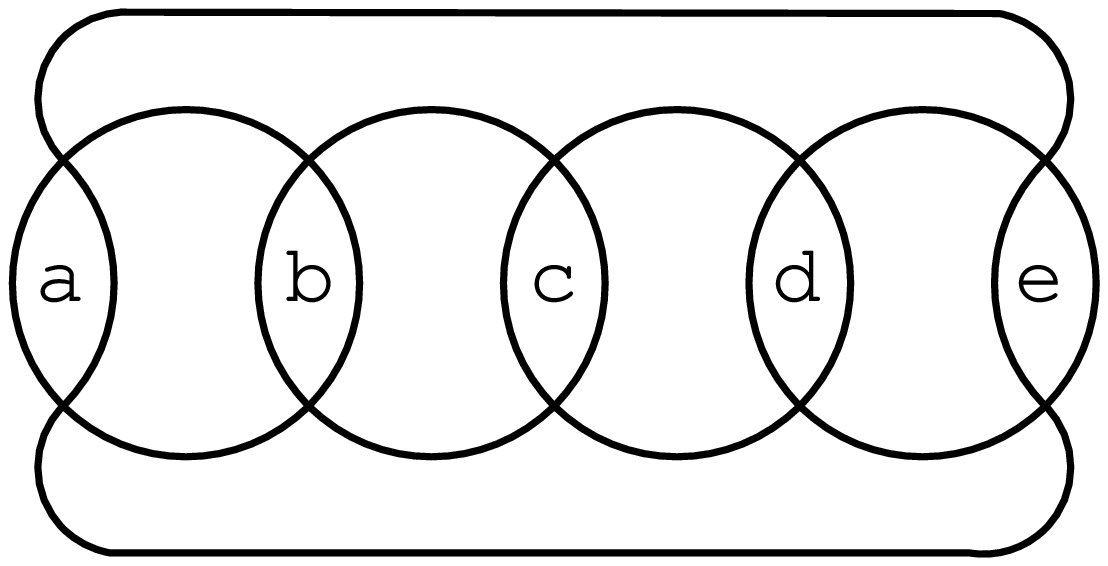}}

\caption{ \ }
$$
(a_1 a_2 a_3,  a_1 a_2 + a_2 a_3 + a_3 a_1) M \left( \begin{array}{c}
a_4 a_5 \\
a_4 + a_5
\end{array} \right) = (a_1 a_2,  a_1 + a_2) M \left( \begin{array}{c}
a_3 a_4 a_5\\
a_3 a_4 + a_4 a_5 + a_5 a_3
\end{array} \right) =
$$
$$
(a_1, 1) M \left( \begin{array}{cc}
0 & a_2 \\
a_2 & 1
\end{array} \right) M \left( \begin{array}{cc}
0 & a_3 \\
a_3 & 1
\end{array} \right) M \left( \begin{array}{cc}
0 & a_4 \\
a_4 & 1
\end{array} \right) M \left( \begin{array}{c}
a_5 \\
1
\end{array} \right)
$$
\end{figure}
\begin{figure}

\centering

\psfrag{a}{\LARGE{$a_1$}}
\psfrag{b}{\LARGE{$a_2$}}
\psfrag{c}{\LARGE{$a_3$}}
\psfrag{d}{\LARGE{$a_4$}}
\psfrag{e}{\LARGE{$a_5$}}
\scalebox{0.50}{\includegraphics{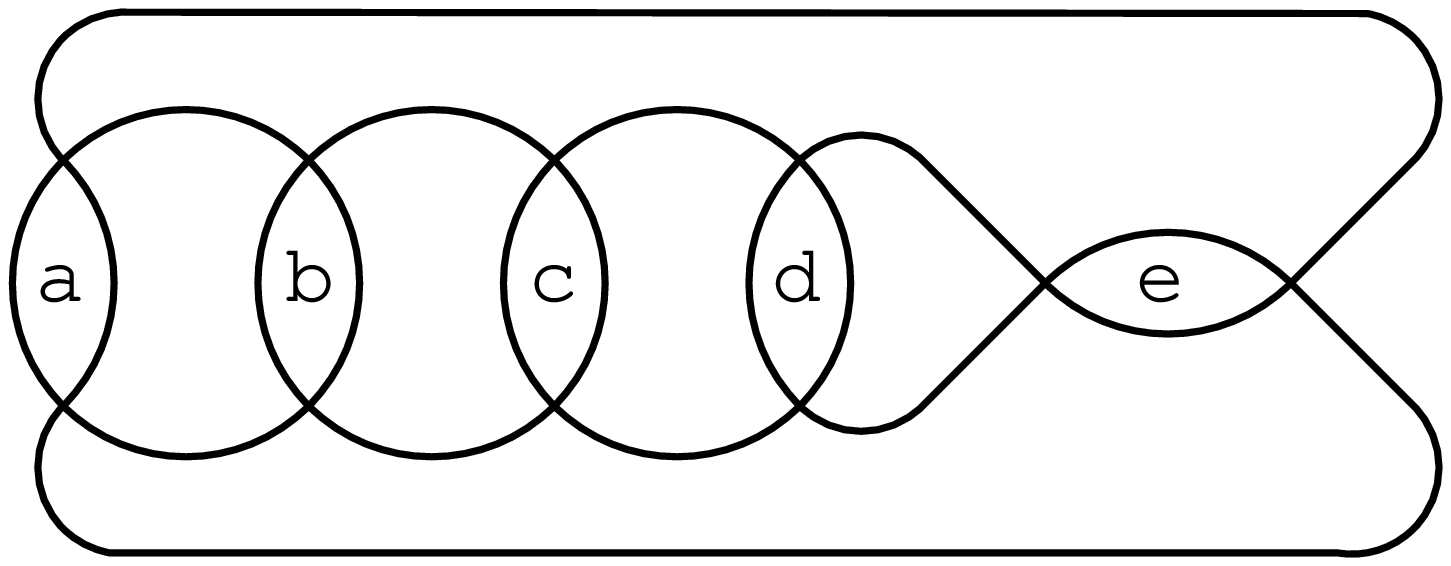}}

\caption{ \ }
$$
(a_1 a_2 a_3,  a_1 a_2 + a_2 a_3 + a_3 a_1) M \left( \begin{array}{c}
a_4 \\
a_4 a_5 + 1
\end{array} \right) = (a_1 a_2,  a_1 + a_2) M \left( \begin{array}{c}
a_3 a_4 \\
a_3 a_4 a_5 +  a_3 + a_4
\end{array} \right) =
$$
$$
(a_1, 1) M \left( \begin{array}{cc}
0 & a_2 \\
a_2 & 1
\end{array} \right) M \left( \begin{array}{cc}
0 & a_3 \\
a_3 & 1
\end{array} \right) M \left( \begin{array}{cc}
0 & a_4 \\
a_4 & 1
\end{array} \right) M \left( \begin{array}{c}
1 \\
a_5
\end{array} \right)
$$
\end{figure}
\begin{figure}

\subsection{The families of knots with seed the knot $5_2$ (four cases)}

\centering

\psfrag{a}{\LARGE{$a_1$}}
\psfrag{b}{\LARGE{$a_2$}}
\psfrag{c}{\LARGE{$a_3$}}
\psfrag{d}{\LARGE{$a_4$}}
\psfrag{e}{\LARGE{$a_5$}}
\scalebox{0.50}{\includegraphics{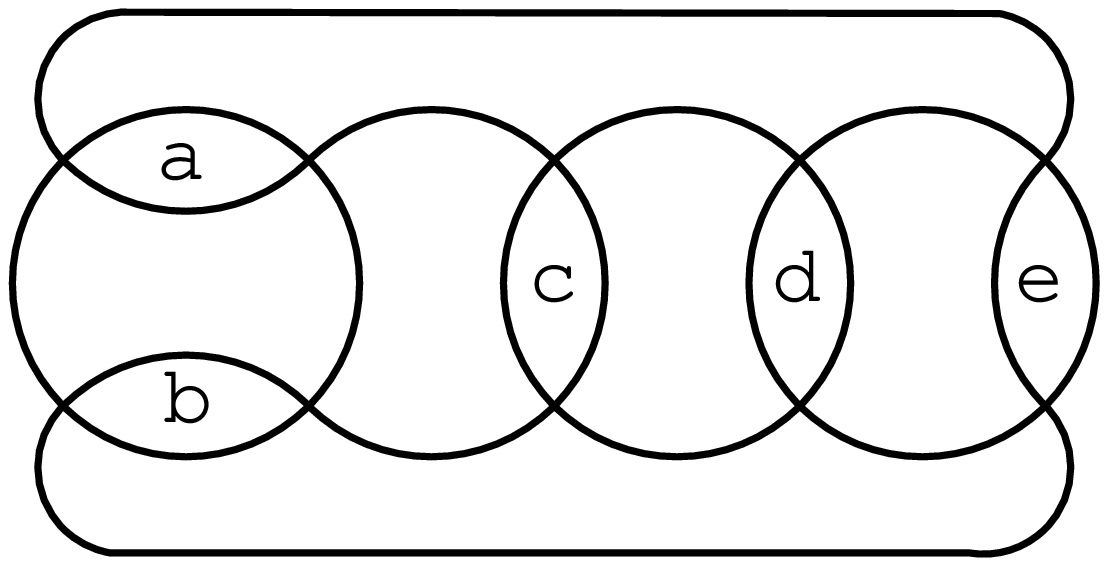}}

\caption{ \ }
$$
(a_1 + a_2,  a_1 a_2) M \left( \begin{array}{c}
a_3 a_4 a_5\\
a_3 a_4 + a_4 a_5 + a_5 a_3
\end{array} \right) = ((a_1 + a_2) a_3,  a_1 a_2 a_3 + a_1 + a_2) M \left( \begin{array}{c}
a_4 a_5 \\
a_4 + a_5
\end{array} \right) =
$$
$$
(1, a_1) M \left( \begin{array}{cc}
1 & a_2 \\
a_2 & 0
\end{array} \right) M \left( \begin{array}{cc}
0 & a_3 \\
a_3 & 1
\end{array} \right) M \left( \begin{array}{cc}
0 & a_4 \\
a_4 & 1
\end{array} \right) M \left( \begin{array}{c}
a_5 \\
1
\end{array} \right)
$$
\end{figure}

\begin{figure}

\centering

\psfrag{a}{\LARGE{$a_1$}}
\psfrag{b}{\LARGE{$a_2$}}
\psfrag{c}{\LARGE{$a_3$}}
\psfrag{d}{\LARGE{$a_4$}}
\psfrag{e}{\LARGE{$a_5$}}
\scalebox{0.50}{\includegraphics{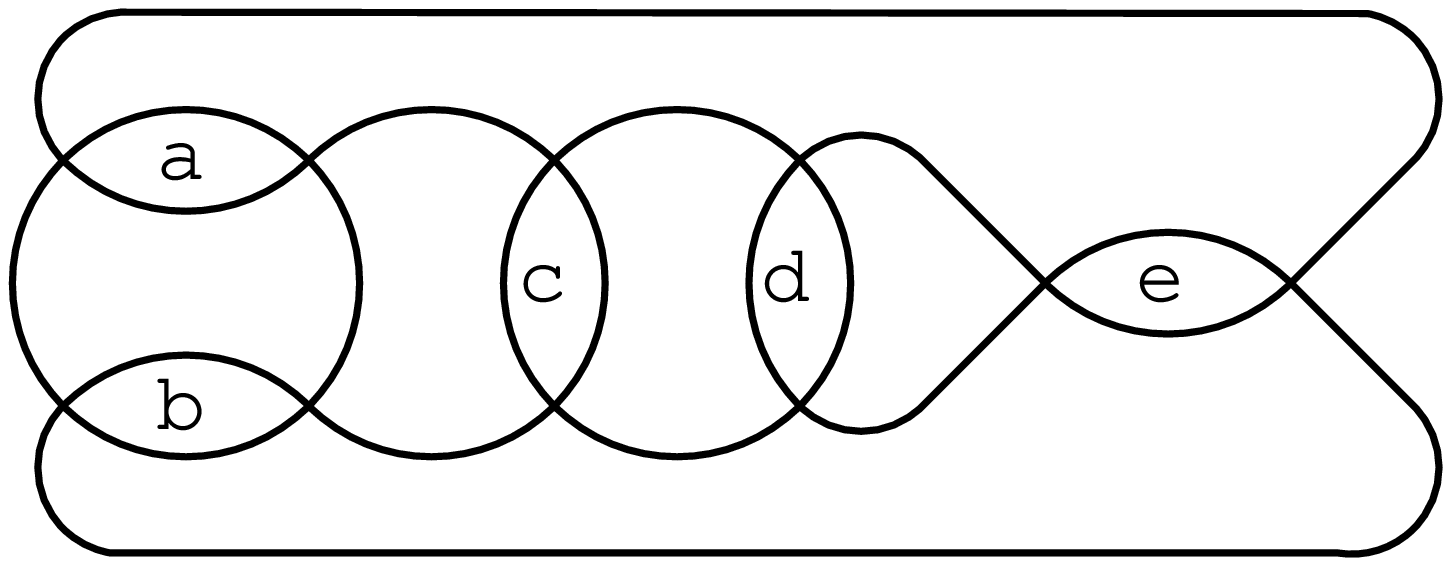}}

\caption{ \ }
$$
((a_1 + a_2) a_3,  a_1 a_2 a_3 + a_1 + a_2) M \left( \begin{array}{c}
a_4 \\
a_4 a_5 + 1
\end{array} \right) = (a_1 + a_2,  a_1 a_2) M \left( \begin{array}{c}
a_3 a_4 \\
a_3 a_4 a_5 + a_3 + a_4
\end{array} \right) =
$$
$$
(1, a_1) M \left( \begin{array}{cc}
1 & a_2 \\
a_2 & 0
\end{array} \right) M \left( \begin{array}{cc}
0 & a_3 \\
a_3 & 1
\end{array} \right) M \left( \begin{array}{cc}
0 & a_4 \\
a_4 & 1
\end{array} \right) M \left( \begin{array}{c}
1 \\
a_5
\end{array} \right)
$$
\end{figure}

\begin{figure}

\centering

\psfrag{a}{\LARGE{$a_1$}}
\psfrag{b}{\LARGE{$a_2$}}
\psfrag{c}{\LARGE{$a_3$}}
\psfrag{d}{\LARGE{$a_4$}}
\psfrag{e}{\LARGE{$a_5$}}
\scalebox{0.50}{\includegraphics{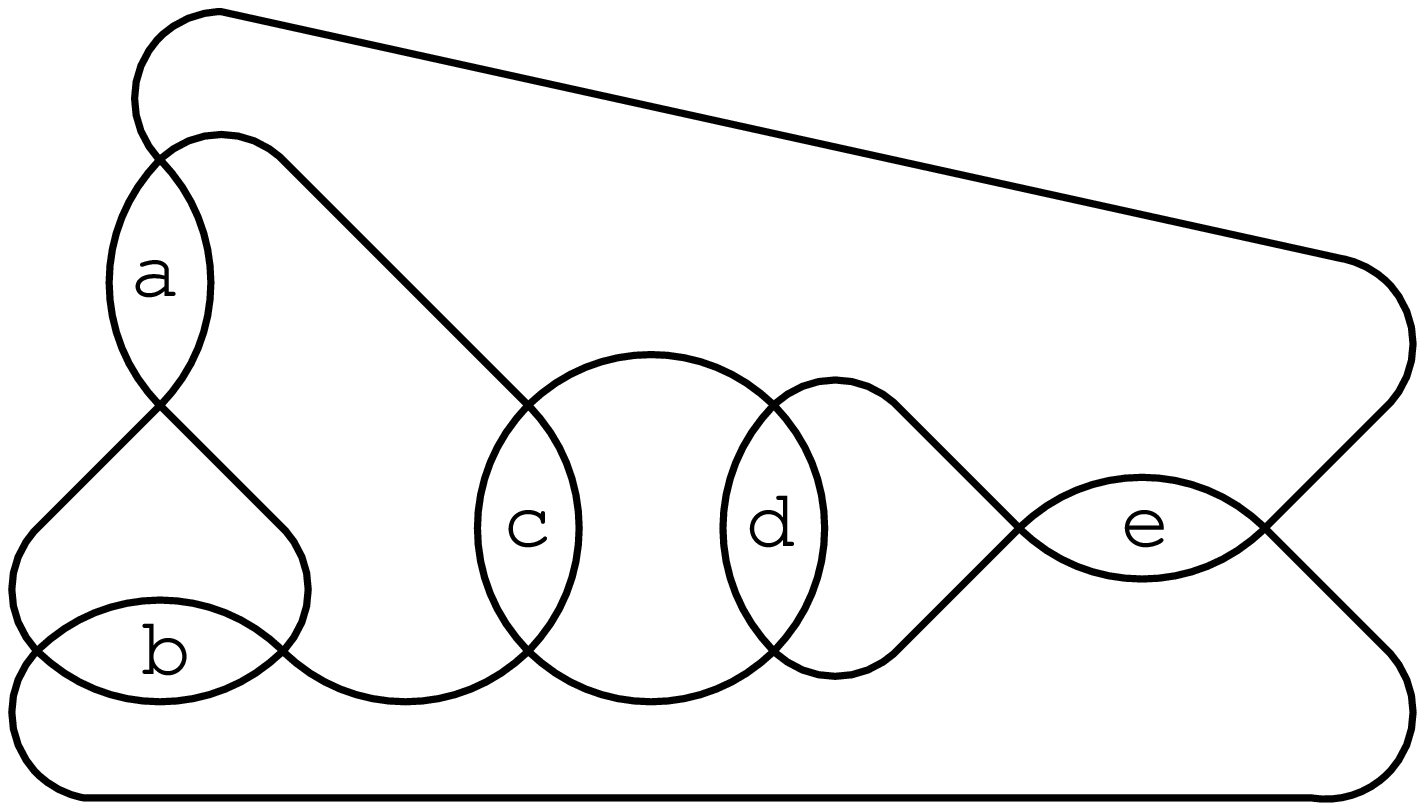}}

\caption{ \ }
$$
(a_1 a_2 + 1, a_2) M \left( \begin{array}{c}
a_3 a_4 \\
a_3 a_4 a_5 + a_3 + a_4
\end{array} \right) = (a_1 a_2 a_3 + a_3, a_2 (a_1 + a_3) + 1) M \left( \begin{array}{c}
a_4 \\
a_4 a_5 + 1
\end{array} \right) =
$$
$$
(a_1, 1) M \left( \begin{array}{cc}
1 & a_2 \\
a_2 & 0
\end{array} \right) M \left( \begin{array}{cc}
0 & a_3 \\
a_3 & 1
\end{array} \right) M \left( \begin{array}{cc}
0 & a_4 \\
a_4 & 1
\end{array} \right) M \left( \begin{array}{c}
1 \\
a_5
\end{array} \right)
$$
\end{figure}

\begin{figure}

\centering

\psfrag{a}{\LARGE{$a_1$}}
\psfrag{b}{\LARGE{$a_2$}}
\psfrag{c}{\LARGE{$a_3$}}
\psfrag{d}{\LARGE{$a_4$}}
\psfrag{e}{\LARGE{$a_5$}}
\scalebox{0.50}{\includegraphics{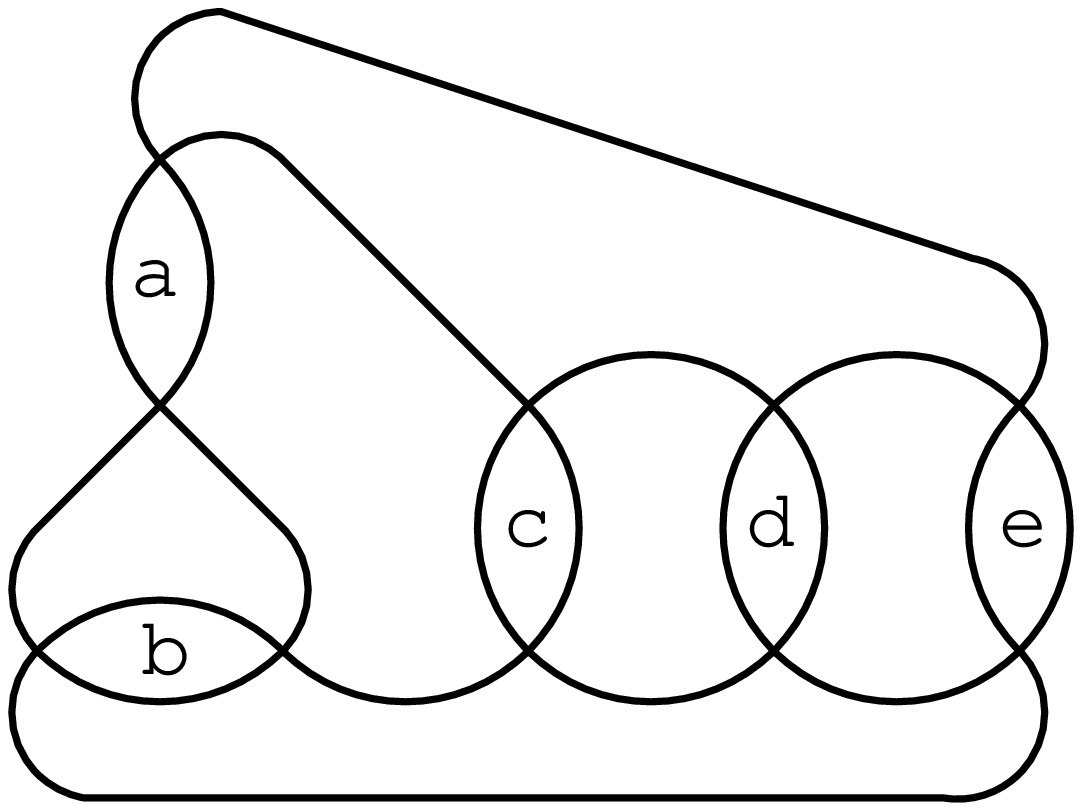}}

\caption{ \ }
$$
(a_1 a_2 + 1, a_2) M \left( \begin{array}{c}
a_3 a_4 a_5 \\
a_3 a_4 + a_4 a_5 + a_5 a_3
\end{array} \right) = (a_1 a_2 a_3 + a_3, a_2 (a_1 + a_3) + 1) M \left( \begin{array}{c}
a_4 a_5\\
a_4 + a_5
\end{array} \right) =
$$
$$
(a_1, 1) M \left( \begin{array}{cc}
1 & a_2 \\
a_2 & 0
\end{array} \right) M \left( \begin{array}{cc}
0 & a_3 \\
a_3 & 1
\end{array} \right) M \left( \begin{array}{cc}
0 & a_4 \\
a_4 & 1
\end{array} \right) M \left( \begin{array}{c}
a_5 \\
1
\end{array} \right)
$$
\end{figure}

\begin{figure}

\subsection{The families of knots with seed the Whitehead link $5_1^2$ (six cases)}

\centering

\psfrag{a}{\LARGE{$a_1$}}
\psfrag{b}{\LARGE{$a_4$}}
\psfrag{c}{\LARGE{$a_3$}}
\psfrag{d}{\LARGE{$a_2$}}
\psfrag{e}{\LARGE{$a_5$}}
\scalebox{0.45}{\includegraphics{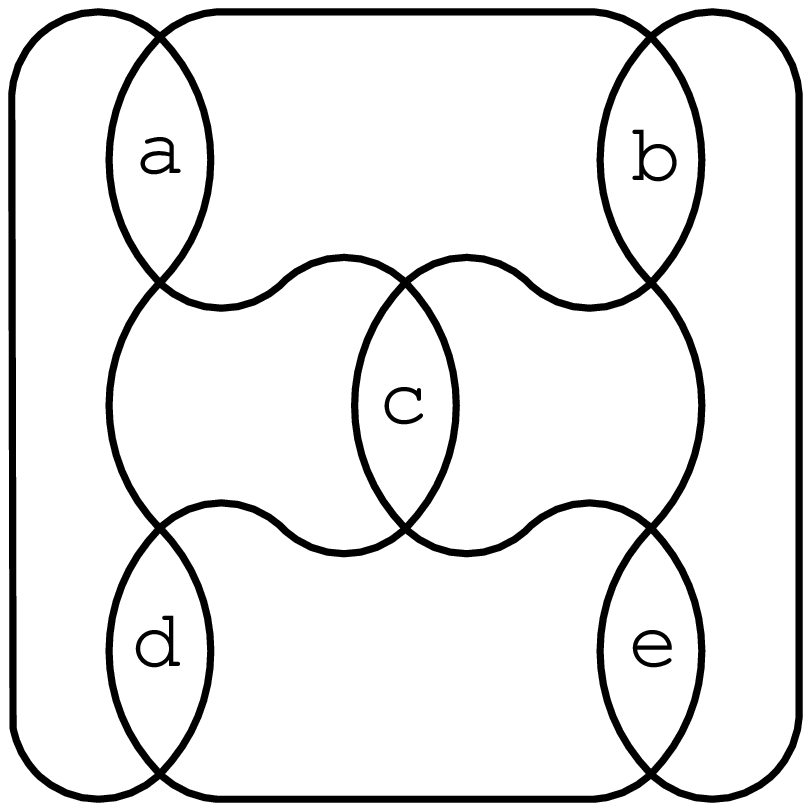}}

\caption{ \ }
$$
(a_1 + a_2, a_1 a_2) M \left( \begin{array}{c}
a_3 (a_4 + a_5) \\
a_3 a_4 a_5 + a_4 + a_5
\end{array} \right) =
$$
$$
(1, a_1) M \left( \begin{array}{cc}
1 & a_2 \\
a_2 & 0
\end{array} \right) M \left( \begin{array}{cc}
0 & a_3 \\
a_3 & 1
\end{array} \right) M \left( \begin{array}{cc}
1 & a_4 \\
a_4 & 0
\end{array} \right) M \left( \begin{array}{c}
1 \\
a_5
\end{array} \right)
$$
\end{figure}

\begin{figure}

\centering

\psfrag{a}{\LARGE{$a_1$}}
\psfrag{b}{\LARGE{$a_4$}}
\psfrag{c}{\LARGE{$a_3$}}
\psfrag{d}{\LARGE{$a_2$}}
\psfrag{e}{\LARGE{$a_5$}}
\scalebox{0.45}{\includegraphics{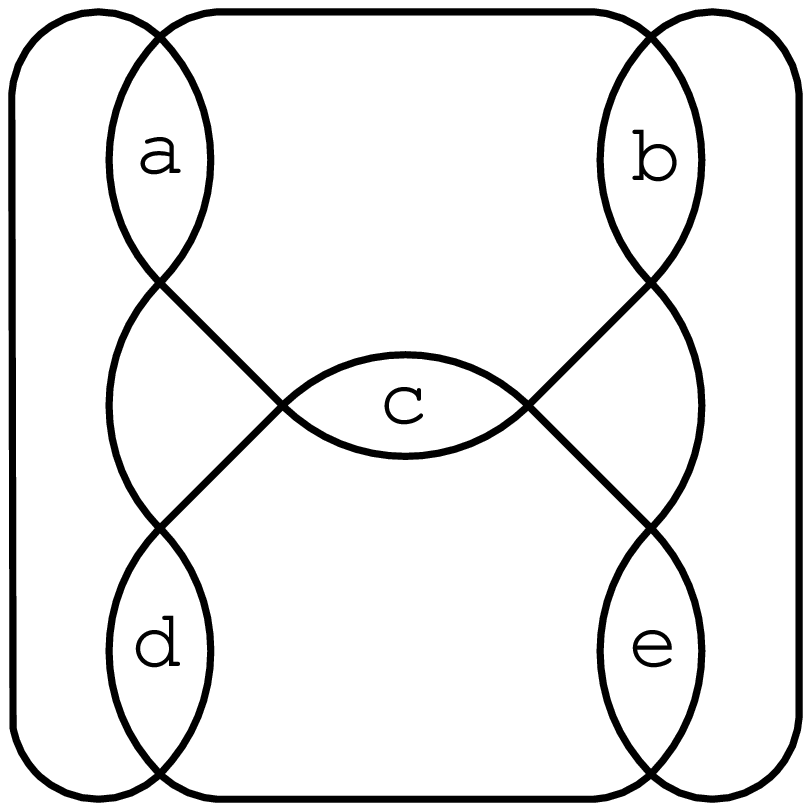}}

\caption{ \ }
$$
(a_1 + a_2, a_1 a_2) M \left( \begin{array}{c}
a_4 + a_5 \\
a_3 a_4 + a_4 a_5 + a_5 a_4
\end{array} \right) =
$$
$$
(1, a_1) M \left( \begin{array}{cc}
1 & a_2 \\
a_2 & 0
\end{array} \right) M \left( \begin{array}{cc}
0 & 1 \\
1 & a_3
\end{array} \right) M \left( \begin{array}{cc}
1 & a_4 \\
a_4 & 0
\end{array} \right) M \left( \begin{array}{c}
1 \\
a_5
\end{array} \right)
$$
\end{figure}

\begin{figure}

\centering

\psfrag{a}{\LARGE{$a_1$}}
\psfrag{b}{\LARGE{$a_4$}}
\psfrag{c}{\LARGE{$a_3$}}
\psfrag{d}{\LARGE{$a_2$}}
\psfrag{e}{\LARGE{$a_5$}}
\scalebox{0.50}{\includegraphics{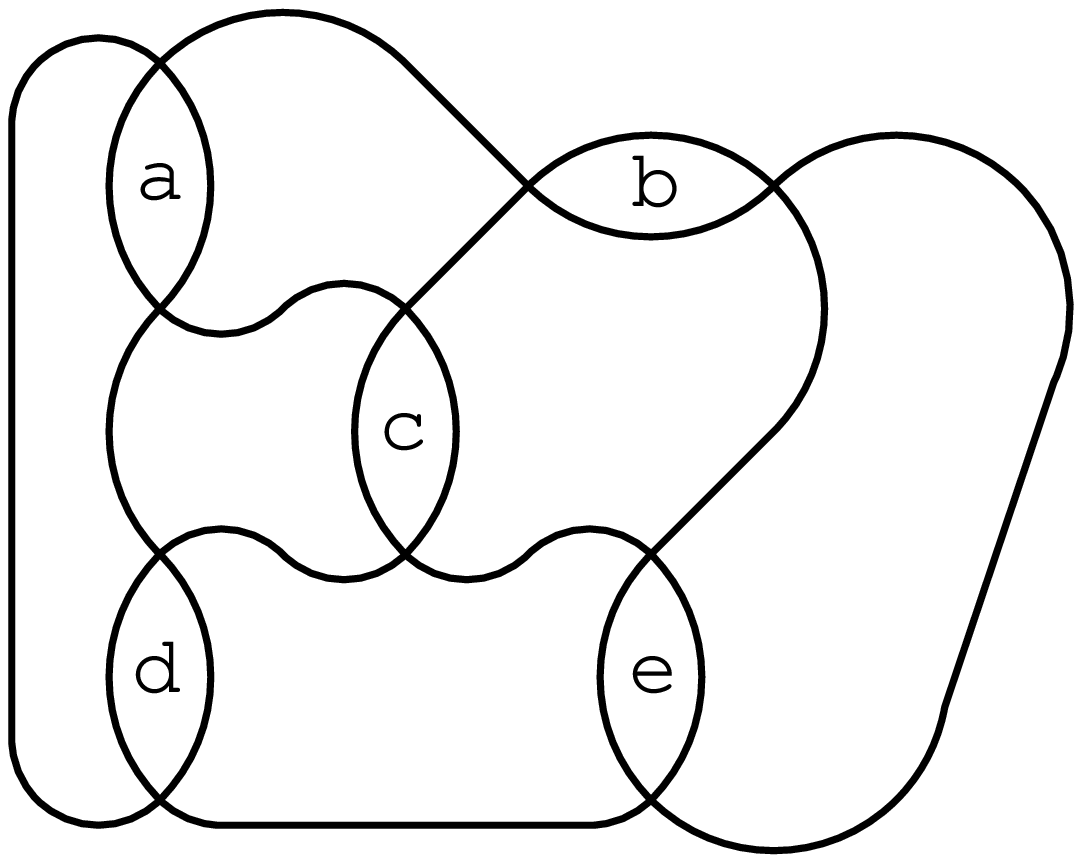}}

\caption{ \ }
$$
(a_1 + a_2, a_1 a_2) M \left( \begin{array}{c}
a_3 (a_4 a_5 + 1) \\
(a_3 (a_4 + a_5) + 1
\end{array} \right) = ((a_1 + a_2) a_3, a_1 a_2 a_3 + a_1 + a_2) M \left( \begin{array}{c}
a_4 a_5 + 1 \\
a_5
\end{array} \right) =
$$
$$
(1, a_1) M \left( \begin{array}{cc}
1 & a_2 \\
a_2 & 0
\end{array} \right) M \left( \begin{array}{cc}
0 & a_3 \\
a_3 & 1
\end{array} \right) M \left( \begin{array}{cc}
a_4 & 1 \\
1 & 0
\end{array} \right) M \left( \begin{array}{c}
1 \\
a_5
\end{array} \right)
$$
\end{figure}

\begin{figure}

\centering

\psfrag{a}{\LARGE{$a_1$}}
\psfrag{b}{\LARGE{$a_4$}}
\psfrag{c}{\LARGE{$a_3$}}
\psfrag{d}{\LARGE{$a_2$}}
\psfrag{e}{\LARGE{$a_5$}}
\scalebox{0.50}{\includegraphics{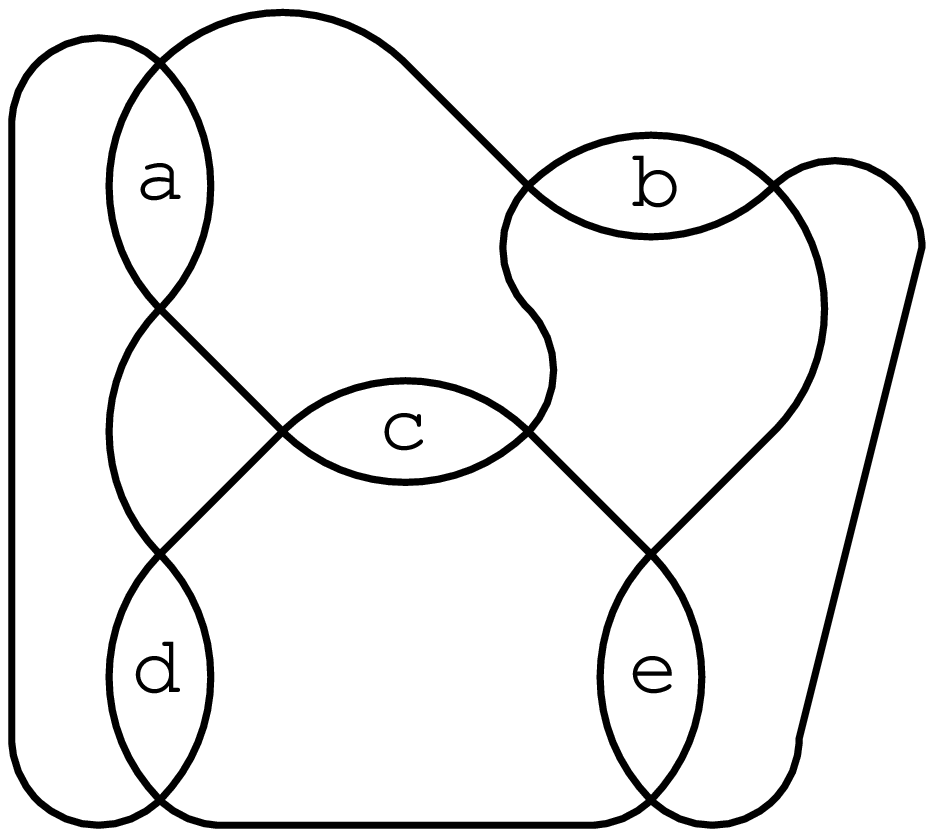}}

\caption{ \ }
$$
(a_1 + a_2, a_1 a_2) M \left( \begin{array}{c}
a_4 a_5 + 1 \\
a_3 a_4 a_5 + a_3 + a_5
\end{array} \right) = (a_1 + a_2, a_1 a_2 + a_2 a_3 + a_3 a_1) M \left( \begin{array}{c}
a_4 a_5 + 1 \\
a_5
\end{array} \right) =
$$
$$
(1, a_1) M \left( \begin{array}{cc}
1 & a_2 \\
a_2 & 0
\end{array} \right) M \left( \begin{array}{cc}
0 & 1 \\
1 & a_3
\end{array} \right) M \left( \begin{array}{cc}
a_4 & 1 \\
1 & 0
\end{array} \right) M \left( \begin{array}{c}
1 \\
a_5
\end{array} \right)
$$
\end{figure}

\begin{figure}

\centering

\psfrag{a}{\LARGE{$a_1$}}
\psfrag{b}{\LARGE{$a_4$}}
\psfrag{c}{\LARGE{$a_3$}}
\psfrag{d}{\LARGE{$a_2$}}
\psfrag{e}{\LARGE{$a_5$}}
\scalebox{0.50}{\includegraphics{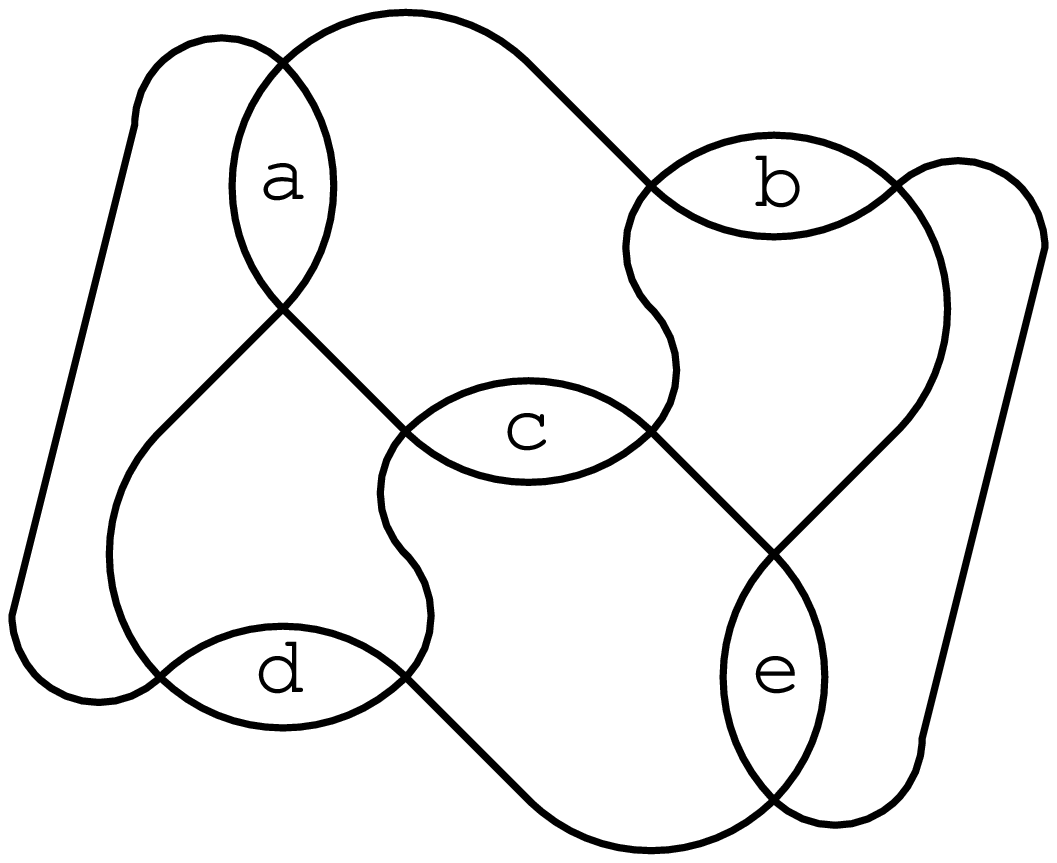}}

\caption{ \ }
$$
(a_1 a_2 + 1, a_1) M \left( \begin{array}{c}
a_4 a_5 + 1 \\
a_3 a_4 a_5 + a_3 + a_5
\end{array} \right) =
$$
$$
(1, a_1) M \left( \begin{array}{cc}
a_2 & 1 \\
1 & 0
\end{array} \right) M \left( \begin{array}{cc}
0 & 1 \\
1 & a_3
\end{array} \right) M \left( \begin{array}{cc}
a_4 & 1 \\
1 & 0
\end{array} \right) M \left( \begin{array}{c}
1 \\
a_5
\end{array} \right)
$$
\end{figure}

\begin{figure}

\centering

\psfrag{a}{\LARGE{$a_1$}}
\psfrag{b}{\LARGE{$a_4$}}
\psfrag{c}{\LARGE{$a_3$}}
\psfrag{d}{\LARGE{$a_2$}}
\psfrag{e}{\LARGE{$a_5$}}
\scalebox{0.50}{\includegraphics{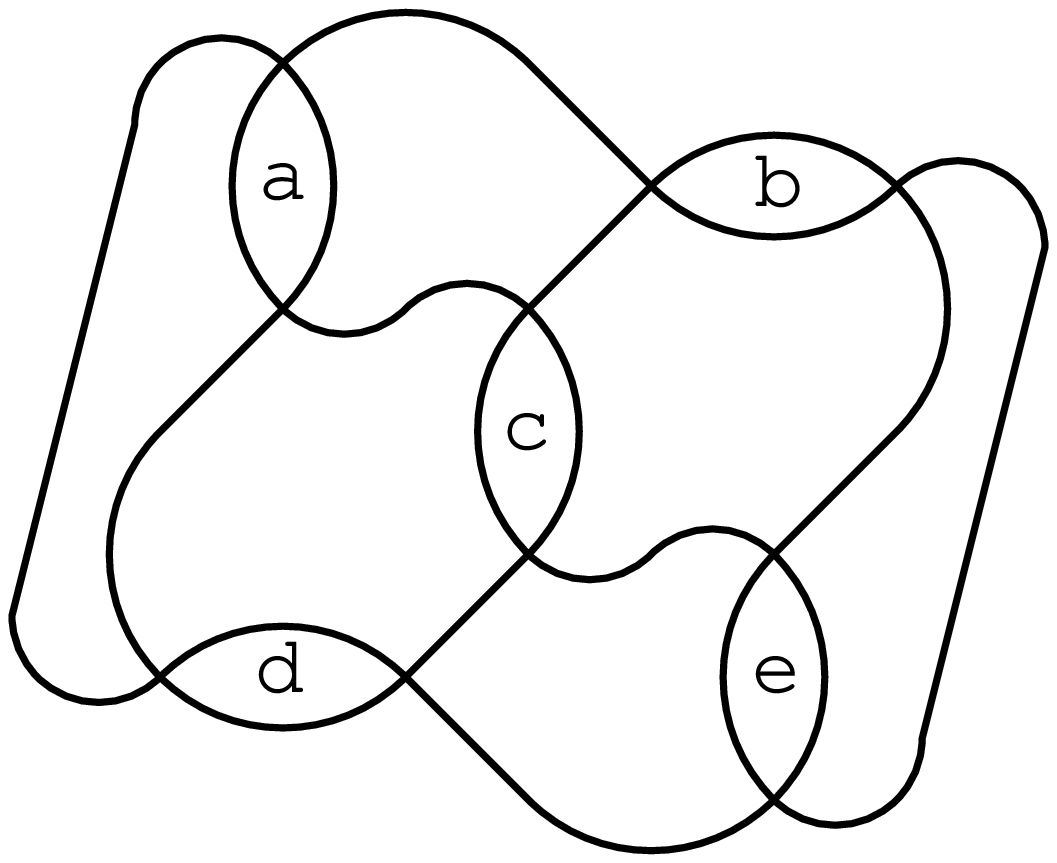}}

\caption{ \ }
$$
(a_1 a_2 + 1, a_1) M \left( \begin{array}{c}
a_3 (a_4 a_5 + 1) \\
a_5 (a_3 + a_4) + 1
\end{array} \right) =
$$
$$
(1, a_1) M \left( \begin{array}{cc}
a_2 & 1 \\
1 & 0
\end{array} \right) M \left( \begin{array}{cc}
0 & a_3 \\
a_3 & 1
\end{array} \right) M \left( \begin{array}{cc}
a_4 & 1 \\
1 & 0
\end{array} \right) M \left( \begin{array}{c}
1 \\
a_5
\end{array} \right)
$$
\end{figure}


\begin{figure}

\section{THE FAMILIES OF KNOTS OF SIX CONWAYS (44 CASES)}

\subsection{The families of knots with seed the torus link $6_2^1$ (two cases)}

\centering

\psfrag{a}{\LARGE{$a_1$}}
\psfrag{b}{\LARGE{$a_2$}}
\psfrag{c}{\LARGE{$a_3$}}
\psfrag{d}{\LARGE{$a_4$}}
\psfrag{e}{\LARGE{$a_5$}}
\psfrag{f}{\LARGE{$a_6$}}
\scalebox{0.50}{\includegraphics{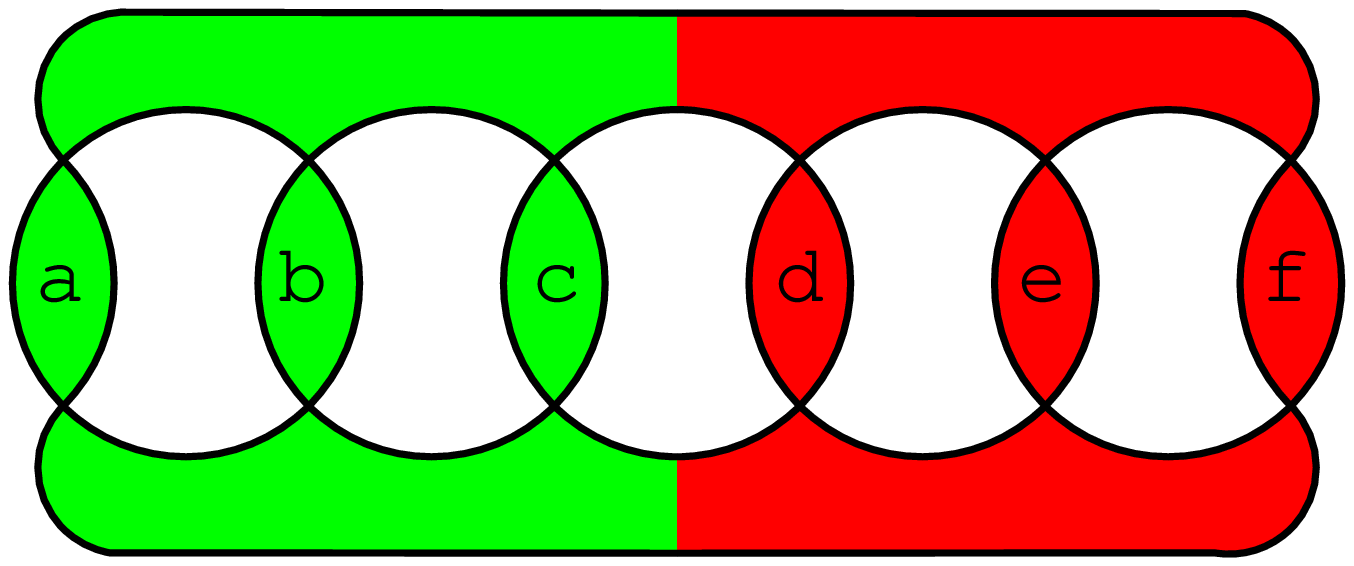}}

\caption{ \ }
$$
(a_1 a_2 a_3, a_1 a_2 + a_2 a_3 + a_3 a_1) \left( \begin{array}{cc}
0 & 1 \\
1 & 0
\end{array} \right) \left( \begin{array}{cc}
a_4 a_5 a_6 \\
a_4 a_5 + a_5 a_6 + a_6 a_4
\end{array} \right) =
$$
$$
(a_1, 1) M \left( \begin{array}{cc}
0 & a_2 \\
a_2 & 1
\end{array} \right) M \left( \begin{array}{cc}
0 & a_3 \\
a_3 & 1
\end{array} \right) M \left( \begin{array}{cc}
0 & a_4 \\
a_4 & 1
\end{array} \right) M \left( \begin{array}{cc}
0 & a_5 \\
a_5 & 1
\end{array} \right) M \left( \begin{array}{c}
a_6 \\
1
\end{array} \right)
$$
\end{figure}

\begin{figure}

\centering

\psfrag{a}{\LARGE{$a_1$}}
\psfrag{b}{\LARGE{$a_2$}}
\psfrag{c}{\LARGE{$a_3$}}
\psfrag{d}{\LARGE{$a_4$}}
\psfrag{e}{\LARGE{$a_5$}}
\psfrag{f}{\LARGE{$a_6$}}
\scalebox{0.5}{\includegraphics{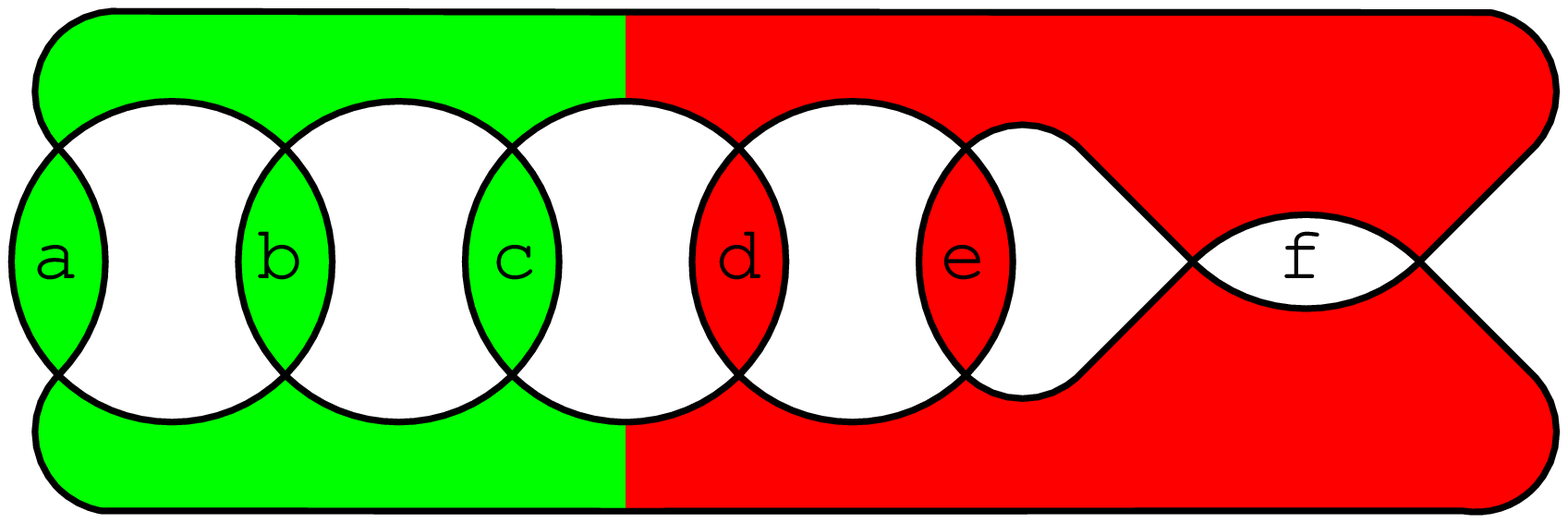}}

\caption{ \ }
$$
(a_1 a_2 a_3, a_1 a_2 + a_2 a_3 + a_3 a_1)  \left( \begin{array}{cc}
0 & 1 \\
1 & 0
\end{array} \right) \left( \begin{array}{cc}
a_4 a_5 \\
a_4 a_5 a_6 + a_4 + a_5
\end{array} \right) =
$$
$$
(a_1, 1) M \left( \begin{array}{cc}
0 & a_2 \\
a_2 & 1
\end{array} \right) M \left( \begin{array}{cc}
0 & a_3 \\
a_3 & 1
\end{array} \right) M \left( \begin{array}{cc}
0 & a_4 \\
a_4 & 1
\end{array} \right) M \left( \begin{array}{cc}
0 & a_5 \\
a_5 & 1
\end{array} \right) M \left( \begin{array}{c}
1 \\
a_6
\end{array} \right)
$$
\end{figure}

\begin{figure}

\subsection{Families of knots with seed the knot $6_1$ (four cases)}

\centering

\psfrag{a}{\LARGE{$a_1$}}
\psfrag{b}{\LARGE{$a_2$}}
\psfrag{c}{\LARGE{$a_3$}}
\psfrag{d}{\LARGE{$a_4$}}
\psfrag{e}{\LARGE{$a_5$}}
\psfrag{f}{\LARGE{$a_6$}}

\scalebox{0.5}{\includegraphics{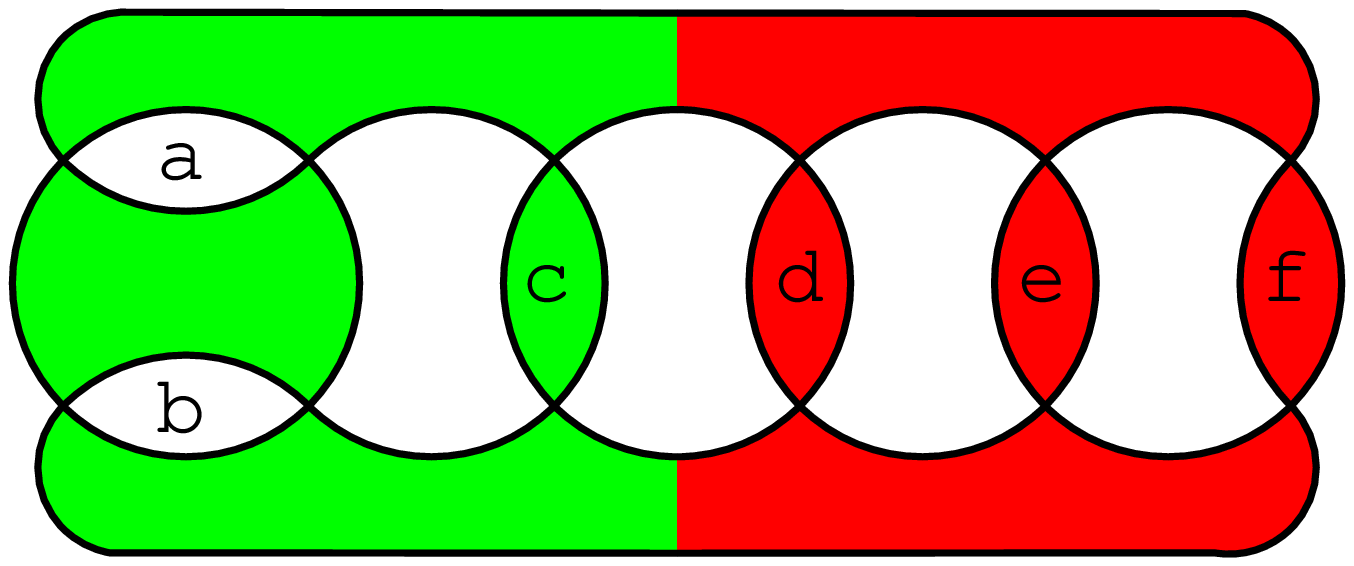}}

\caption{ \ }
$$
((a_1 + a_2) a_3, a_1 a_2 a_3 + a_1 + a_2) \left( \begin{array}{cc}
0 & 1 \\
1 & 0
\end{array} \right)  \left( \begin{array}{cc}
a_4 a_5 a_6 \\
a_4 a_5 + a_5 a_6 + a_6 a_4
\end{array} \right) =
$$
$$
(1, a_1) M \left( \begin{array}{cc}
1 & a_2 \\
a_2 & 0
\end{array} \right) M \left( \begin{array}{cc}
0 & a_3 \\
a_3 & 1
\end{array} \right) M \left( \begin{array}{cc}
0 & a_4 \\
a_4 & 1
\end{array} \right) M \left( \begin{array}{cc}
0 & a_5 \\
a_5 & 1
\end{array} \right) M \left( \begin{array}{c}
a_6 \\
1
\end{array} \right)
$$
\end{figure}

\begin{figure}

\centering

\psfrag{a}{\LARGE{$a_1$}}
\psfrag{b}{\LARGE{$a_2$}}
\psfrag{c}{\LARGE{$a_3$}}
\psfrag{d}{\LARGE{$a_4$}}
\psfrag{e}{\LARGE{$a_5$}}
\psfrag{f}{\LARGE{$a_6$}}

\scalebox{0.5}{\includegraphics{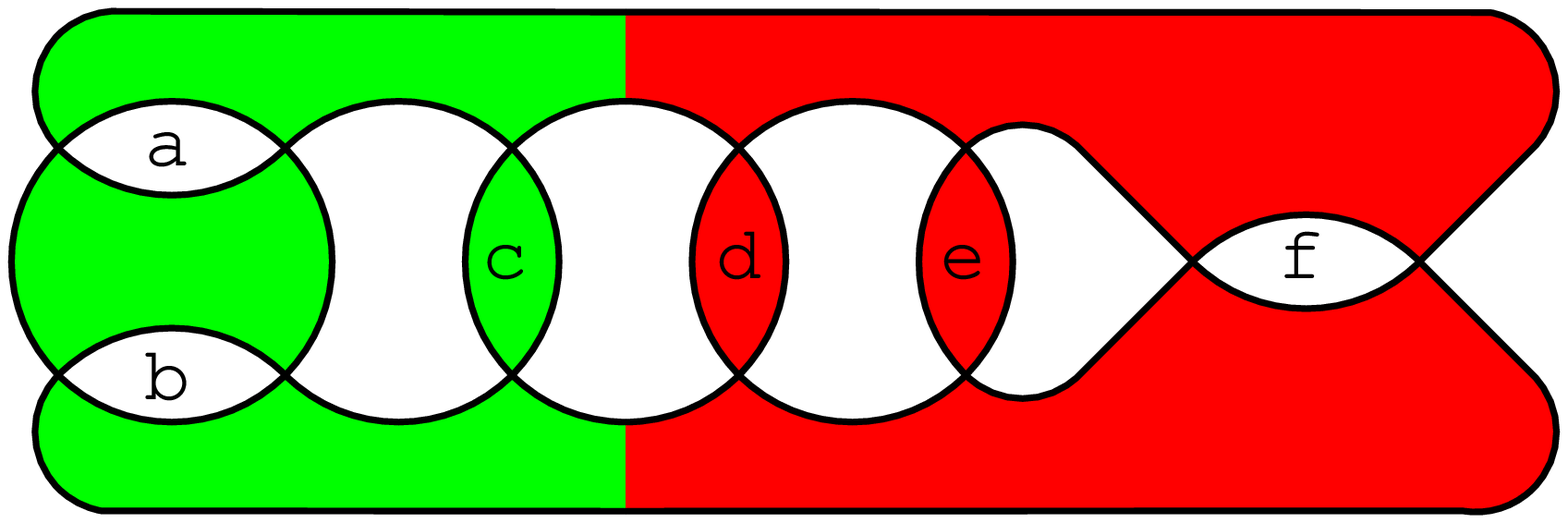}}

\caption{ \ }
$$
(a_1 + a_2) a_3, (a_1 a_2 a_3 + a_1 + a_2) \left( \begin{array}{cc}
0 & 1 \\
1 & 0
\end{array} \right)  \left( \begin{array}{cc}
a_4 a_5 \\
a_4 a_5 a_6 + a_4 + a_5
\end{array} \right) =
$$
$$
(1, a_1) M \left( \begin{array}{cc}
1 & a_2 \\
a_2 & 0
\end{array} \right) M \left( \begin{array}{cc}
0 & a_3 \\
a_3 & 1
\end{array} \right) M \left( \begin{array}{cc}
0 & a_4 \\
a_4 & 1
\end{array} \right) M \left( \begin{array}{cc}
0 & a_5 \\
a_5 & 1
\end{array} \right) M \left( \begin{array}{c}
1 \\
a_6
\end{array} \right)
$$
\end{figure}

\begin{figure}

\centering

\psfrag{a}{\LARGE{$a_1$}}
\psfrag{b}{\LARGE{$a_2$}}
\psfrag{c}{\LARGE{$a_3$}}
\psfrag{d}{\LARGE{$a_4$}}
\psfrag{e}{\LARGE{$a_5$}}
\psfrag{f}{\LARGE{$a_6$}}

\scalebox{0.5}{\includegraphics{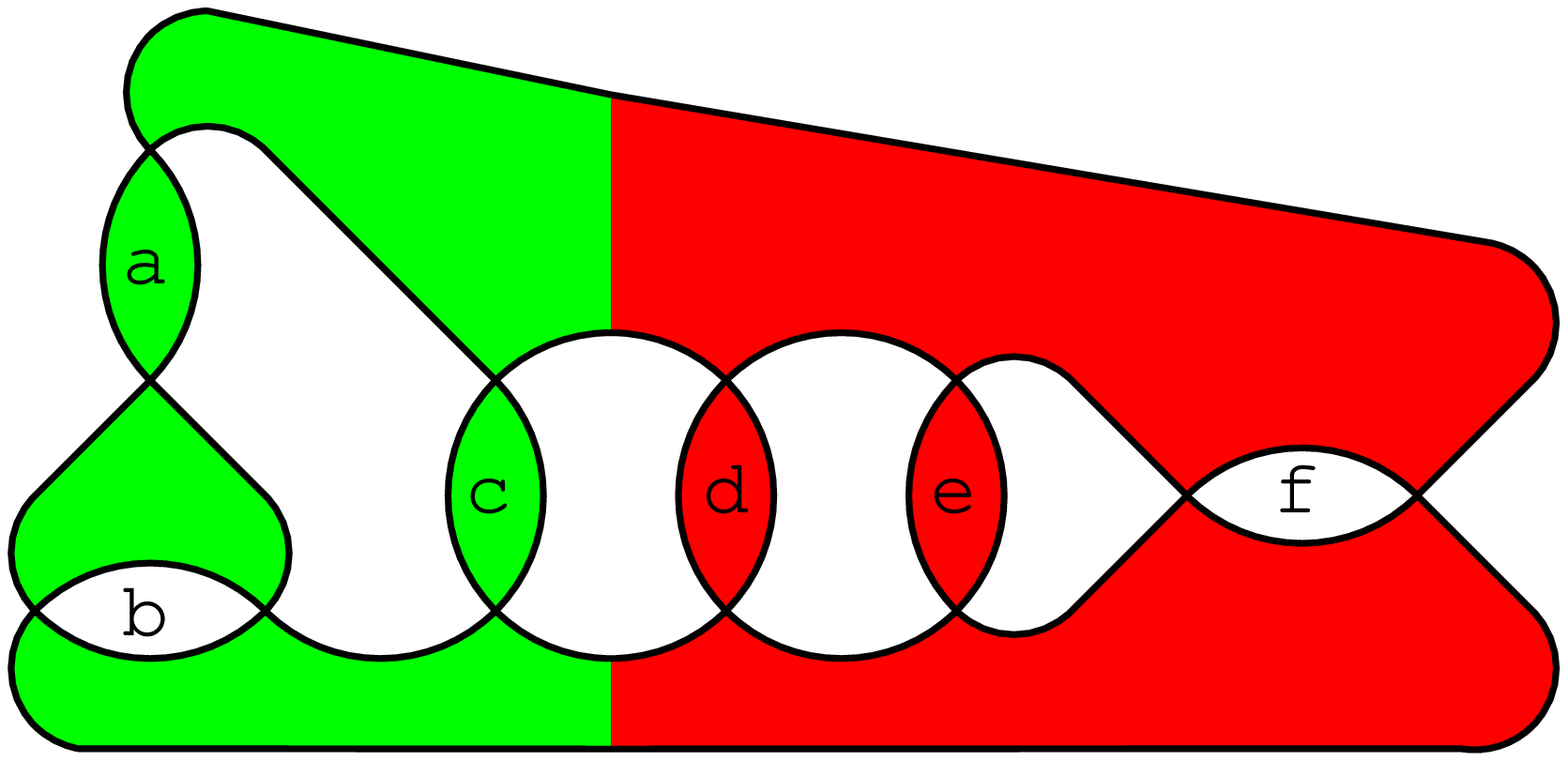}}
\caption{ \ }

$$
((a_1 a_2 + 1) a_3, a_2 (a_1 + a_3) + 1) \left( \begin{array}{cc}
0 & 1 \\
1 & 0
\end{array} \right)   \left( \begin{array}{cc}
a_4 a_5 \\
a_4 a_5 a_6 + a_4 + a_5
\end{array} \right) =
$$
$$
(a_1, 1) M \left( \begin{array}{cc}
1 & a_2 \\
a_2 & 0
\end{array} \right) M \left( \begin{array}{cc}
0 & a_3 \\
a_3 & 1
\end{array} \right) M \left( \begin{array}{cc}
0 & a_4 \\
a_4 & 1
\end{array} \right) M \left( \begin{array}{cc}
0 & a_5 \\
a_5 & 1
\end{array} \right) M \left( \begin{array}{c}
1 \\
a_6
\end{array} \right)
$$
\end{figure}

\begin{figure}

\centering

\psfrag{a}{\LARGE{$a_1$}}
\psfrag{b}{\LARGE{$a_2$}}
\psfrag{c}{\LARGE{$a_3$}}
\psfrag{d}{\LARGE{$a_4$}}
\psfrag{e}{\LARGE{$a_5$}}
\psfrag{f}{\LARGE{$a_6$}}

\scalebox{0.5}{\includegraphics{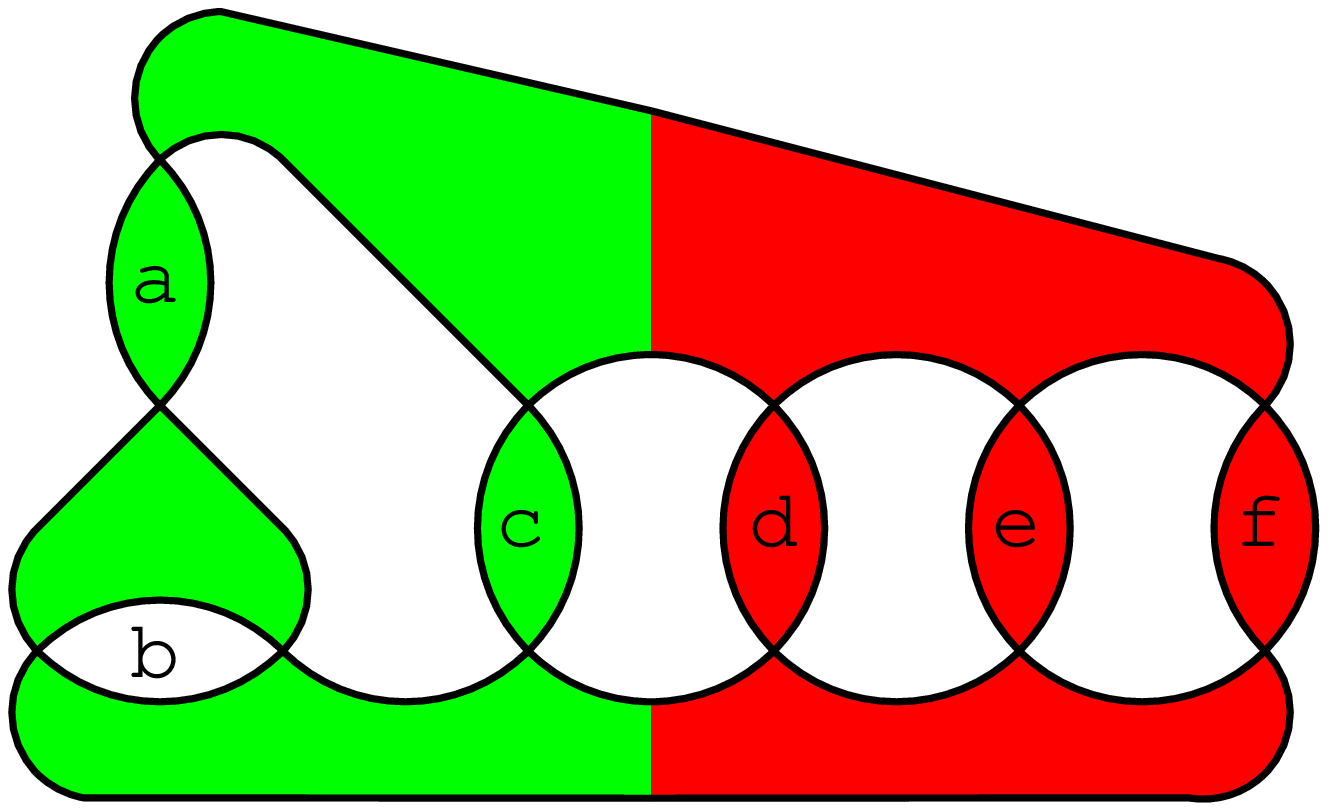}}
\caption{ \ }
$$
((a_1 a_2 + 1) a_3, a_2 (a_1 + a_3) + 1) \left( \begin{array}{cc}
0 & 1 \\
1 & 0
\end{array} \right)   \left( \begin{array}{cc}
a_4 a_5 a_6 \\
a_4 a_5 + a_5 a_6 + a_6 a_4
\end{array} \right) =
$$
$$
(a_1, 1) M \left( \begin{array}{cc}
1 & a_2 \\
a_2 & 0
\end{array} \right) M \left( \begin{array}{cc}
0 & a_3 \\
a_3 & 1
\end{array} \right) M \left( \begin{array}{cc}
0 & a_4 \\
a_4 & 1
\end{array} \right) M \left( \begin{array}{cc}
0 & a_5 \\
a_5 & 1
\end{array} \right) M \left( \begin{array}{c}
a_6 \\
1
\end{array} \right)
$$
\end{figure}

\begin{figure}

\subsection{Families of knots with seed the link $6_2^2$ (three cases)}

\centering

\psfrag{a}{\LARGE{$a_1$}}
\psfrag{b}{\LARGE{$a_2$}}
\psfrag{c}{\LARGE{$a_3$}}
\psfrag{d}{\LARGE{$a_4$}}
\psfrag{e}{\LARGE{$a_5$}}
\psfrag{f}{\LARGE{$a_6$}}

\scalebox{0.5}{\includegraphics{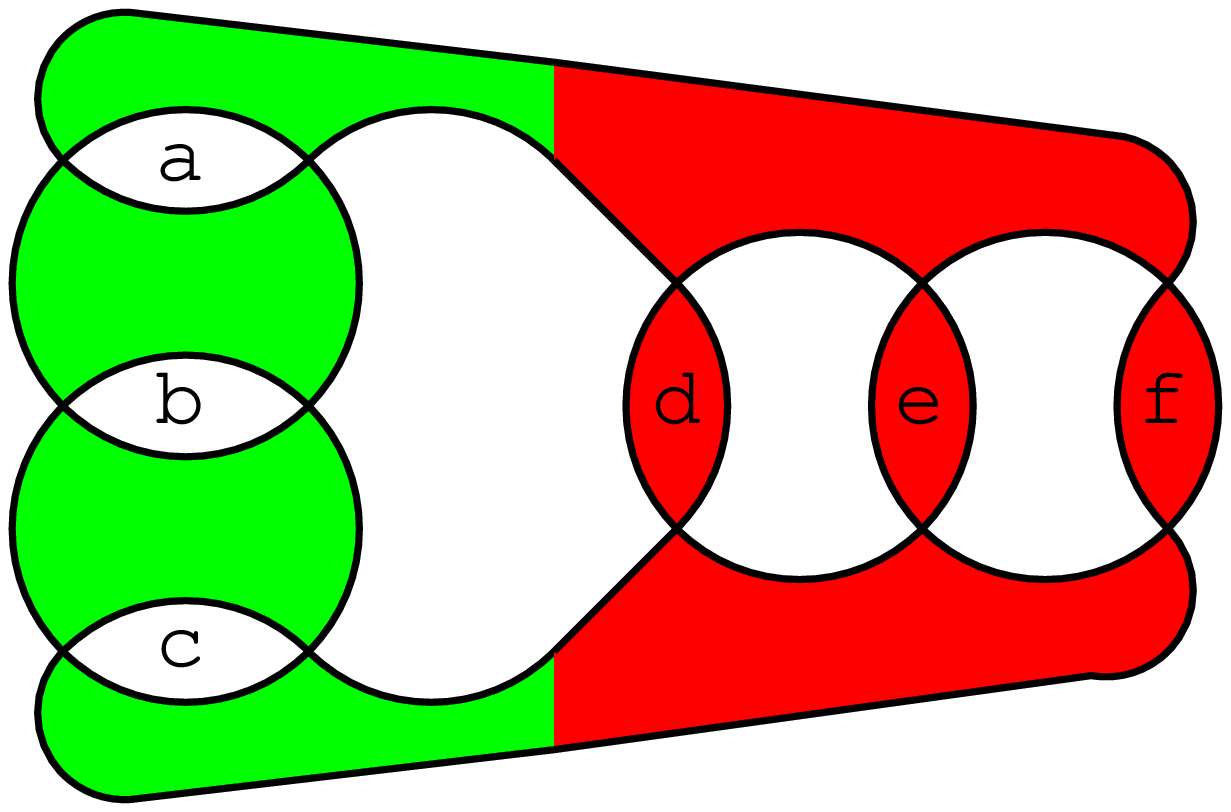}}
\caption{ \ }
$$
(a_1 a_2 + a_2 a_3 + a_3 a_1, a_1 a_2 a_3)  \left( \begin{array}{cc}
0 & 1 \\
1 & 0
\end{array} \right) \left( \begin{array}{cc}
a_4 a_5 a_6 \\
a_4 a_5 + a_5 a_6 + a_6 a_4
\end{array} \right) =
$$
$$
(1, a_1) M \left( \begin{array}{cc}
1 & a_2 \\
a_2 & 0
\end{array} \right) M \left( \begin{array}{cc}
1 & a_3 \\
a_3 & 0
\end{array} \right) M \left( \begin{array}{cc}
0 & a_4 \\
a_4 & 1
\end{array} \right) M \left( \begin{array}{cc}
0 & a_5 \\
a_5 & 1
\end{array} \right) M \left( \begin{array}{c}
a_6 \\
1
\end{array} \right)
$$
\end{figure}

\begin{figure}

\centering

\psfrag{a}{\LARGE{$a_1$}}
\psfrag{b}{\LARGE{$a_2$}}
\psfrag{c}{\LARGE{$a_3$}}
\psfrag{d}{\LARGE{$a_4$}}
\psfrag{e}{\LARGE{$a_5$}}
\psfrag{f}{\LARGE{$a_6$}}

\scalebox{0.5}{\includegraphics{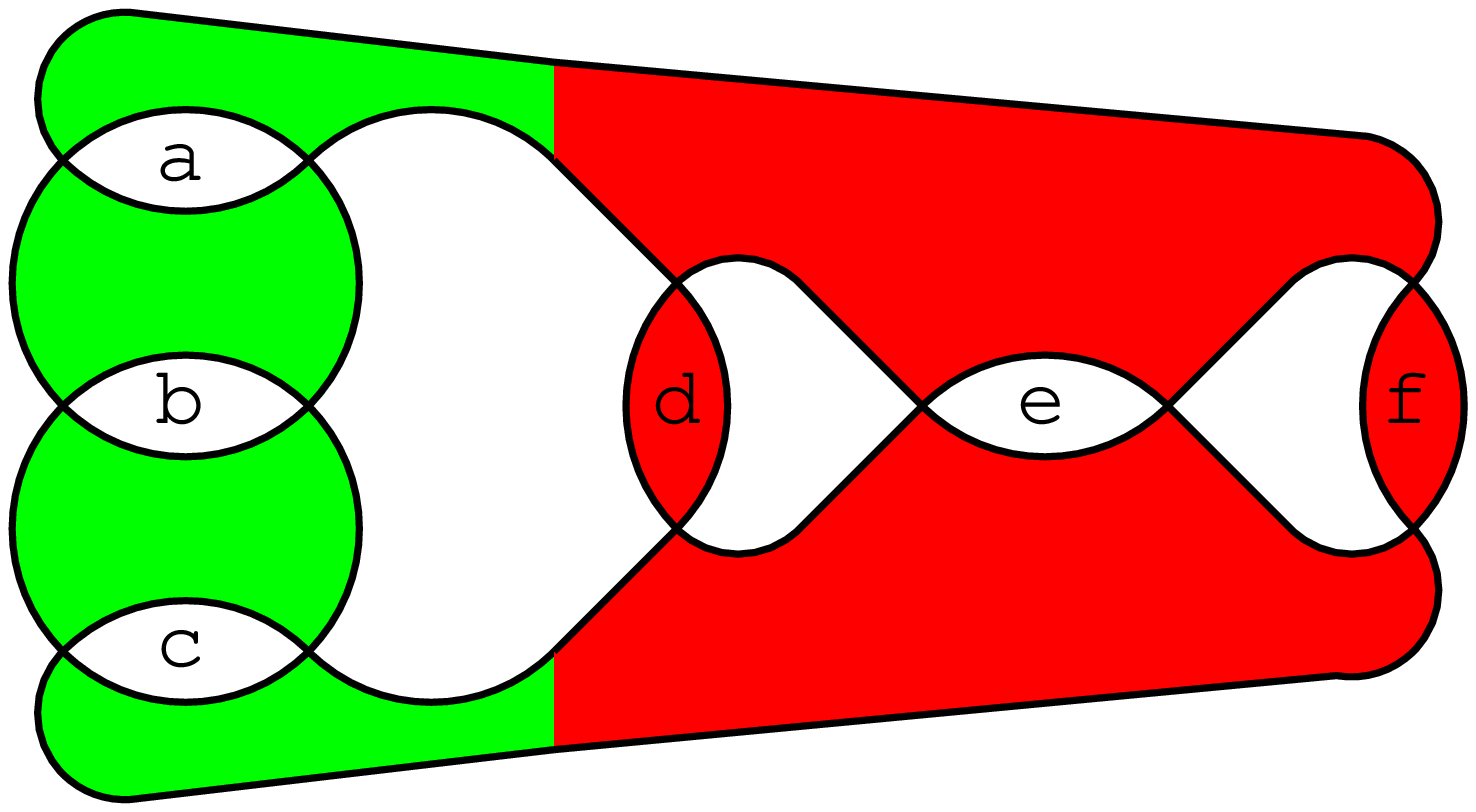}}
\caption{ \ }
$$
(a_1 a_2 + a_2 a_3 + a_3 a_1, a_1 a_2 a_3)  \left( \begin{array}{cc}
0 & 1 \\
1 & 0
\end{array} \right)   \left( \begin{array}{cc}
a_4 a_6 \\
a_4 a_5 a_6 + a_4 + a_6
\end{array} \right) =
$$
$$
(1, a_1) M \left( \begin{array}{cc}
1 & a_2 \\
a_2 & 0
\end{array} \right) M \left( \begin{array}{cc}
1 & a_3 \\
a_3 & 0
\end{array} \right) M \left( \begin{array}{cc}
0 & a_4 \\
a_4 & 1
\end{array} \right) M \left( \begin{array}{cc}
0 & 1 \\
1 & a_5
\end{array} \right) M \left( \begin{array}{c}
a_6 \\
1
\end{array} \right)
$$
\end{figure}

\begin{figure}

\centering

\psfrag{a}{\LARGE{$a_1$}}
\psfrag{b}{\LARGE{$a_2$}}
\psfrag{c}{\LARGE{$a_3$}}
\psfrag{d}{\LARGE{$a_4$}}
\psfrag{e}{\LARGE{$a_5$}}
\psfrag{f}{\LARGE{$a_6$}}

\scalebox{0.5}{\includegraphics{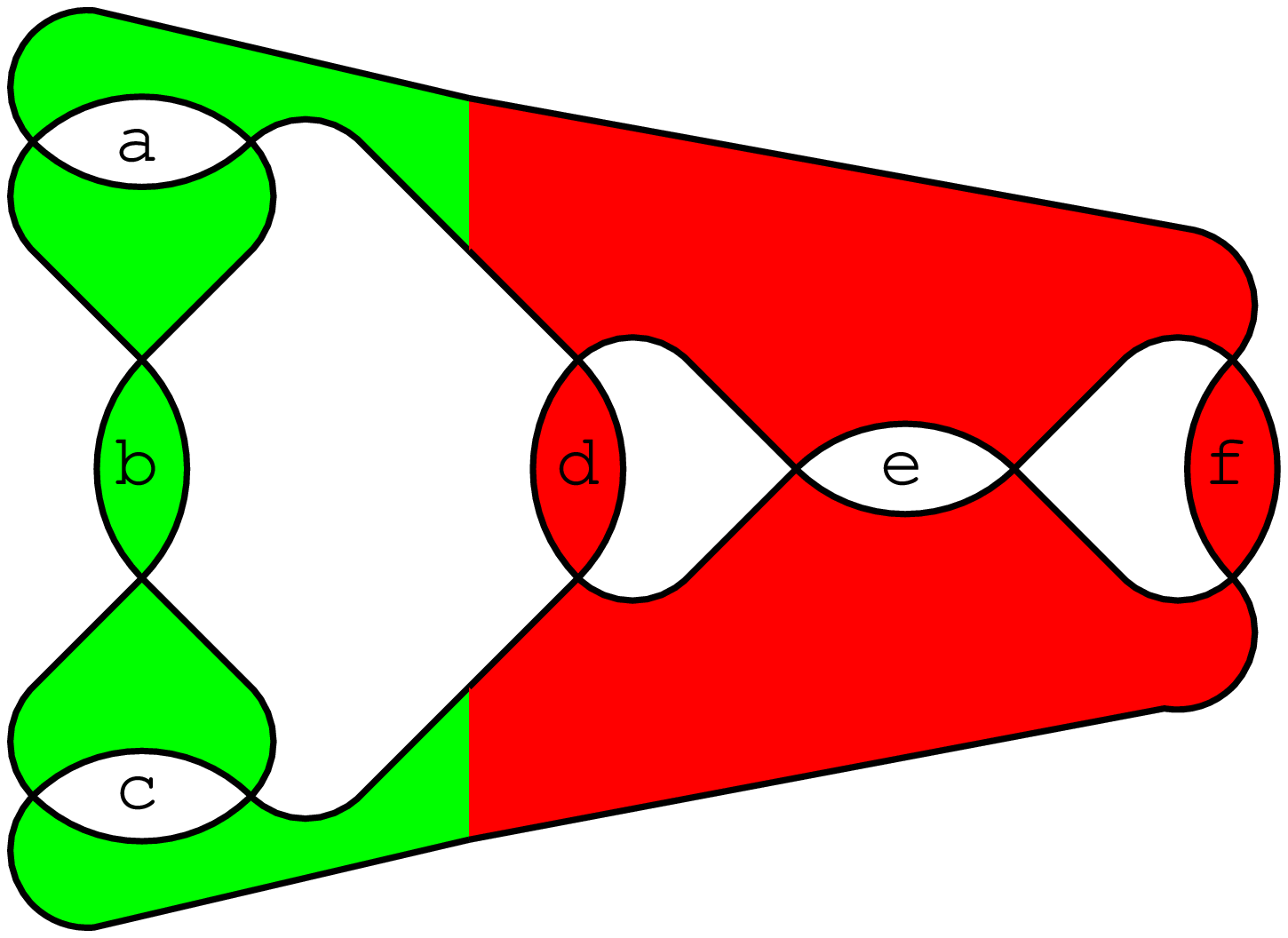}}

\caption{ \ }
$$
(a_1 a_2 a_3 + a_1 + a_3, a_1 a_3)  \left( \begin{array}{cc}
0 & 1 \\
1 & 0
\end{array} \right)   \left( \begin{array}{cc}
a_4 a_6 \\
a_4 a_5 a_6 + a_4 + a_6
\end{array} \right) =
$$
$$
(1, a_1) M \left( \begin{array}{cc}
a_2 & 1 \\
1 & 0
\end{array} \right) M \left( \begin{array}{cc}
1 & a_3 \\
a_3 & 0
\end{array} \right) M \left( \begin{array}{cc}
0 & a_4 \\
a_4 & 1
\end{array} \right) M \left( \begin{array}{cc}
0 & 1 \\
1 & a_5
\end{array} \right) M \left( \begin{array}{c}
a_6 \\
1
\end{array} \right)
$$
\end{figure}

\begin{figure}
\setcounter{subsection}{3}
\subsection{Families of knots with seed the link $6_3^2$ with twelve terms (six cases)}

\centering

\psfrag{a}{\LARGE{$a_1$}}
\psfrag{b}{\LARGE{$a_2$}}
\psfrag{c}{\LARGE{$a_3$}}
\psfrag{d}{\LARGE{$a_4$}}
\psfrag{e}{\LARGE{$a_5$}}
\psfrag{f}{\LARGE{$a_6$}}

\scalebox{0.5}{\includegraphics{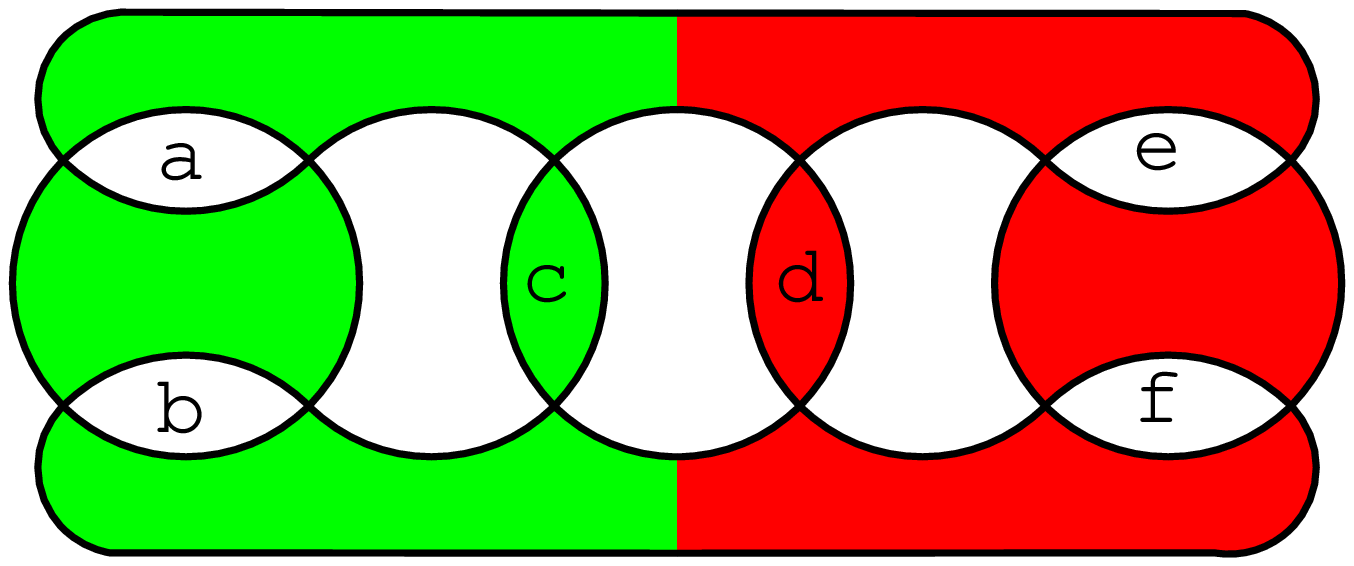}}

\caption{ \ }
$$
((a_1 + a_2) a_3, a_1 a_2 a_3 + a_1 + a_2) \left( \begin{array}{cc}
0 & 1 \\
1 & 0
\end{array} \right)  \left( \begin{array}{cc}
a_4 (a_5 + a_6) \\
a_4 a_5 a_6 + a_5 + a_6
\end{array} \right) =
$$
$$
(1, a_1) M \left( \begin{array}{cc}
1 & a_2 \\
a_2 & 0
\end{array} \right) M \left( \begin{array}{cc}
0 & a_3 \\
a_3 & 1
\end{array} \right) M \left( \begin{array}{cc}
0 & a_4 \\
a_4 & 1
\end{array} \right) M \left( \begin{array}{cc}
1 & a_5 \\
a_5 & 0
\end{array} \right) M \left( \begin{array}{c}
1 \\
a_6
\end{array} \right)
$$
\end{figure}

\begin{figure}

\centering

\psfrag{a}{\LARGE{$a_1$}}
\psfrag{b}{\LARGE{$a_2$}}
\psfrag{c}{\LARGE{$a_3$}}
\psfrag{d}{\LARGE{$a_4$}}
\psfrag{e}{\LARGE{$a_5$}}
\psfrag{f}{\LARGE{$a_6$}}

\scalebox{0.5}{\includegraphics{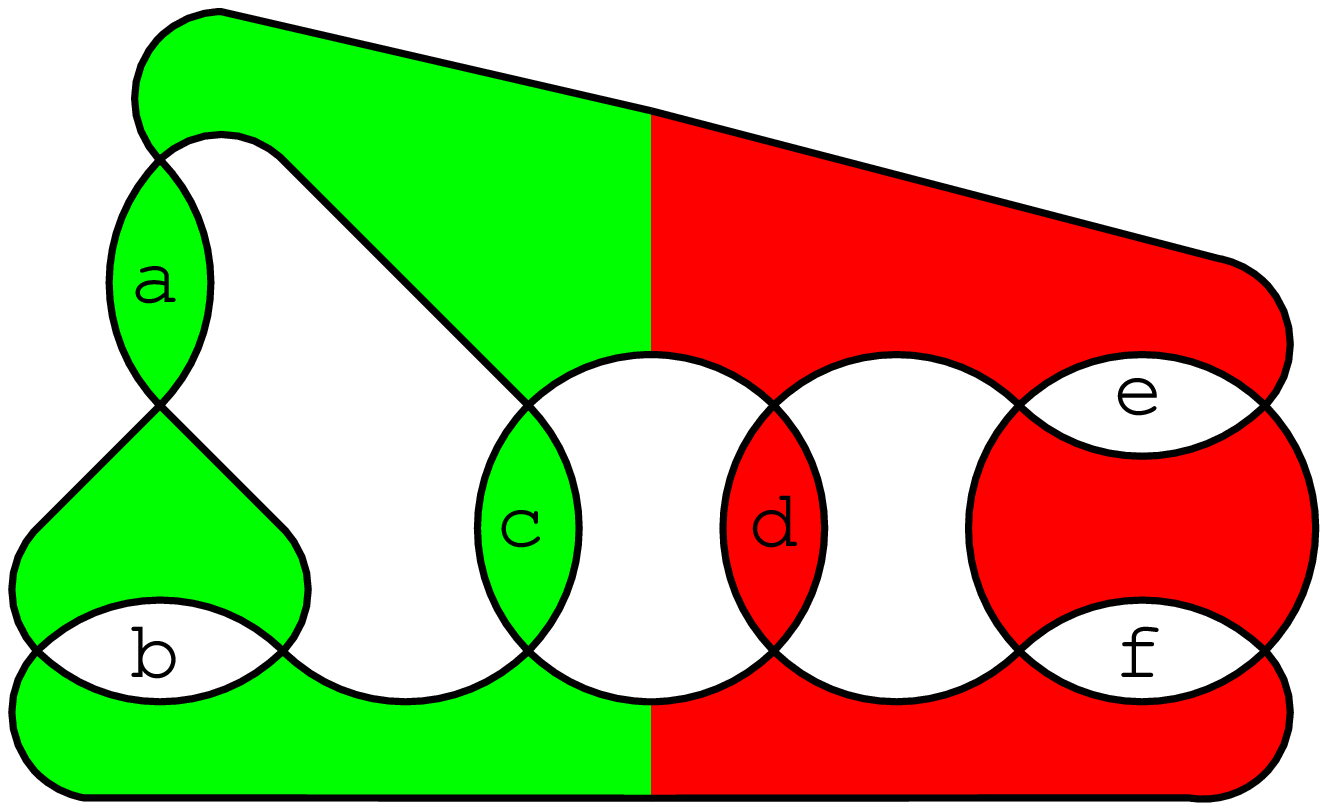}}
\caption{ \ }
$$
((a_1 a_2 + 1) a_3, a_2 (a_1 + a_3) + 1) \left( \begin{array}{cc}
0 & 1 \\
1 & 0
\end{array} \right)  \left( \begin{array}{cc}
a_4 (a_5 + a_6) \\
a_4 a_5 a_6 + a_5 + a_6
\end{array} \right) =
$$
$$
(a_1, 1) M \left( \begin{array}{cc}
1 & a_2 \\
a_2 & 0
\end{array} \right) M \left( \begin{array}{cc}
0 & a_3 \\
a_3 & 1
\end{array} \right) M \left( \begin{array}{cc}
0 & a_4 \\
a_4 & 1
\end{array} \right) M \left( \begin{array}{cc}
1 & a_5 \\
a_5 & 0
\end{array} \right) M \left( \begin{array}{c}
1 \\
a_6
\end{array} \right)
$$
\end{figure}

\begin{figure}

\centering

\psfrag{a}{\LARGE{$a_1$}}
\psfrag{b}{\LARGE{$a_2$}}
\psfrag{c}{\LARGE{$a_3$}}
\psfrag{d}{\LARGE{$a_4$}}
\psfrag{e}{\LARGE{$a_5$}}
\psfrag{f}{\LARGE{$a_6$}}

\scalebox{0.5}{\includegraphics{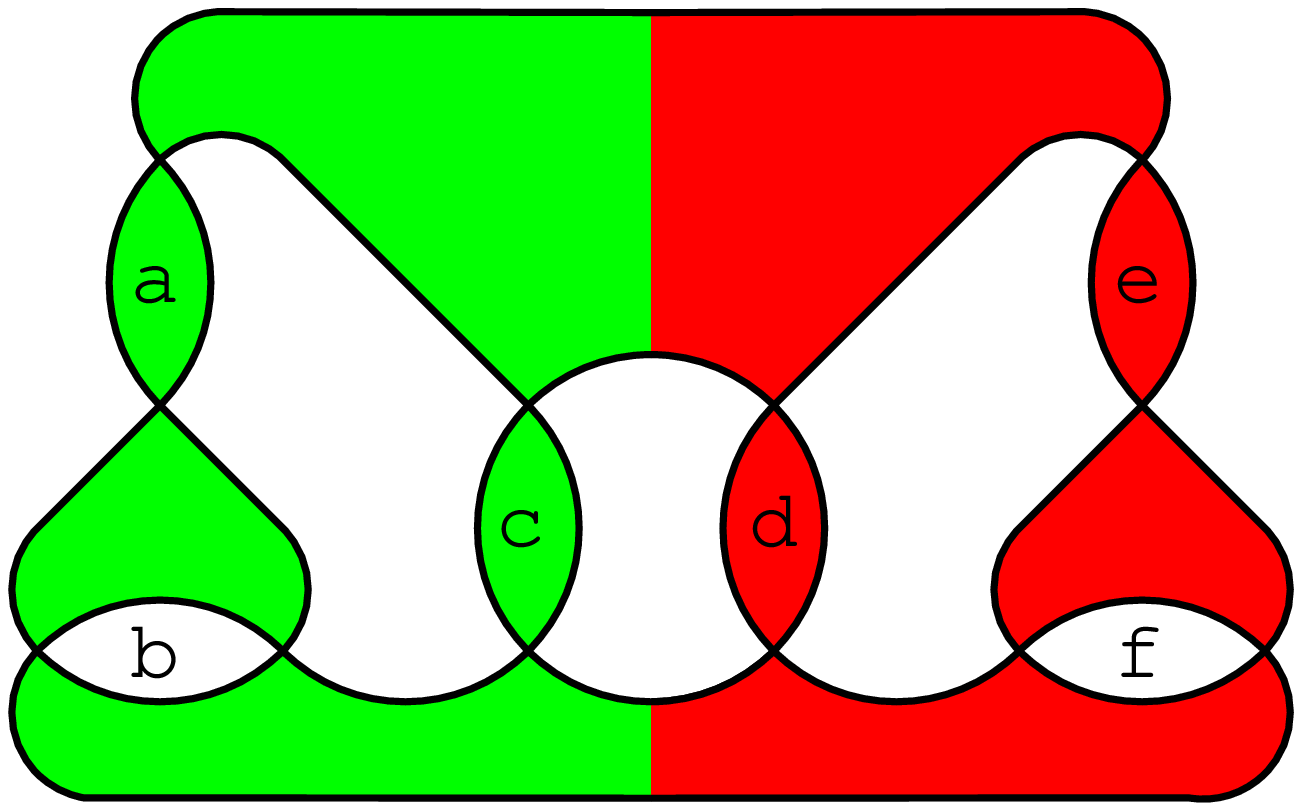}}

\caption{ \ }
$$
((a_1 a_2 + 1) a_3, a_2 (a_1 + a_3) + 1) \left( \begin{array}{cc}
0 & 1 \\
1 & 0
\end{array} \right)  \left( \begin{array}{cc}
a_4 (a_5 a_6 + 1) \\
(a_4 + a_5) a_6 + 1
\end{array} \right) =
$$
$$
(a_1, 1) M \left( \begin{array}{cc}
1 & a_2 \\
a_2 & 0
\end{array} \right) M \left( \begin{array}{cc}
0 & a_3 \\
a_3 & 1
\end{array} \right) M \left( \begin{array}{cc}
0 & a_4 \\
a_4 & 1
\end{array} \right) M \left( \begin{array}{cc}
a_5 & 1 \\
1 & 0
\end{array} \right) M \left( \begin{array}{c}
1 \\
a_6
\end{array} \right)
$$
\end{figure}

\begin{figure}

\centering

\psfrag{a}{\LARGE{$a_1$}}
\psfrag{b}{\LARGE{$a_2$}}
\psfrag{c}{\LARGE{$a_3$}}
\psfrag{d}{\LARGE{$a_4$}}
\psfrag{e}{\LARGE{$a_5$}}
\psfrag{f}{\LARGE{$a_6$}}

\scalebox{0.5}{\includegraphics{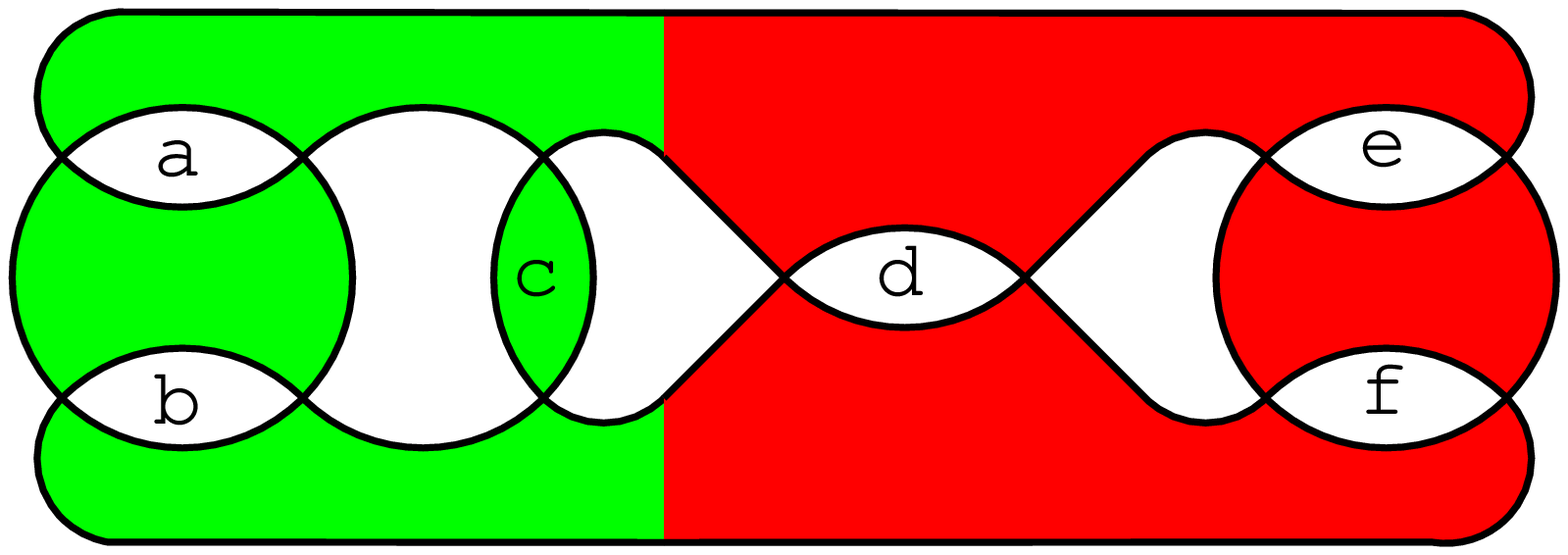}}

\caption{ \ }
$$
((a_1 + a_2) a_3, a_1 a_2 a_3 + a_1 + a_2) \left( \begin{array}{cc}
0 & 1 \\
1 & 0
\end{array} \right)  \left( \begin{array}{cc}
a_5 + a_6 \\
a_4 a_5 + a_5 a_6 + a_6 a_4
\end{array} \right) =
$$
$$
(1, a_1) M \left( \begin{array}{cc}
1 & a_2 \\
a_2 & 0
\end{array} \right) M \left( \begin{array}{cc}
0 & a_3 \\
a_3 & 1
\end{array} \right) M \left( \begin{array}{cc}
0 & 1 \\
1 & a_4
\end{array} \right) M \left( \begin{array}{cc}
1 & a_5 \\
a_5 & 0
\end{array} \right) M \left( \begin{array}{c}
1 \\
a_6
\end{array} \right)
$$
\end{figure}

\begin{figure}

\centering

\psfrag{a}{\LARGE{$a_1$}}
\psfrag{b}{\LARGE{$a_2$}}
\psfrag{c}{\LARGE{$a_3$}}
\psfrag{d}{\LARGE{$a_4$}}
\psfrag{e}{\LARGE{$a_5$}}
\psfrag{f}{\LARGE{$a_6$}}

\scalebox{0.5}{\includegraphics{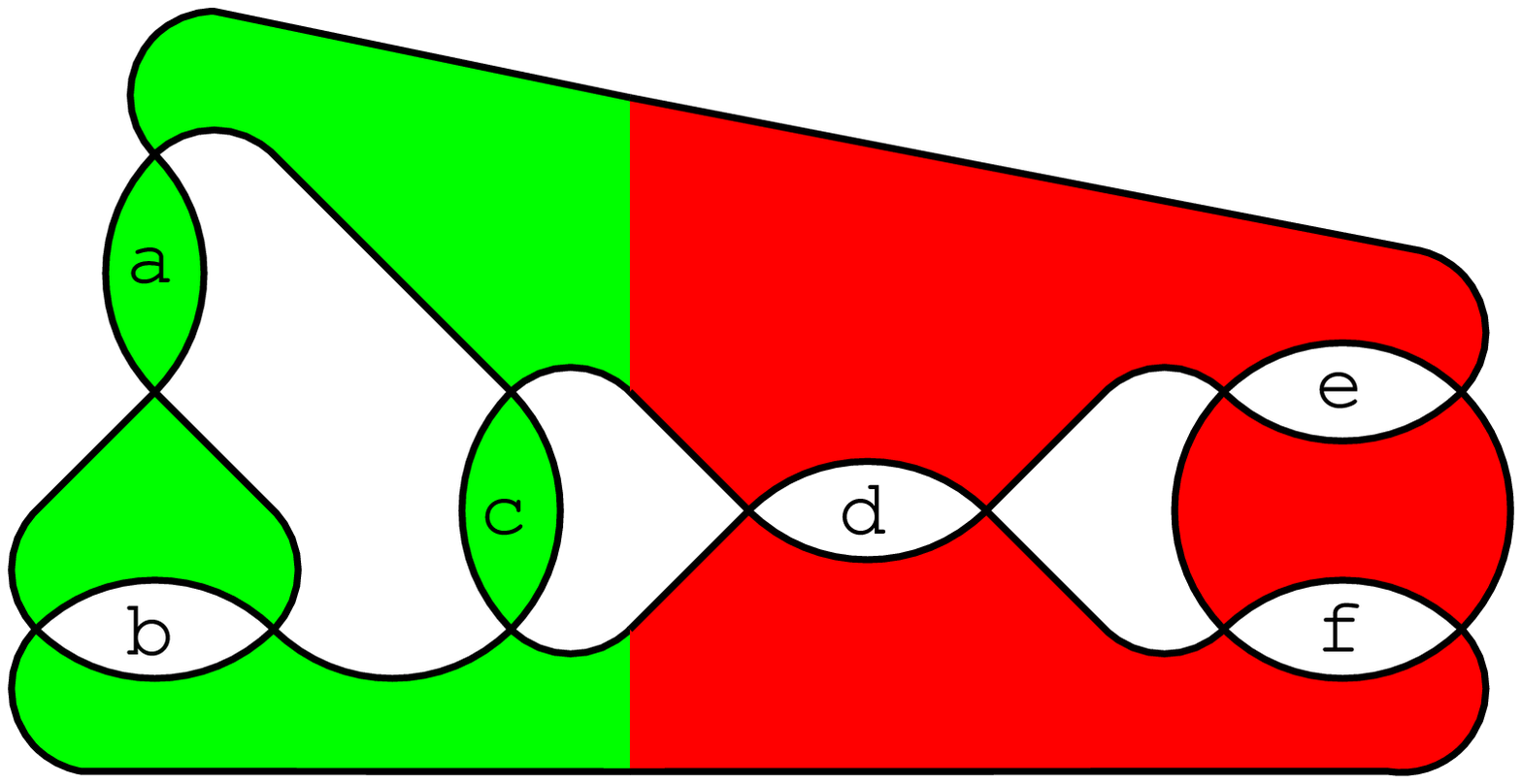}}

\caption{ \ }
$$
((a_1 a_2 + 1) a_3, a_2 (a_1 + a_3) + 1) \left( \begin{array}{cc}
0 & 1 \\
1 & 0
\end{array} \right)  \left( \begin{array}{cc}
a_5 + a_6 \\
a_4 a_5 + a_5 a_6 + a_6 a_4
\end{array} \right) =
$$
$$
(a_1, 1) M \left( \begin{array}{cc}
1 & a_2 \\
a_2 & 0
\end{array} \right) M \left( \begin{array}{cc}
0 & a_3 \\
a_3 & 1
\end{array} \right) M \left( \begin{array}{cc}
0 & 1 \\
1 & a_4
\end{array} \right) M \left( \begin{array}{cc}
1 & a_5 \\
a_5 & 0
\end{array} \right) M \left( \begin{array}{c}
1 \\
a_6
\end{array} \right)
$$
\end{figure}

\begin{figure}

\centering

\psfrag{a}{\LARGE{$a_1$}}
\psfrag{b}{\LARGE{$a_2$}}
\psfrag{c}{\LARGE{$a_3$}}
\psfrag{d}{\LARGE{$a_4$}}
\psfrag{e}{\LARGE{$a_5$}}
\psfrag{f}{\LARGE{$a_6$}}

\scalebox{0.5}{\includegraphics{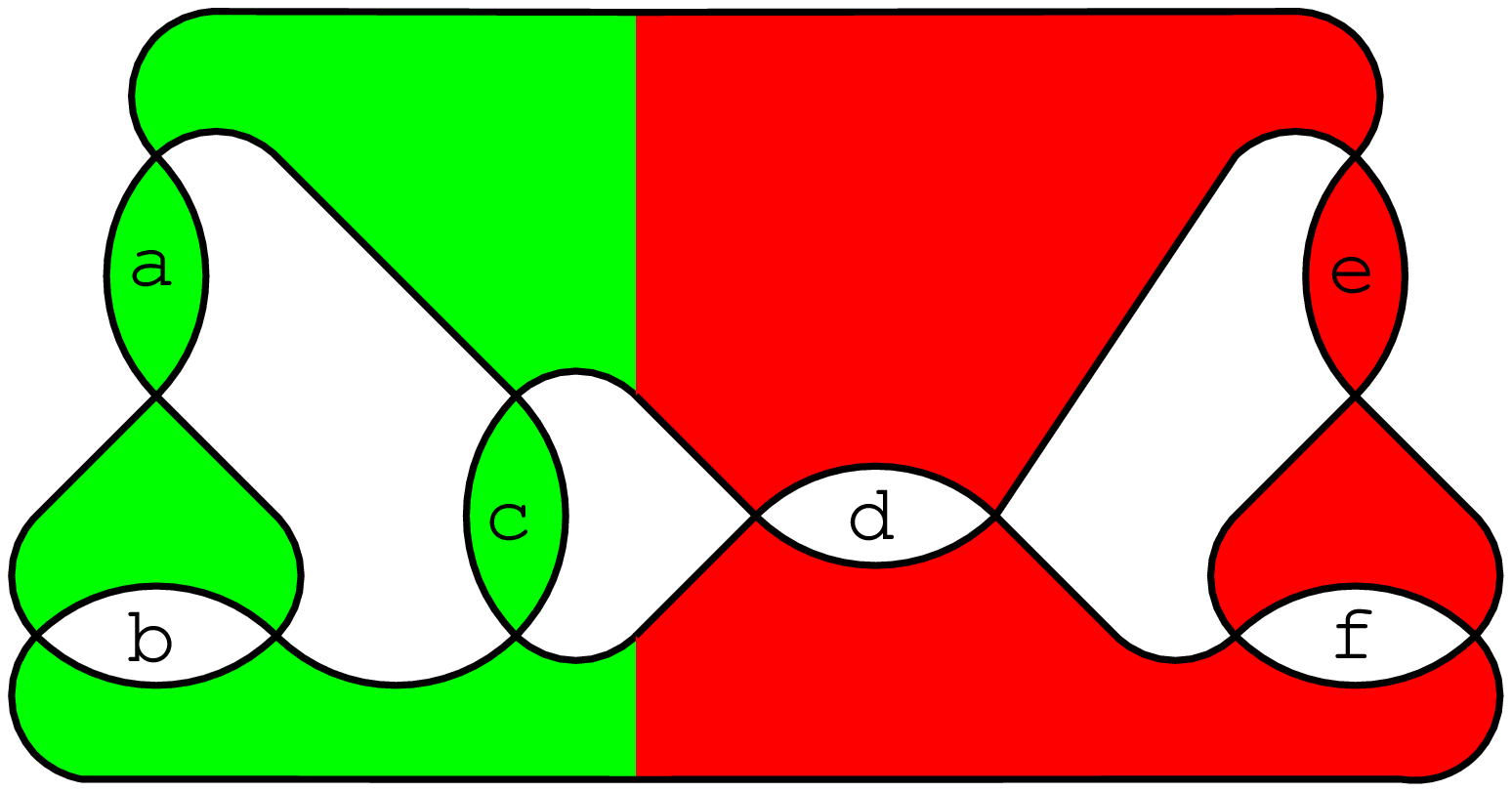}}

\caption{ \ }
$$
((a_1 a_2 + 1) a_3, a_2 (a_1 + a_3) + 1) \left( \begin{array}{cc}
0 & 1 \\
1 & 0
\end{array} \right)  \left( \begin{array}{cc}
a_5 a_6 + 1 \\
a_4 a_5 a_6 + a_4 + a_6
\end{array} \right) =
$$
$$
(a_1, 1) M \left( \begin{array}{cc}
1 & a_2 \\
a_2 & 0
\end{array} \right) M \left( \begin{array}{cc}
0 & a_3 \\
a_3 & 1
\end{array} \right) M \left( \begin{array}{cc}
0 & 1 \\
1 & a_4
\end{array} \right) M \left( \begin{array}{cc}
a_5 & 1 \\
1 & 0
\end{array} \right) M \left( \begin{array}{c}
1 \\
a_6
\end{array} \right)
$$
\end{figure}

\begin{figure}
\subsection{Families of knots with seed the knot $6_2$, the functions have eleven terms (eight cases)}

\centering

\psfrag{a}{\LARGE{$a_1$}}
\psfrag{b}{\LARGE{$a_2$}}
\psfrag{c}{\LARGE{$a_3$}}
\psfrag{d}{\LARGE{$a_4$}}
\psfrag{e}{\LARGE{$a_5$}}
\psfrag{f}{\LARGE{$a_6$}}

\scalebox{0.45}{\includegraphics{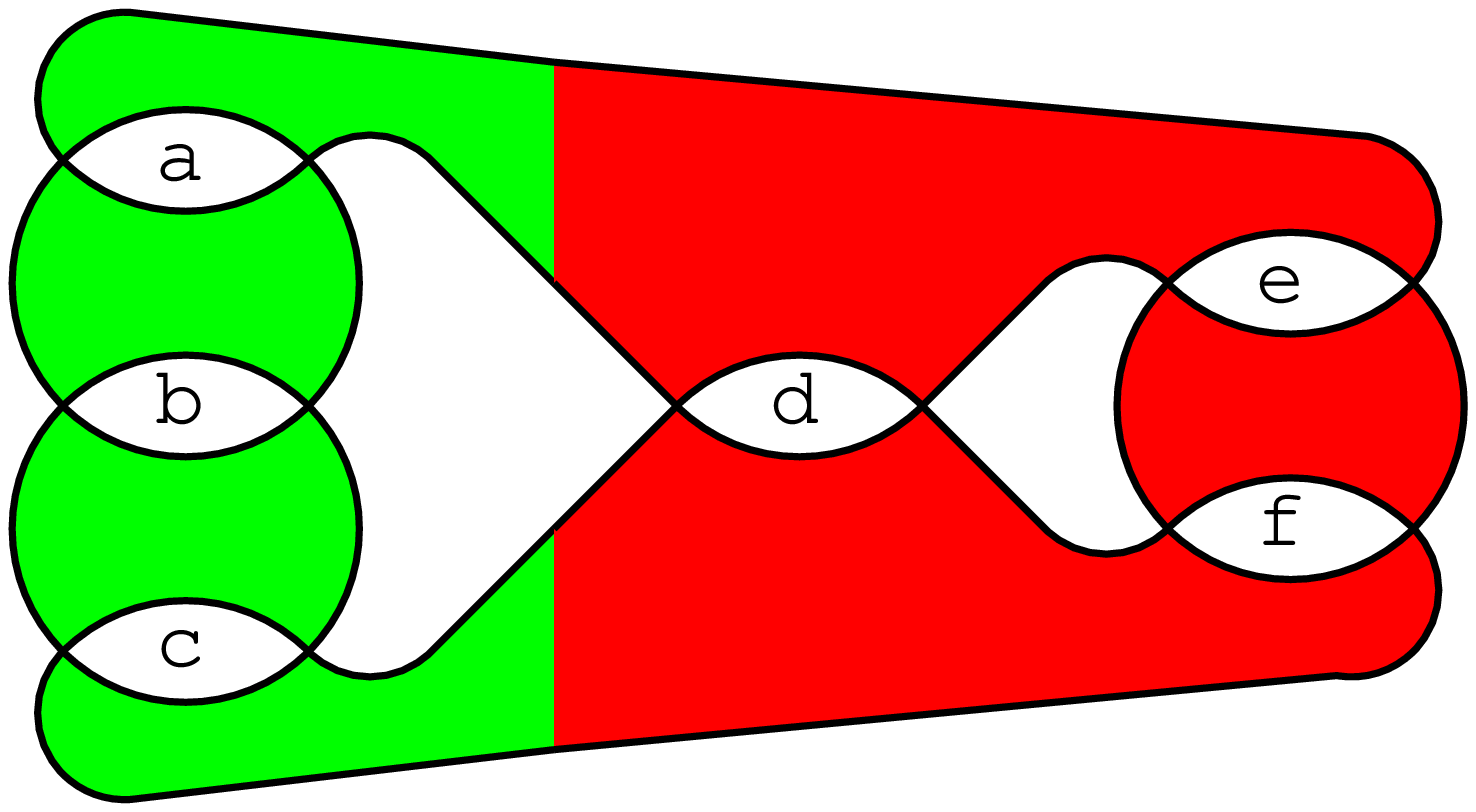}}

\caption{ \ }
$$
(a_1 a_2 + a_2 a_3 + a_3 a_1, a_1 a_2 a_3) \left( \begin{array}{cc}
0 & 1 \\
1 & 0
\end{array} \right)  \left( \begin{array}{cc}
a_5 + a_6 \\
a_4 a_5 + a_5 a_6 + a_6 a_4
\end{array} \right) =
$$
$$
(1, a_1) M \left( \begin{array}{cc}
1 & a_2 \\
a_2 & 0
\end{array} \right) M \left( \begin{array}{cc}
1 & a_3 \\
a_3 & 0
\end{array} \right) M \left( \begin{array}{cc}
0 & 1 \\
1 & a_4
\end{array} \right) M \left( \begin{array}{cc}
1 & a_5 \\
a_5 & 0
\end{array} \right) M \left( \begin{array}{c}
1 \\
a_6
\end{array} \right)
$$
\end{figure}

\begin{figure}

\centering

\psfrag{a}{\LARGE{$a_1$}}
\psfrag{b}{\LARGE{$a_2$}}
\psfrag{c}{\LARGE{$a_3$}}
\psfrag{d}{\LARGE{$a_4$}}
\psfrag{e}{\LARGE{$a_5$}}
\psfrag{f}{\LARGE{$a_6$}}

\scalebox{0.45}{\includegraphics{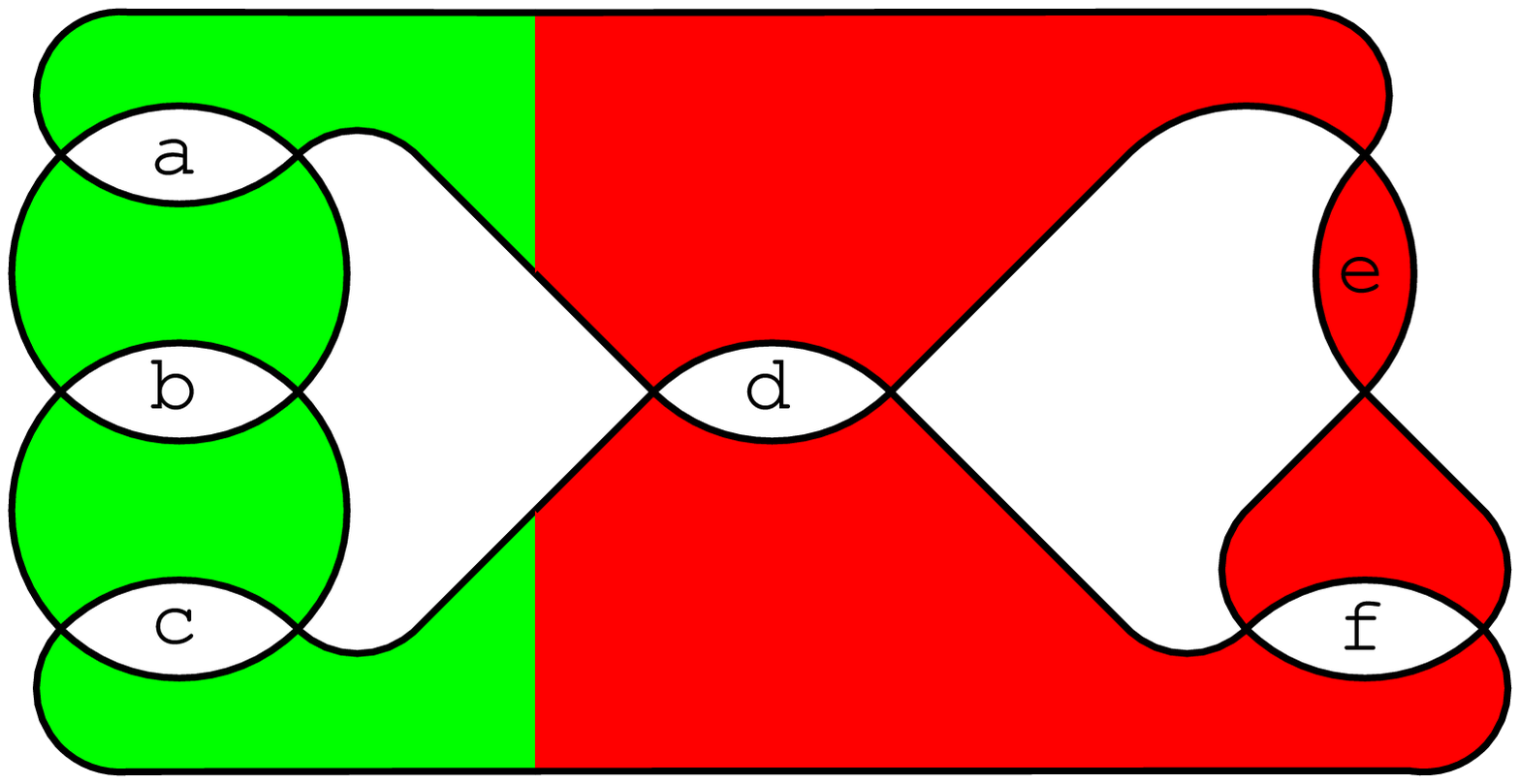}}
\caption{ \ }
$$
(a_1 a_2 + a_2 a_3 + a_3 a_1, a_1 a_2 a_3) \left( \begin{array}{cc}
0 & 1 \\
1 & 0
\end{array} \right)  \left( \begin{array}{cc}
a_5 a_6 + 1 \\
a_4 a_5 a_6 + a_4 + a_6
\end{array} \right) =
$$
$$
(1, a_1) M \left( \begin{array}{cc}
1 & a_2 \\
a_2 & 0
\end{array} \right) M \left( \begin{array}{cc}
1 & a_3 \\
a_3 & 0
\end{array} \right) M \left( \begin{array}{cc}
0 & 1 \\
1 & a_4
\end{array} \right) M \left( \begin{array}{cc}
a_5 & 1 \\
1 & 0
\end{array} \right) M \left( \begin{array}{c}
1 \\
a_6
\end{array} \right)
$$
\end{figure}

\begin{figure}38

\centering

\psfrag{a}{\LARGE{$a_1$}}
\psfrag{b}{\LARGE{$a_2$}}
\psfrag{c}{\LARGE{$a_3$}}
\psfrag{d}{\LARGE{$a_4$}}
\psfrag{e}{\LARGE{$a_5$}}
\psfrag{f}{\LARGE{$a_6$}}

\scalebox{0.4}{\includegraphics{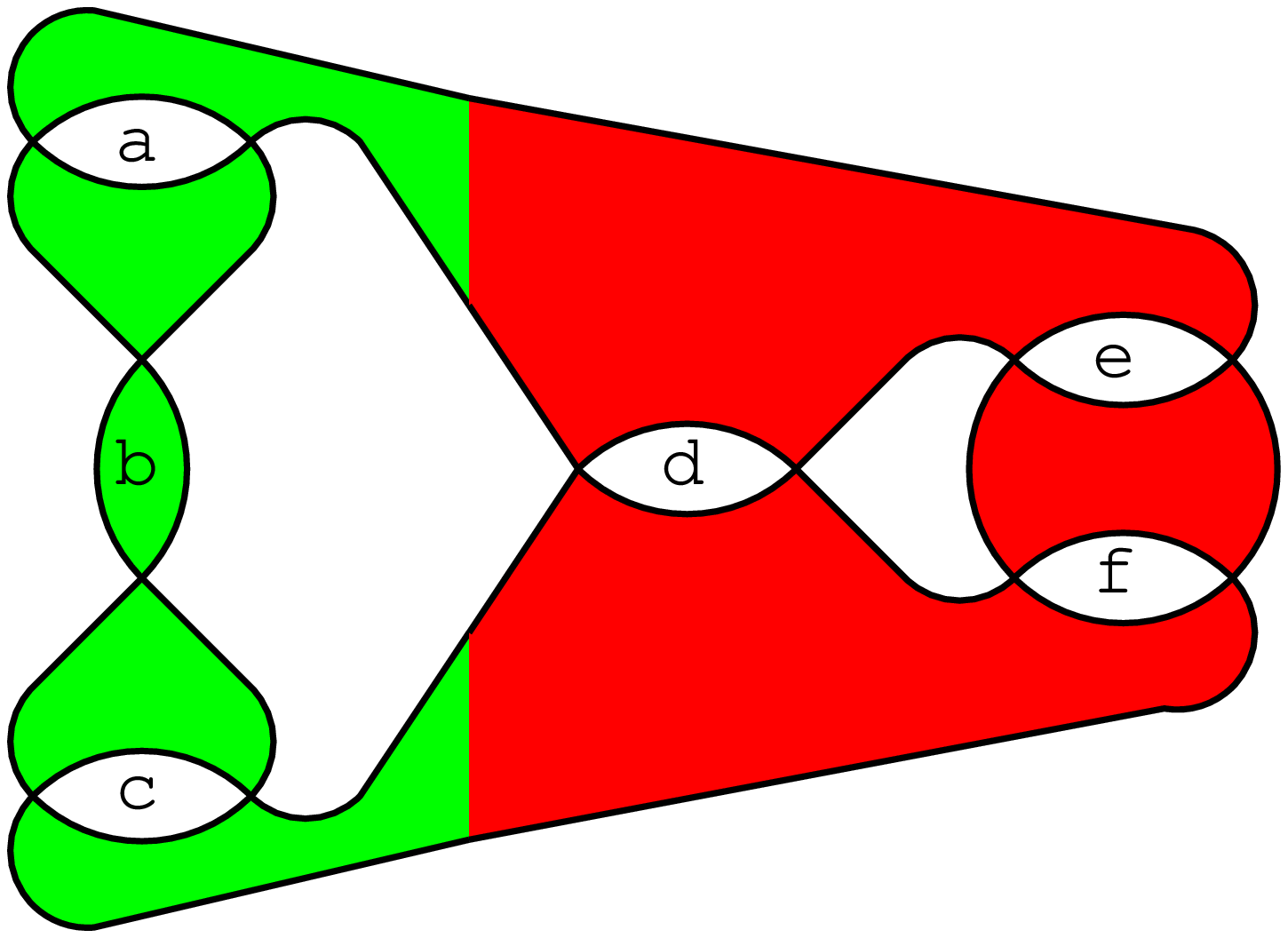}}
\caption{ \ }
$$
(a_1 a_2 a_3 + a_1 + a_3, a_1 a_3) \left( \begin{array}{cc}
0 & 1 \\
1 & 0
\end{array} \right)  \left( \begin{array}{cc}
a_5 + a_6 \\
a_4 a_5 + a_5 a_6 + a_6 a_4
\end{array} \right) =
$$
$$
(1, a_1) M \left( \begin{array}{cc}
a_2 & 1 \\
1 & 0
\end{array} \right) M \left( \begin{array}{cc}
1 & a_3 \\
a_3 & 0
\end{array} \right) M \left( \begin{array}{cc}
0 & 1 \\
1 & a_4
\end{array} \right) M \left( \begin{array}{cc}
1 & a_5 \\
a_5 & 0
\end{array} \right) M \left( \begin{array}{c}
1 \\
a_6
\end{array} \right)
$$
\end{figure}

\begin{figure}

\centering

\psfrag{a}{\LARGE{$a_1$}}
\psfrag{b}{\LARGE{$a_2$}}
\psfrag{c}{\LARGE{$a_3$}}
\psfrag{d}{\LARGE{$a_4$}}
\psfrag{e}{\LARGE{$a_5$}}
\psfrag{f}{\LARGE{$a_6$}}

\scalebox{0.4}{\includegraphics{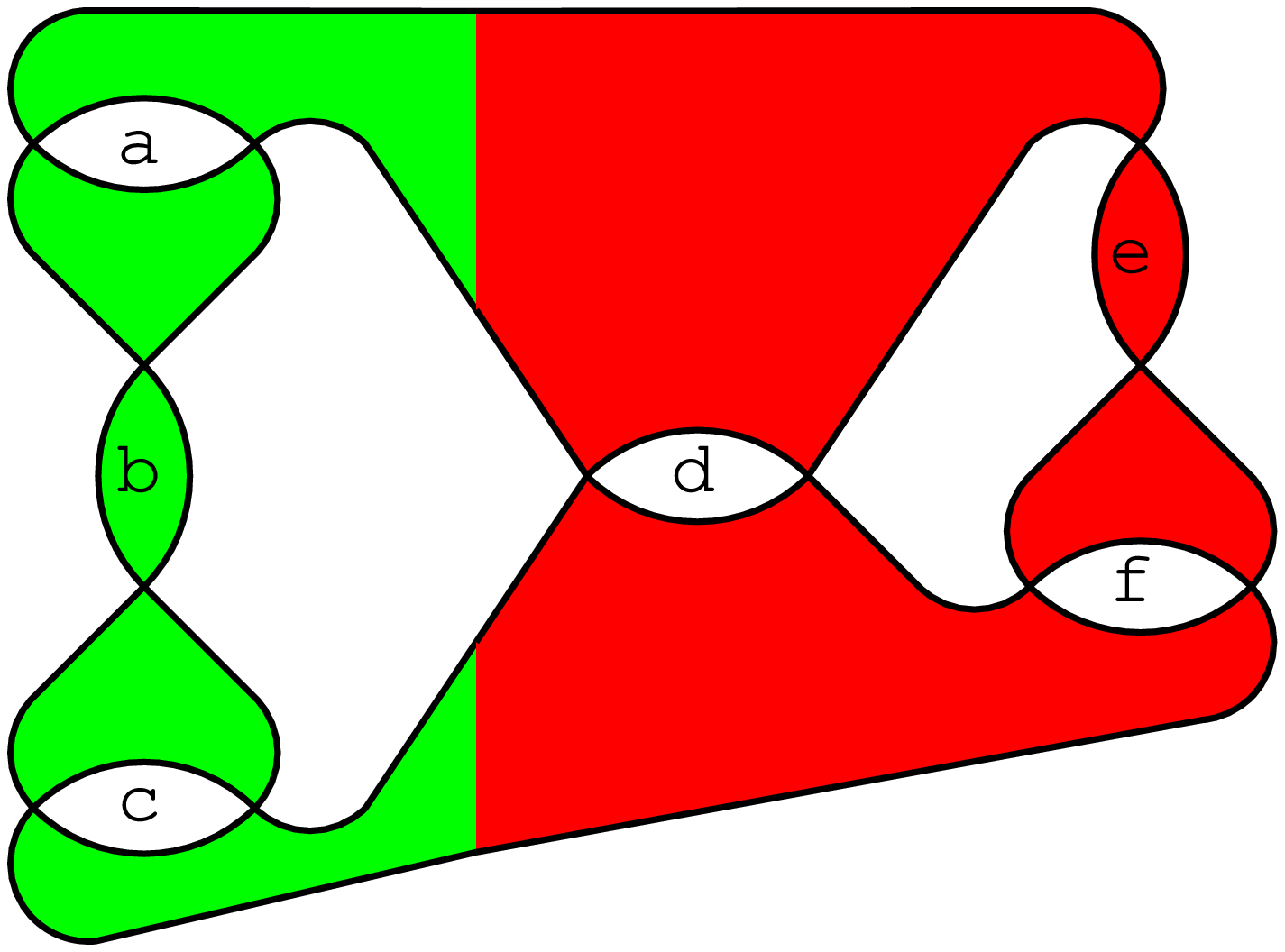}}
\caption{ \ }
$$
(a_1 a_2 a_3 + a_1 + a_3, a_1 a_3) \left( \begin{array}{cc}
0 & 1 \\
1 & 0
\end{array} \right)  \left( \begin{array}{cc}
a_5 a_6 + 1 \\
a_4 a_5 a_6 + a_4 + a_6
\end{array} \right) =
$$
$$
(1, a_1) M \left( \begin{array}{cc}
a_2 & 1 \\
1 & 0
\end{array} \right) M \left( \begin{array}{cc}
1 & a_3 \\
a_3 & 0
\end{array} \right) M \left( \begin{array}{cc}
0 & 1 \\
1 & a_4
\end{array} \right) M \left( \begin{array}{cc}
a_5 & 1 \\
1 & 0
\end{array} \right) M \left( \begin{array}{c}
1 \\
a_6
\end{array} \right)
$$
\end{figure}

\begin{figure}

\centering

\psfrag{a}{\LARGE{$a_1$}}
\psfrag{b}{\LARGE{$a_2$}}
\psfrag{c}{\LARGE{$a_3$}}
\psfrag{d}{\LARGE{$a_4$}}
\psfrag{e}{\LARGE{$a_5$}}
\psfrag{f}{\LARGE{$a_6$}}

\scalebox{0.5}{\includegraphics{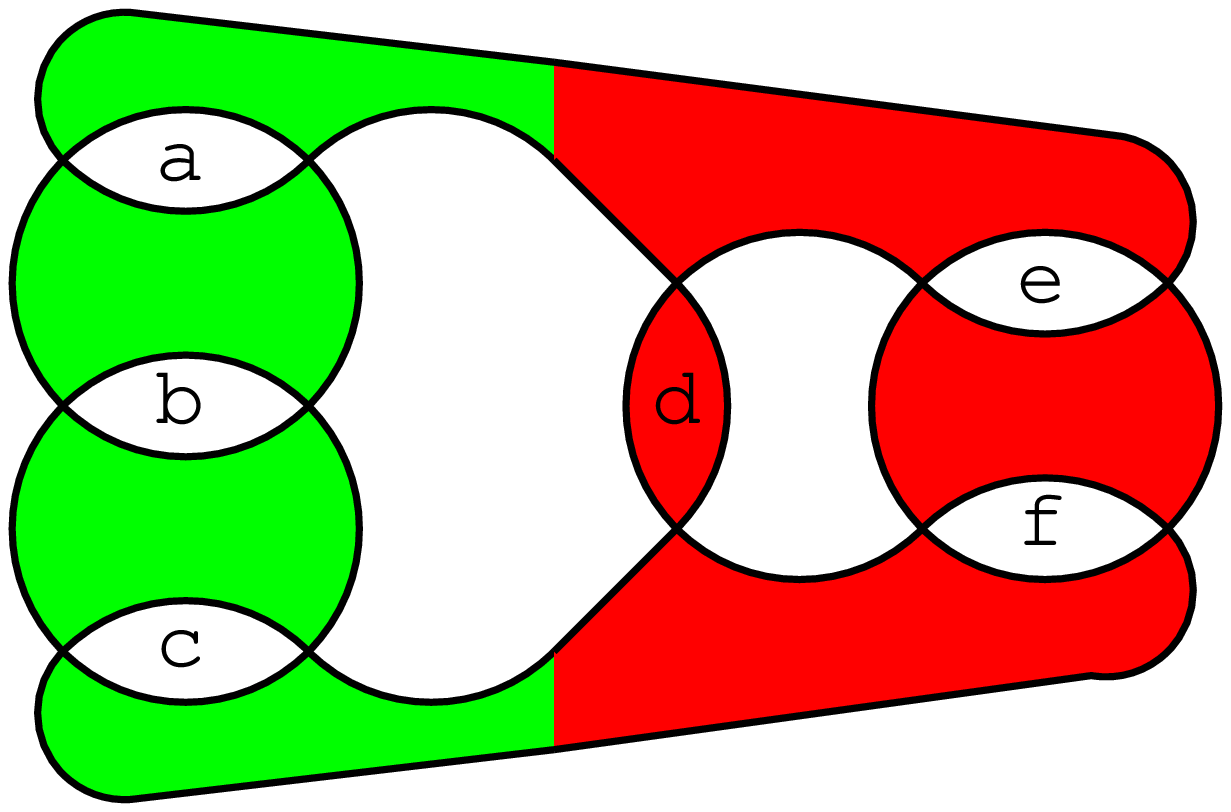}}
\caption{ \ }
$$
(a_1 a_2 + a_2 a_3 + a_3 a_1, a_1 a_2 a_3) \left( \begin{array}{cc}
0 & 1 \\
1 & 0
\end{array} \right)  \left( \begin{array}{cc}
a_4 (a_5 + a_6) \\
a_4 a_5 a_6 + a_5 + a_6
\end{array} \right) =
$$
$$
(1, a_1) M \left( \begin{array}{cc}
1 & a_2 \\
a_2 & 0
\end{array} \right) M \left( \begin{array}{cc}
1 & a_3 \\
a_3 & 0
\end{array} \right) M \left( \begin{array}{cc}
0 & a_4 \\
a_4 & 1
\end{array} \right) M \left( \begin{array}{cc}
1 & a_5 \\
a_5 & 0
\end{array} \right) M \left( \begin{array}{c}
1 \\
a_6
\end{array} \right)
$$
\end{figure}

\begin{figure}

\centering

\psfrag{a}{\LARGE{$a_1$}}
\psfrag{b}{\LARGE{$a_2$}}
\psfrag{c}{\LARGE{$a_3$}}
\psfrag{d}{\LARGE{$a_4$}}
\psfrag{e}{\LARGE{$a_5$}}
\psfrag{f}{\LARGE{$a_6$}}

\scalebox{0.5}{\includegraphics{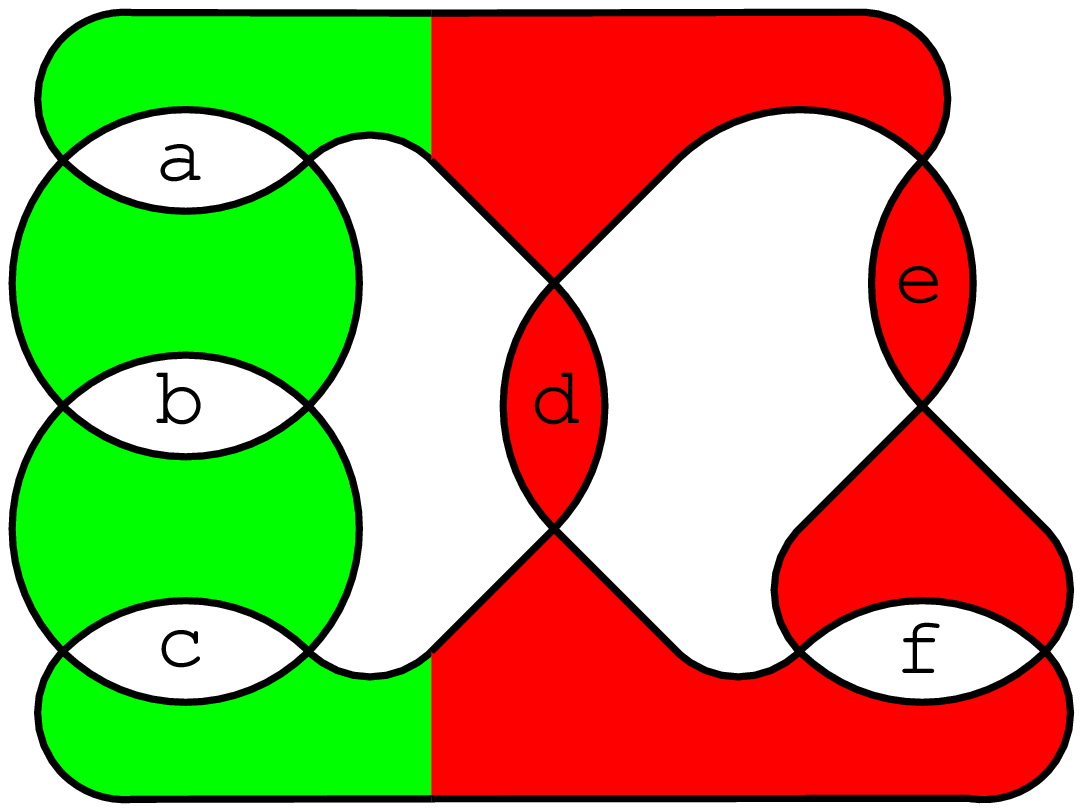}}
\caption{ \ }
$$
(a_1 a_2 + a_2 a_3 + a_3 a_1, a_1 a_2 a_3) \left( \begin{array}{cc}
0 & 1 \\
1 & 0
\end{array} \right)  \left( \begin{array}{cc}
a_4 (a_5 a_6 + 1) \\
(a_4 + a_5) a_6 + 1
\end{array} \right) =
$$
$$
(1, a_1) M \left( \begin{array}{cc}
1 & a_2 \\
a_2 & 0
\end{array} \right) M \left( \begin{array}{cc}
1 & a_3 \\
a_3 & 0
\end{array} \right) M \left( \begin{array}{cc}
0 & a_4 \\
a_4 & 1
\end{array} \right) M \left( \begin{array}{cc}
a_5 & 1 \\
1 & 0
\end{array} \right) M \left( \begin{array}{c}
1 \\
a_6
\end{array} \right)
$$
\end{figure}

\begin{figure}

\centering

\psfrag{a}{\LARGE{$a_1$}}
\psfrag{b}{\LARGE{$a_2$}}
\psfrag{c}{\LARGE{$a_3$}}
\psfrag{d}{\LARGE{$a_4$}}
\psfrag{e}{\LARGE{$a_5$}}
\psfrag{f}{\LARGE{$a_6$}}

\scalebox{0.4}{\includegraphics{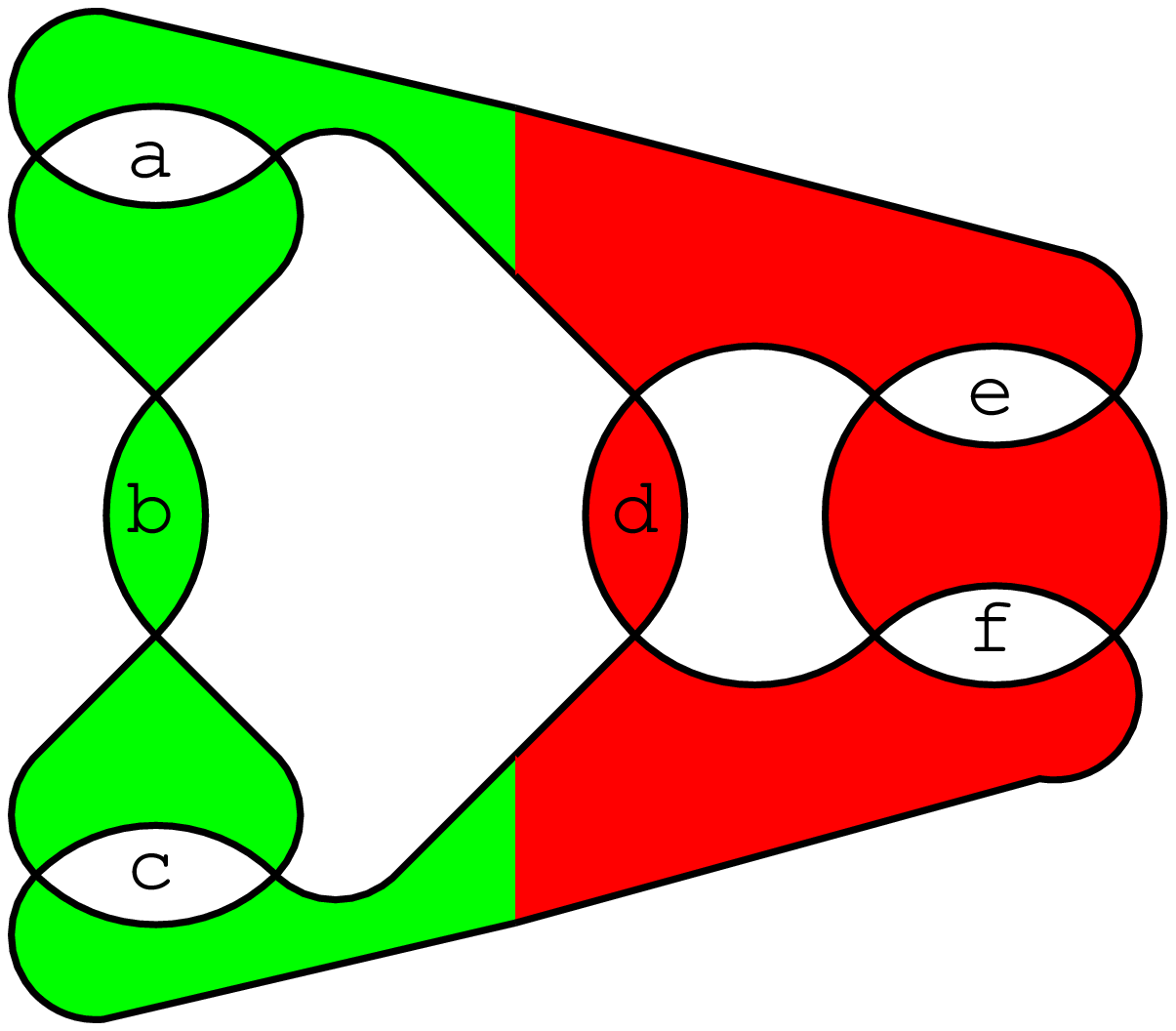}}
\caption{ \ }
$$
(a_1 a_2 a_3 + a_1 + a_3, a_1 a_3) \left( \begin{array}{cc}
0 & 1 \\
1 & 0
\end{array} \right)  \left( \begin{array}{cc}
a_4 (a_5 + a_6) \\
a_4 a_5 a_6 + a_5 + a_6
\end{array} \right) =
$$
$$
(1, a_1) M \left( \begin{array}{cc}
a_2 & 1 \\
1 & 0
\end{array} \right) M \left( \begin{array}{cc}
1 & a_3 \\
a_3 & 0
\end{array} \right) M \left( \begin{array}{cc}
0 & a_4 \\
a_4 & 1
\end{array} \right) M \left( \begin{array}{cc}
1 & a_5 \\
a_5 & 0
\end{array} \right) M \left( \begin{array}{c}
1 \\
a_6
\end{array} \right)
$$
\end{figure}

\begin{figure}

\centering

\psfrag{a}{\LARGE{$a_1$}}
\psfrag{b}{\LARGE{$a_2$}}
\psfrag{c}{\LARGE{$a_3$}}
\psfrag{d}{\LARGE{$a_4$}}
\psfrag{e}{\LARGE{$a_5$}}
\psfrag{f}{\LARGE{$a_6$}}

\scalebox{0.4}{\includegraphics{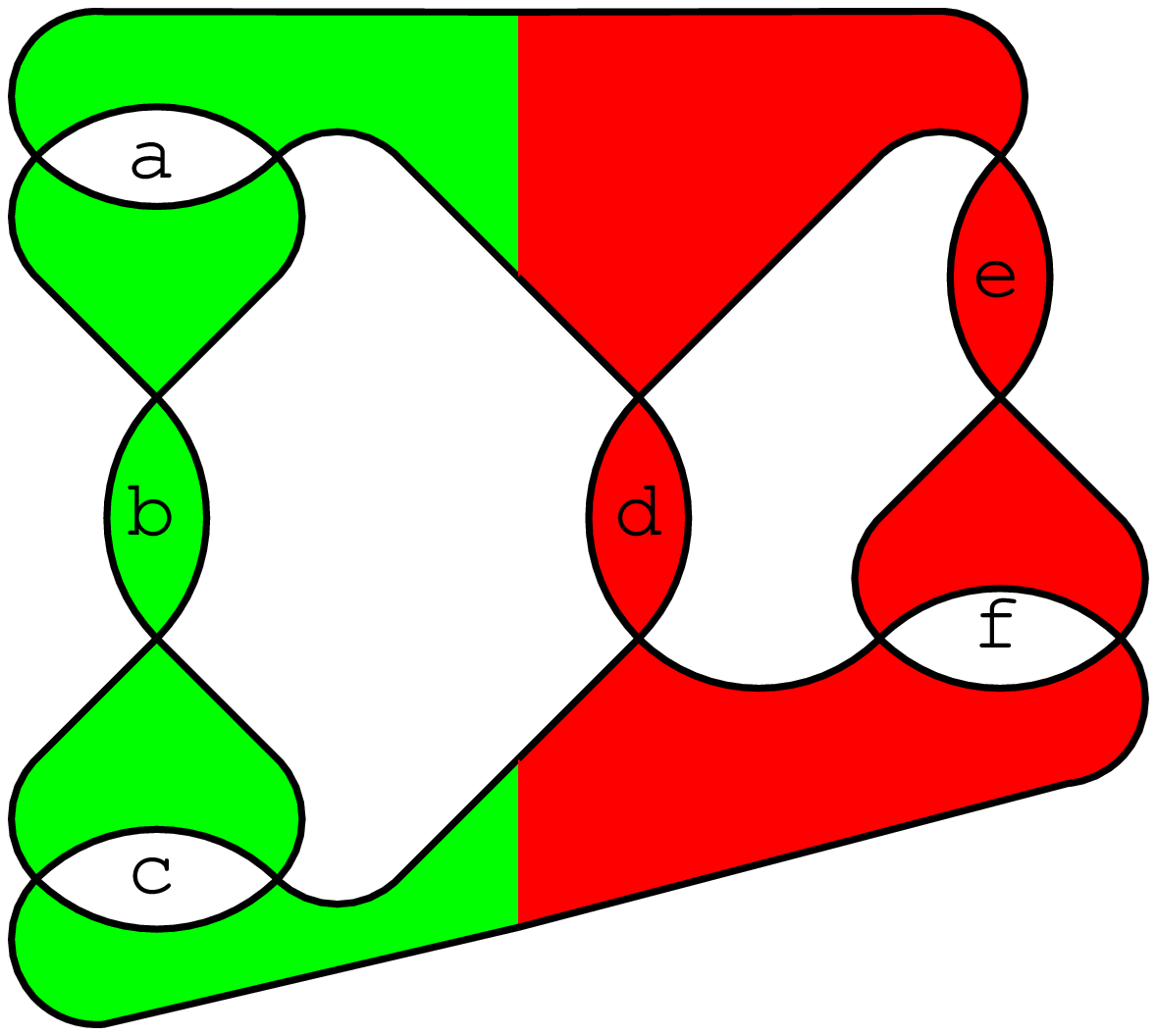}}
\caption{ \ }
$$
(a_1 a_2 a_3 + a_1 + a_3, a_1 a_3) \left( \begin{array}{cc}
0 & 1 \\
1 & 0
\end{array} \right)  \left( \begin{array}{cc}
a_4 (a_5 a_6 + 1) \\
(a_4 + a_5) a_6 + 1
\end{array} \right) =
$$
$$
(1, a_1) M \left( \begin{array}{cc}
a_2 & 1 \\
1 & 0
\end{array} \right) M \left( \begin{array}{cc}
1 & a_3 \\
a_3 & 0
\end{array} \right) M \left( \begin{array}{cc}
0 & a_4 \\
a_4 & 1
\end{array} \right) M \left( \begin{array}{cc}
a_5 & 1 \\
1 & 0
\end{array} \right) M \left( \begin{array}{c}
1 \\
a_6
\end{array} \right)
$$
\end{figure}


\begin{figure}
\setcounter{subsection}{5}
\subsection{Families of knots with seed the knot $6_3$, the function of Conway has thirteen terms (ten cases)}

\centering

\psfrag{b}{\LARGE{$a_1$}}
\psfrag{c}{\LARGE{$a_2$}}
\psfrag{a}{\LARGE{$a_3$}}
\psfrag{d}{\LARGE{$a_4$}}
\psfrag{e}{\LARGE{$a_5$}}
\psfrag{f}{\LARGE{$a_6$}}

\scalebox{0.35}{\includegraphics{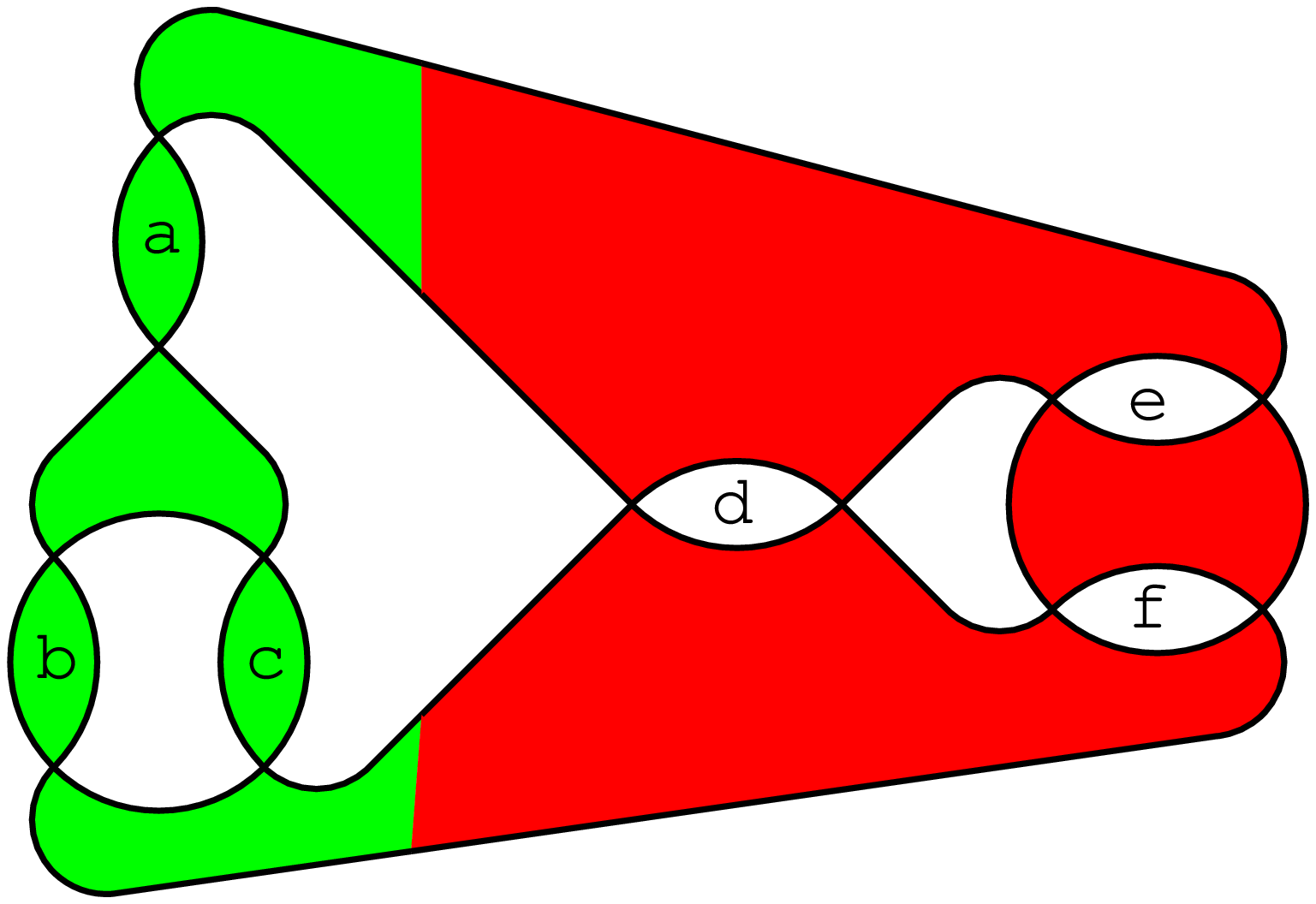}}
\caption{ \ }
$$
(a_1 a_2 + a_2 a_3 + a_3 a_1, a_1 + a_2) \left( \begin{array}{cc}
0 & 1 \\
1 & 0
\end{array} \right)  \left( \begin{array}{cc}
a_5 + a_6 \\
a_4 a_5 + a_5 a_6 + a_6 a_4
\end{array} \right) =
$$
$$
(a_1, 1) M \left( \begin{array}{cc}
0 & a_2 \\
a_2 & 1
\end{array} \right) M \left( \begin{array}{cc}
a_3 & 1 \\
1 & 0
\end{array} \right) M \left( \begin{array}{cc}
0 & 1 \\
1 & a_4
\end{array} \right) M \left( \begin{array}{cc}
1 & a_5 \\
a_5 & 0
\end{array} \right) M \left( \begin{array}{c}
1 \\
a_6
\end{array} \right)
$$
\end{figure}

\begin{figure}

\centering

\psfrag{b}{\LARGE{$a_1$}}
\psfrag{c}{\LARGE{$a_2$}}
\psfrag{a}{\LARGE{$a_3$}}
\psfrag{d}{\LARGE{$a_4$}}
\psfrag{e}{\LARGE{$a_5$}}
\psfrag{f}{\LARGE{$a_6$}}

\scalebox{0.35}{\includegraphics{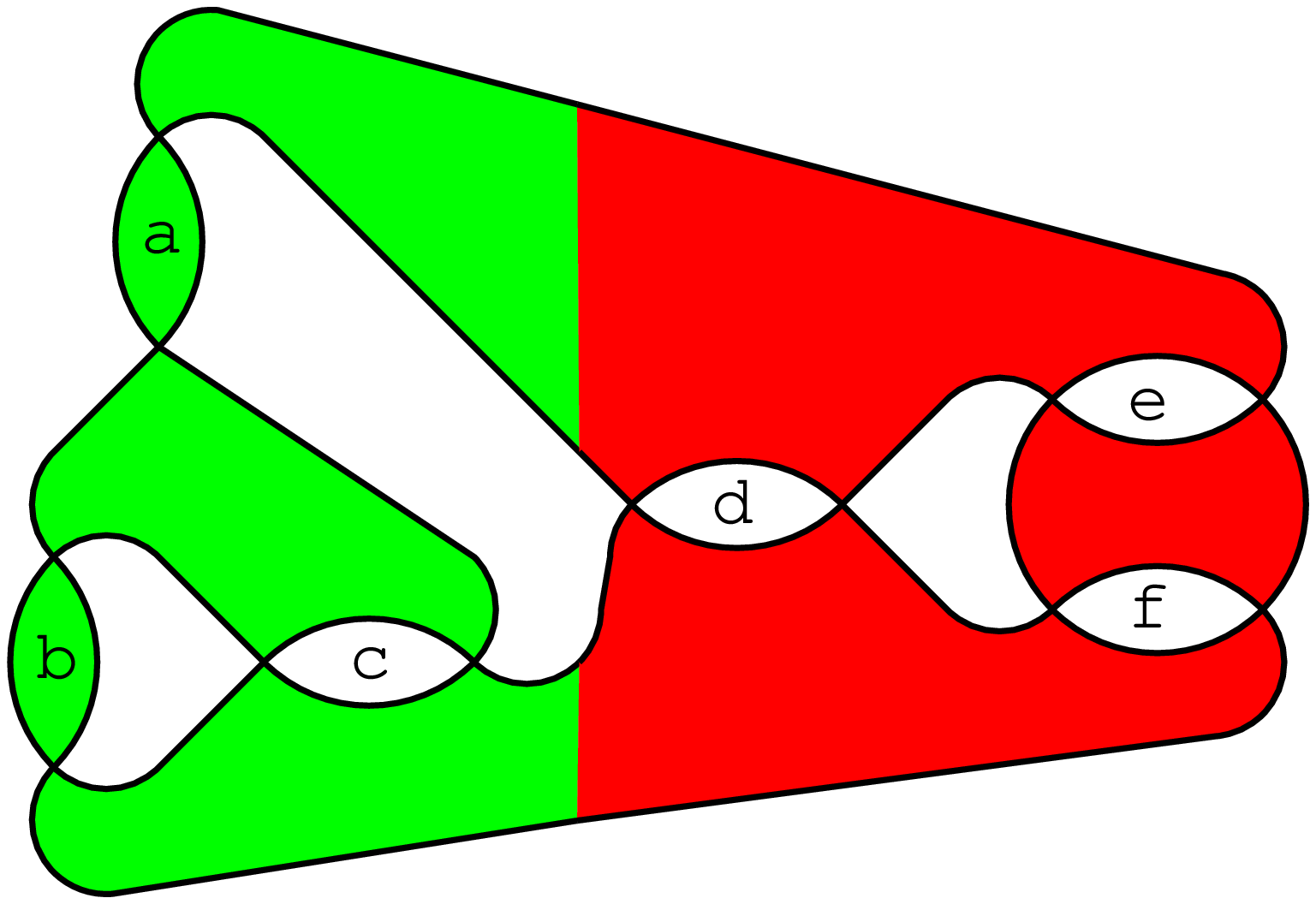}}
\caption{ \ }
$$
(a_1 a_2 a_3 + a_1 + a_3, a_1 a_2 + 1) \left( \begin{array}{cc}
0 & 1 \\
1 & 0
\end{array} \right)  \left( \begin{array}{cc}
a_5 + a_6 \\
a_4 a_5 + a_5 a_6 + a_6 a_4
\end{array} \right) =
$$
$$
(a_1, 1) M \left( \begin{array}{cc}
0 & 1 \\
1 & a_2
\end{array} \right) M \left( \begin{array}{cc}
a_3 & 1 \\
1 & 0
\end{array} \right) M \left( \begin{array}{cc}
0 & 1 \\
1 & a_4
\end{array} \right) M \left( \begin{array}{cc}
1 & a_5 \\
a_5 & 0
\end{array} \right) M \left( \begin{array}{c}
1 \\
a_6
\end{array} \right)
$$
\end{figure}

\begin{figure}

\centering

\psfrag{b}{\LARGE{$a_1$}}
\psfrag{c}{\LARGE{$a_2$}}
\psfrag{a}{\LARGE{$a_3$}}
\psfrag{d}{\LARGE{$a_4$}}
\psfrag{e}{\LARGE{$a_6$}}
\psfrag{f}{\LARGE{$a_5$}}

\scalebox{0.4}{\includegraphics{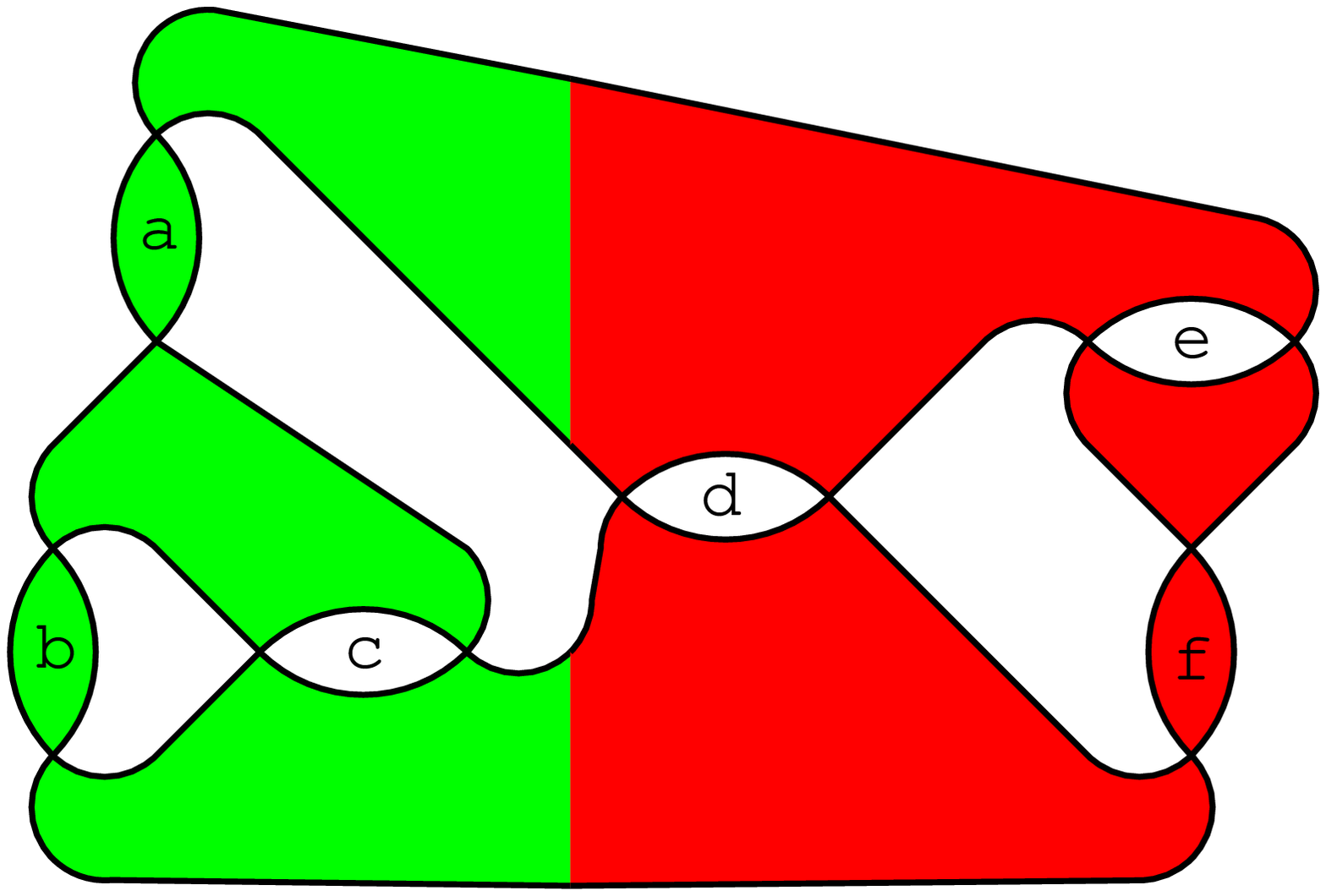}}

\caption{ \ }
$$
(a_1 a_2 a_3 + a_1 + a_3, a_1 a_2 + 1) \left( \begin{array}{cc}
0 & 1 \\
1 & 0
\end{array} \right)  \left( \begin{array}{cc}
a_5 a_6 + 1 \\
a_4 a_5 a_6 + a_4 + a_6
\end{array} \right) =
$$
$$
(a_1, 1) M \left( \begin{array}{cc}
0 & 1 \\
1 & a_2
\end{array} \right) M \left( \begin{array}{cc}
a_3 & 1 \\
1 & 0
\end{array} \right) M \left( \begin{array}{cc}
0 & 1 \\
1 & a_4
\end{array} \right) M \left( \begin{array}{cc}
a_5 & 1 \\
1 & 0
\end{array} \right) M \left( \begin{array}{c}
1 \\
a_6
\end{array} \right)
$$
\end{figure}

\begin{figure}

\centering

\psfrag{a}{\LARGE{$a_1$}}
\psfrag{b}{\LARGE{$a_2$}}
\psfrag{c}{\LARGE{$a_3$}}
\psfrag{d}{\LARGE{$a_4$}}
\psfrag{e}{\LARGE{$a_5$}}
\psfrag{f}{\LARGE{$a_6$}}

\scalebox{0.4}{\includegraphics{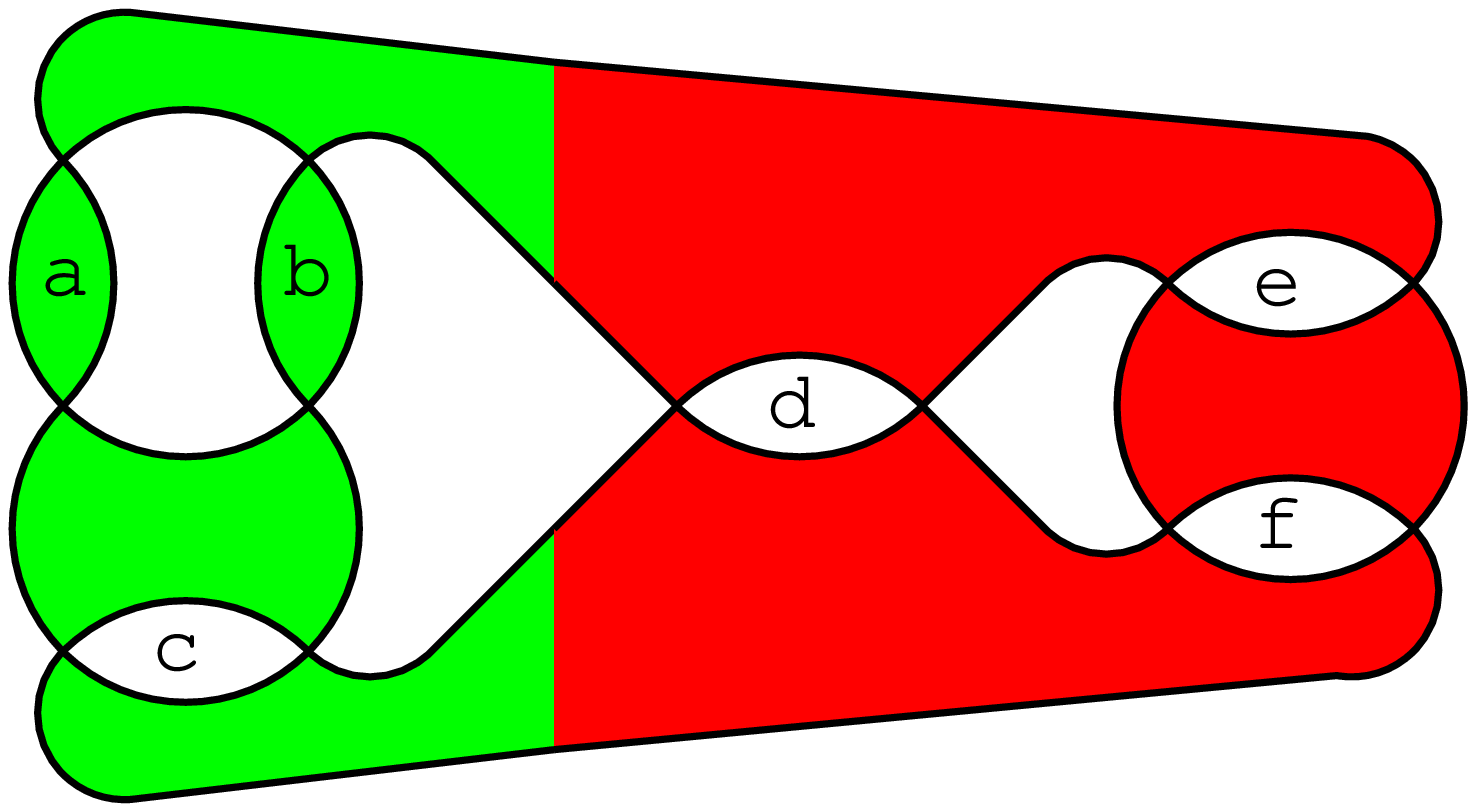}}

\caption{ \ }
$$
(a_1 a_2 a_3 + a_1 + a_2, a_3 (a_1 + a_2)) \left( \begin{array}{cc}
0 & 1 \\
1 & 0
\end{array} \right)  \left( \begin{array}{cc}
a_5 + a_6 \\
a_4 a_5 + a_5 a_6 + a_6 a_4
\end{array} \right) =
$$
$$
(a_1, 1) M \left( \begin{array}{cc}
0 & a_2 \\
a_2 & 1
\end{array} \right) M \left( \begin{array}{cc}
1 & a_3 \\
a_3 & 0
\end{array} \right) M \left( \begin{array}{cc}
0 & 1 \\
1 & a_4
\end{array} \right) M \left( \begin{array}{cc}
1 & a_5 \\
a_5 & 0
\end{array} \right) M \left( \begin{array}{c}
1 \\
a_6
\end{array} \right)
$$
\end{figure}

\begin{figure}

\centering

\psfrag{a}{\LARGE{$a_1$}}
\psfrag{b}{\LARGE{$a_2$}}
\psfrag{c}{\LARGE{$a_3$}}
\psfrag{d}{\LARGE{$a_4$}}
\psfrag{e}{\LARGE{$a_5$}}
\psfrag{f}{\LARGE{$a_6$}}

\scalebox{0.5}{\includegraphics{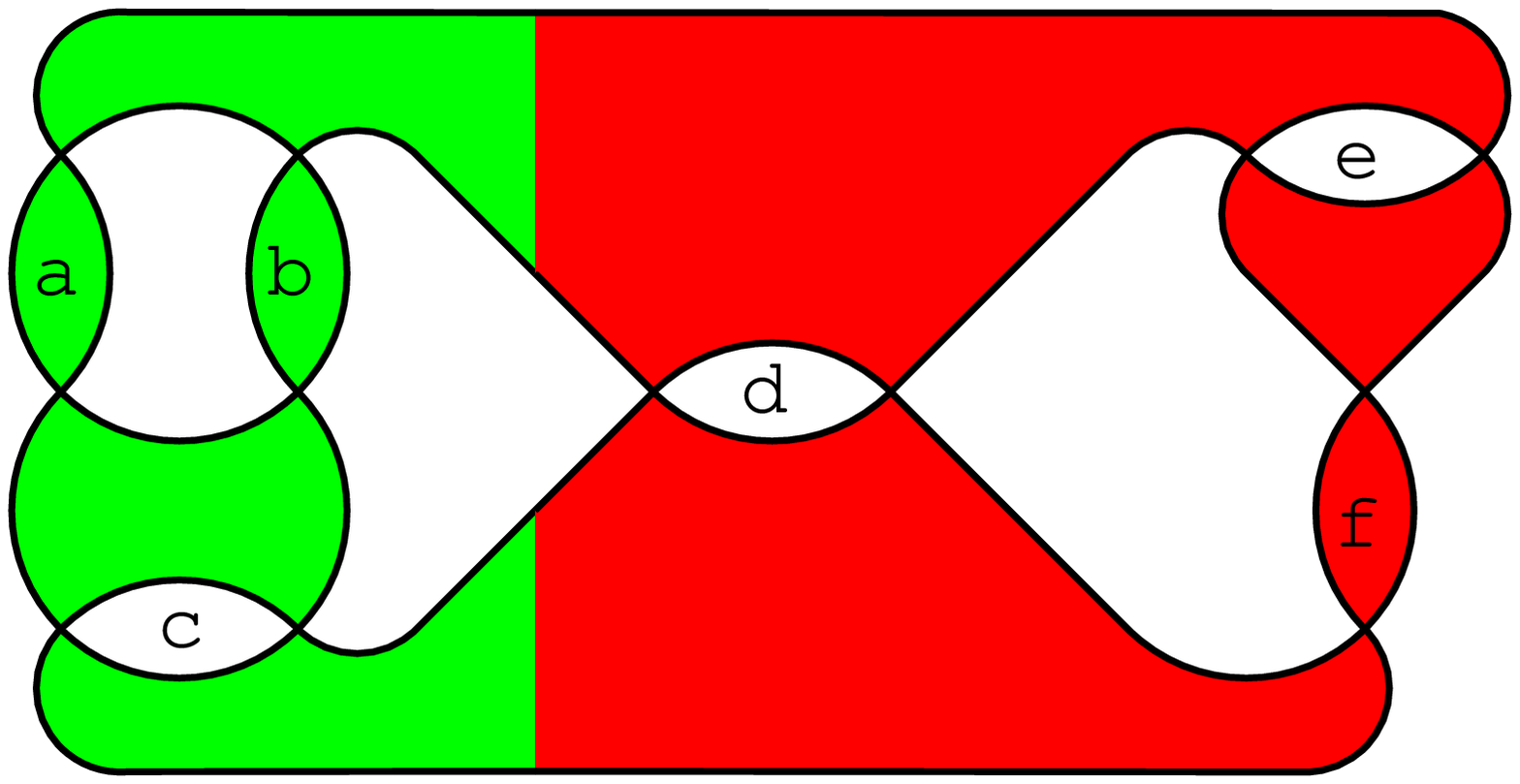}}

\caption{ \ }
$$
(a_1 a_2 a_3 + a_1 + a_2, a_3 (a_1 + a_2)) \left( \begin{array}{cc}
0 & 1 \\
1 & 0
\end{array} \right)  \left( \begin{array}{cc}
a_5 a_6 + 1 \\
a_4 a_5 a_6 + a_4 + a_5
\end{array} \right) =
$$
$$
(a_1, 1) M \left( \begin{array}{cc}
0 & a_2 \\
a_2 & 1
\end{array} \right) M \left( \begin{array}{cc}
1 & a_3 \\
a_3 & 0
\end{array} \right) M \left( \begin{array}{cc}
0 & 1 \\
1 & a_4
\end{array} \right) M \left( \begin{array}{cc}
1 & a_5 \\
a_5 & 0
\end{array} \right) M \left( \begin{array}{c}
a_6 \\
1
\end{array} \right)
$$
\end{figure}

\begin{figure}

\centering

\psfrag{a}{\LARGE{$a_1$}}
\psfrag{b}{\LARGE{$a_2$}}
\psfrag{c}{\LARGE{$a_3$}}
\psfrag{d}{\LARGE{$a_4$}}
\psfrag{e}{\LARGE{$a_5$}}
\psfrag{f}{\LARGE{$a_6$}}

\scalebox{0.5}{\includegraphics{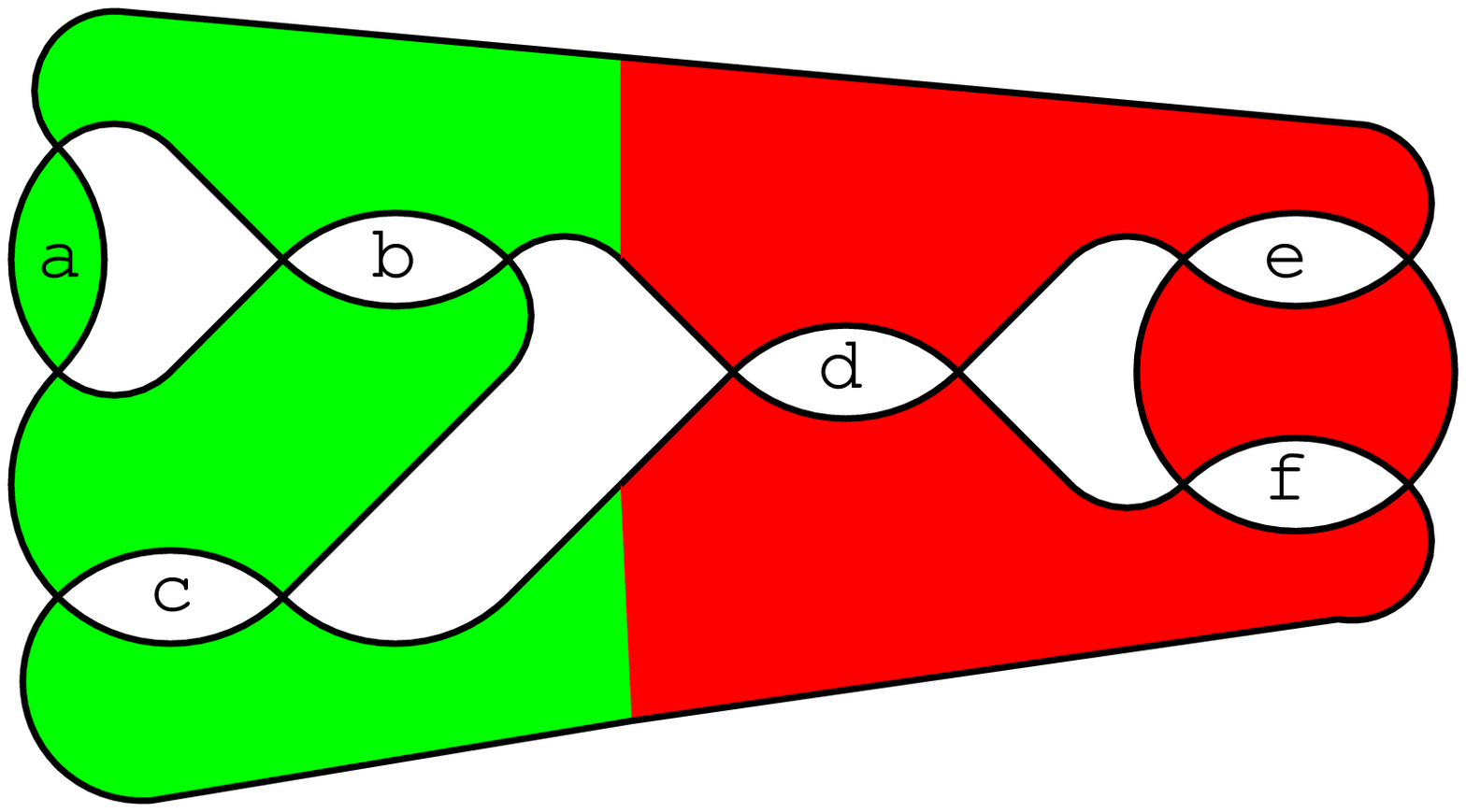}}

\caption{ \ }
$$
(a_1 (a_2 + a_3) + 1, a_3 (a_1 a_2 + 1)) \left( \begin{array}{cc}
0 & 1 \\
1 & 0
\end{array} \right)  \left( \begin{array}{cc}
a_5 + a_6 \\
a_4 a_5 + a_5 a_6 + a_6 a_4
\end{array} \right) =
$$
$$
(a_1, 1) M \left( \begin{array}{cc}
0 & 1 \\
1 & a_2
\end{array} \right) M \left( \begin{array}{cc}
1 & a_3 \\
a_3 & 0
\end{array} \right) M \left( \begin{array}{cc}
0 & 1 \\
1 & a_4
\end{array} \right) M \left( \begin{array}{cc}
1 & a_5 \\
a_5 & 0
\end{array} \right) M \left( \begin{array}{c}
1 \\
a_6
\end{array} \right)
$$
\end{figure}

\begin{figure}

\centering

\psfrag{a}{\LARGE{$a_1$}}
\psfrag{b}{\LARGE{$a_2$}}
\psfrag{c}{\LARGE{$a_3$}}
\psfrag{d}{\LARGE{$a_4$}}
\psfrag{e}{\LARGE{$a_5$}}
\psfrag{f}{\LARGE{$a_6$}}

\scalebox{0.5}{\includegraphics{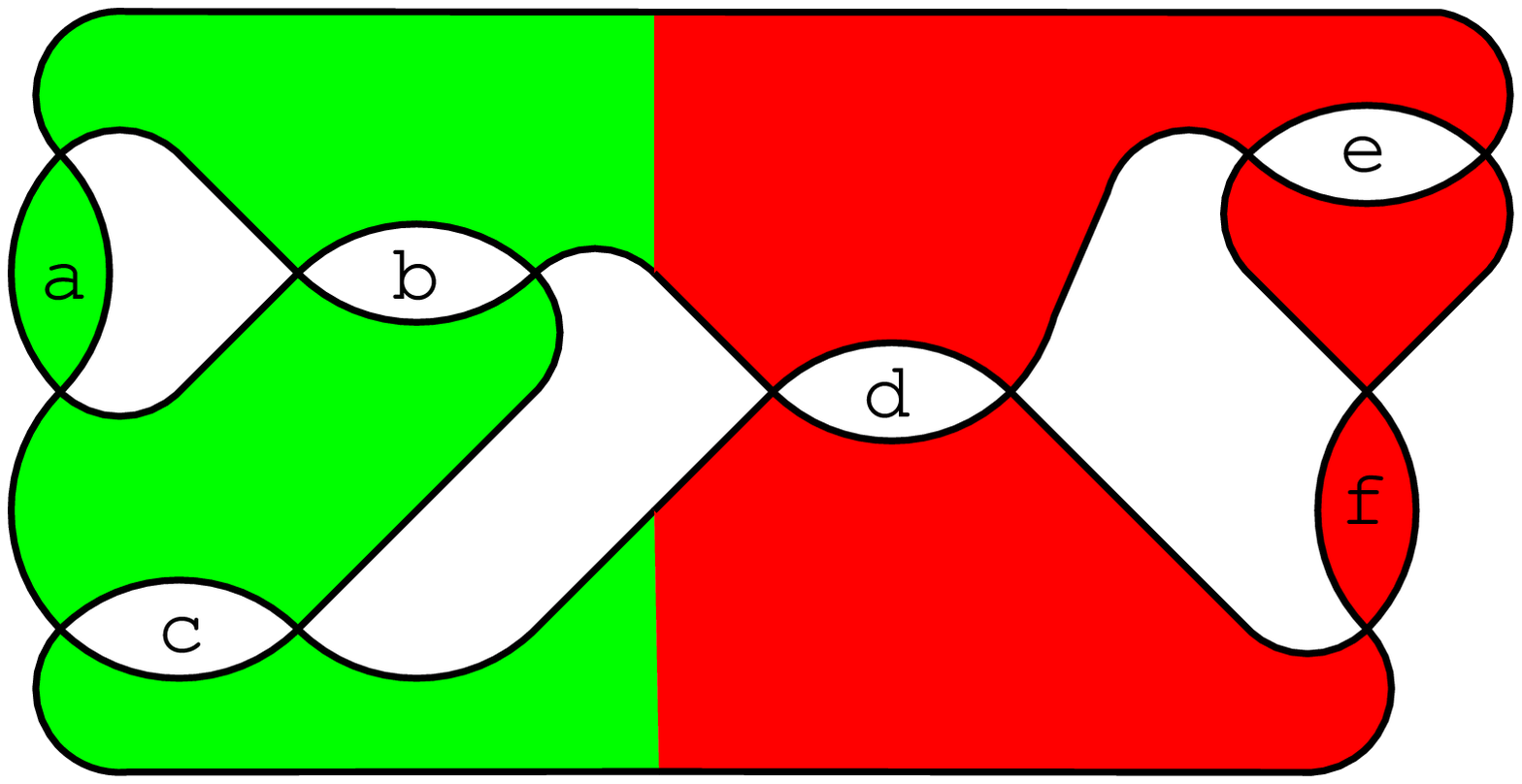}}

\caption{ \ }
$$
(a_1 (a_2 + a_3) + 1, a_3 (a_1 a_2 + 1)) \left( \begin{array}{cc}
0 & 1 \\
1 & 0
\end{array} \right)  \left( \begin{array}{cc}
a_5 a_6 + 1 \\
a_4 a_5 a_6 + a_4 + a_5
\end{array} \right) =
$$
$$
(a_1, 1) M \left( \begin{array}{cc}
0 & 1 \\
1 & a_2
\end{array} \right) M \left( \begin{array}{cc}
1 & a_3 \\
a_3 & 0
\end{array} \right) M \left( \begin{array}{cc}
0 & 1 \\
1 & a_4
\end{array} \right) M \left( \begin{array}{cc}
1 & a_5 \\
a_5 & 0
\end{array} \right) M \left( \begin{array}{c}
a_6 \\
1
\end{array} \right)
$$
\end{figure}

\begin{figure}

\centering

\psfrag{a}{\LARGE{$a_1$}}
\psfrag{b}{\LARGE{$a_2$}}
\psfrag{c}{\LARGE{$a_3$}}
\psfrag{d}{\LARGE{$a_4$}}
\psfrag{e}{\LARGE{$a_5$}}
\psfrag{f}{\LARGE{$a_6$}}

\scalebox{0.5}{\includegraphics{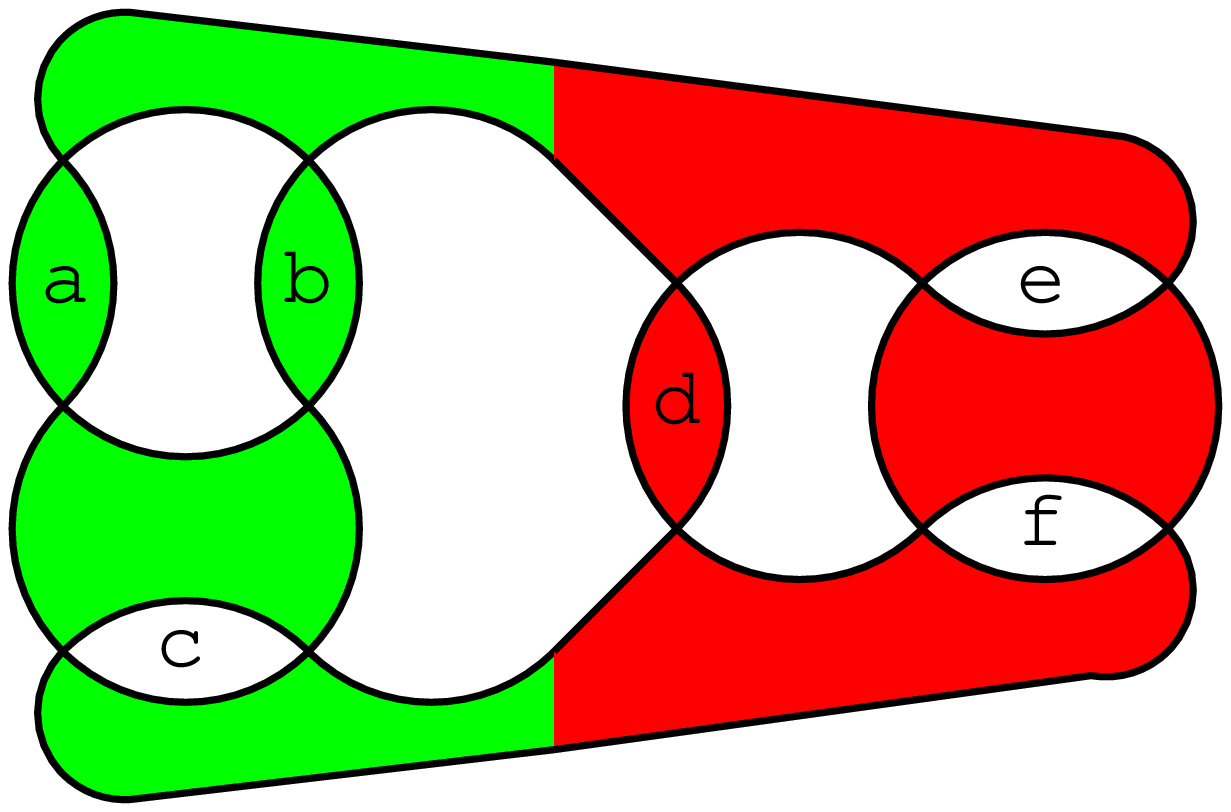}}

\caption{ \ }
$$
(a_1 a_2 a_3 + a_1 + a_2, a_3 (a_1 + a_2)) \left( \begin{array}{cc}
0 & 1 \\
1 & 0
\end{array} \right)  \left( \begin{array}{cc}
a_4 (a_5 + a_6) \\
a_4 a_5 a_6 + a_5 + a_6
\end{array} \right) =
$$
$$
(a_1, 1) M \left( \begin{array}{cc}
0 & a_2 \\
a_2 & 1
\end{array} \right) M \left( \begin{array}{cc}
1 & a_3 \\
a_3 & 0
\end{array} \right) M \left( \begin{array}{cc}
0 & a_4 \\
a_4 & 1
\end{array} \right) M \left( \begin{array}{cc}
1 & a_5 \\
a_5 & 0
\end{array} \right) M \left( \begin{array}{c}
1 \\
a_6
\end{array} \right)
$$
\end{figure}

\begin{figure}

\centering

\psfrag{a}{\LARGE{$a_1$}}
\psfrag{b}{\LARGE{$a_2$}}
\psfrag{c}{\LARGE{$a_3$}}
\psfrag{d}{\LARGE{$a_4$}}
\psfrag{e}{\LARGE{$a_5$}}
\psfrag{f}{\LARGE{$a_6$}}

\scalebox{0.5}{\includegraphics{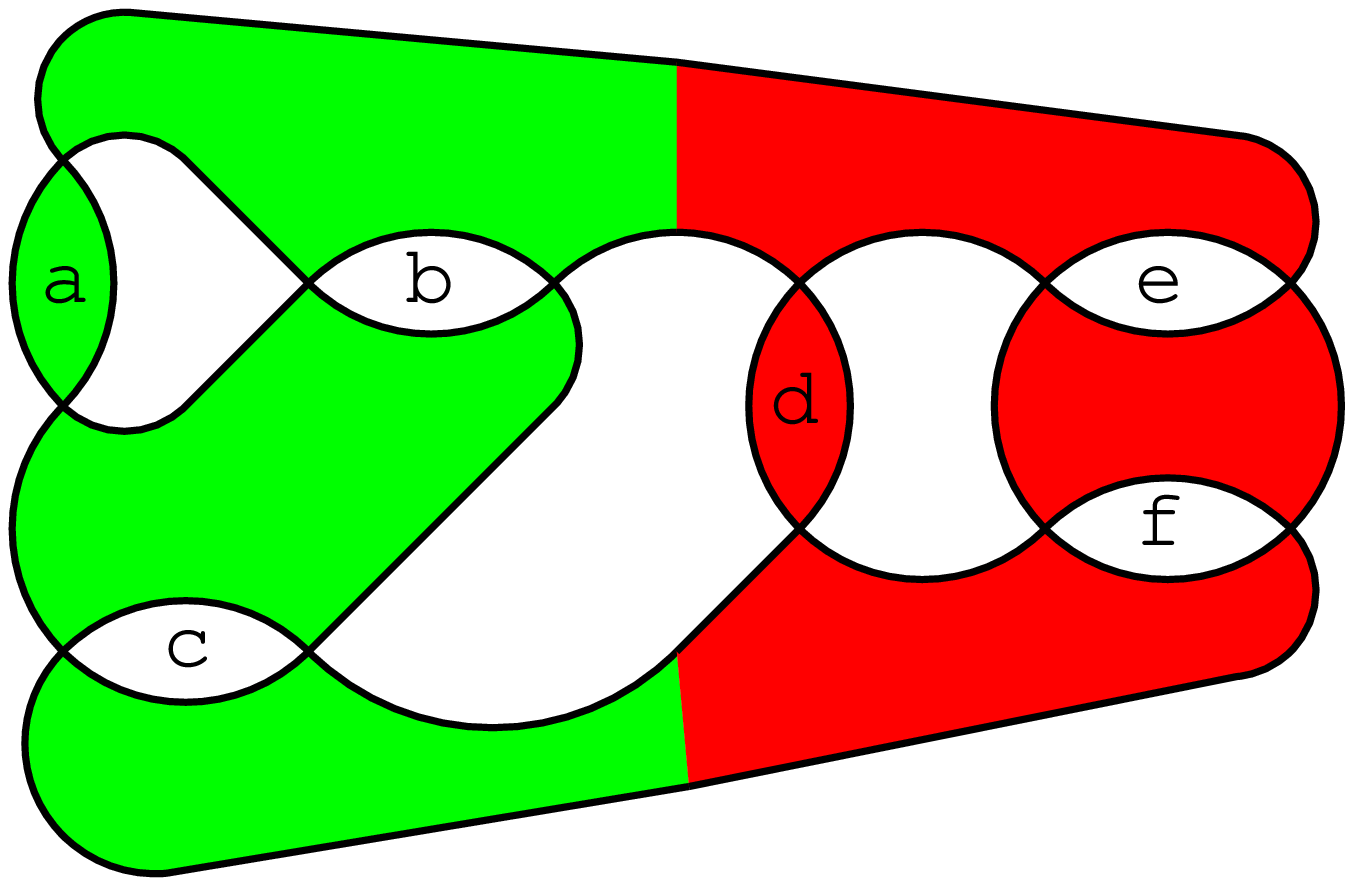}}

\caption{ \ }
$$
(a_1 (a_2 + a_3) + 1, a_3 (a_1 a_2 + 1)) \left( \begin{array}{cc}
0 & 1 \\
1 & 0
\end{array} \right)  \left( \begin{array}{cc}
a_4 (a_5 + a_6) \\
a_4 a_5 a_6 + a_5 + a_6
\end{array} \right) =
$$
$$
(a_1, 1) M \left( \begin{array}{cc}
0 & 1 \\
1 & a_2
\end{array} \right) M \left( \begin{array}{cc}
1 & a_3 \\
a_3 & 0
\end{array} \right) M \left( \begin{array}{cc}
0 & a_4 \\
a_4 & 1
\end{array} \right) M \left( \begin{array}{cc}
1 & a_5 \\
a_5 & 0
\end{array} \right) M \left( \begin{array}{c}
1 \\
a_6
\end{array} \right)
$$
\end{figure}

\begin{figure}

\centering

\psfrag{a}{\LARGE{$a_1$}}
\psfrag{b}{\LARGE{$a_2$}}
\psfrag{c}{\LARGE{$a_3$}}
\psfrag{d}{\LARGE{$a_4$}}
\psfrag{e}{\LARGE{$a_5$}}
\psfrag{f}{\LARGE{$a_6$}}

\scalebox{0.5}{\includegraphics{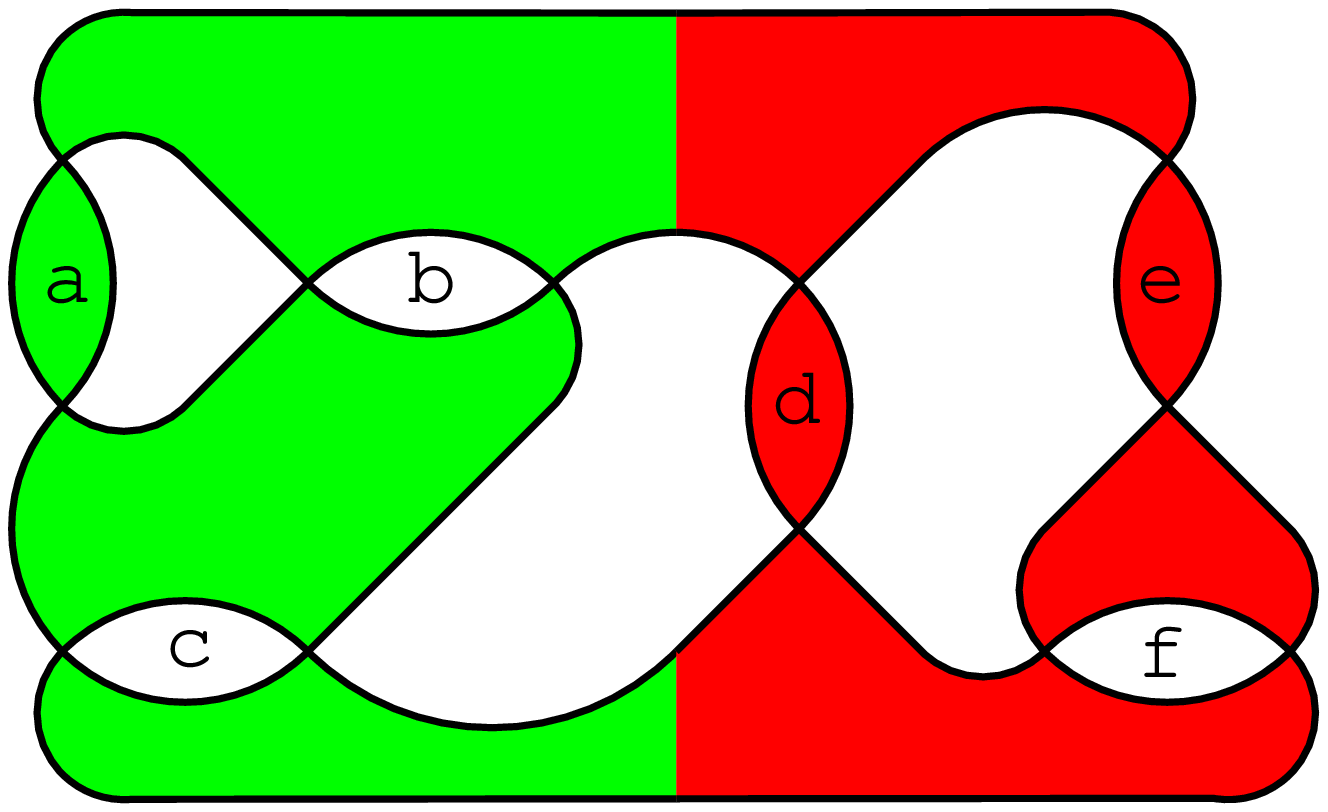}}

\caption{ \ }
$$
(a_1 (a_2 + a_3) + 1, a_3 (a_1 a_2 + 1)) \left( \begin{array}{cc}
0 & 1 \\
1 & 0
\end{array} \right)  \left( \begin{array}{cc}
a_4 (a_5 a_6 + 1) \\
(a_4 + a_5) a_6 + 1
\end{array} \right) =
$$
$$
(a_1, 1) M \left( \begin{array}{cc}
0 & 1 \\
1 & a_2
\end{array} \right) M \left( \begin{array}{cc}
1 & a_3 \\
a_3 & 0
\end{array} \right) M \left( \begin{array}{cc}
0 & a_4 \\
a_4 & 1
\end{array} \right) M \left( \begin{array}{cc}
a_5 & 1 \\
1 & 0
\end{array} \right) M \left( \begin{array}{c}
1 \\
a_6
\end{array} \right)
$$
\end{figure}

\begin{figure}

\subsection{Families of knots with seed the link  $6_1^3$ (four cases)}
\centering

\psfrag{a}{\LARGE{$a_2$}}
\psfrag{b}{\LARGE{$a_1$}}
\psfrag{c}{\LARGE{$a_4$}}
\psfrag{d}{\LARGE{$a_3$}}
\psfrag{e}{\LARGE{$a_6$}}
\psfrag{f}{\LARGE{$a_5$}}

\scalebox{0.5}{\includegraphics{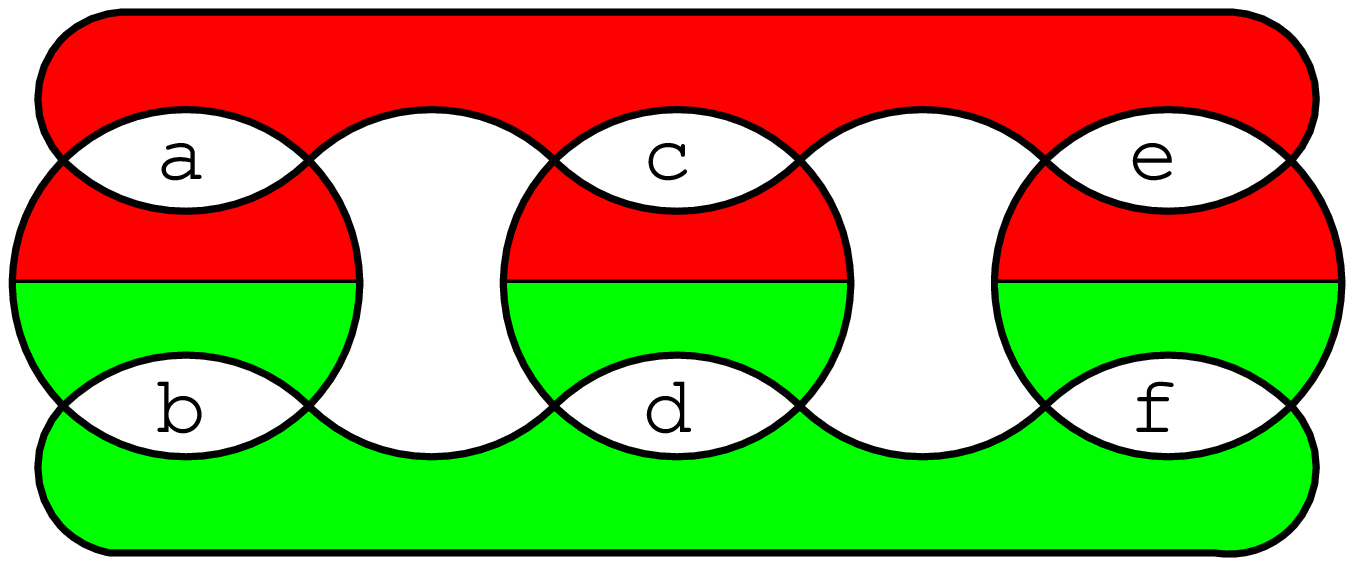}}
\caption{ \ }
$$
\left( \begin{array}{c}
a_1 a_3 a_5 \\
a_3 a_5 \\
a_5 a_1 \\
a_1 a_3 \\
a_1 + a_3 + a_5
\end{array} \right)^{\rm T}
\left( \begin{array}{ccccc}
0 & 0 & 0 & 0 & 1 \\
0 & 0 & 1 & 1 & 0 \\
0 & 1 & 0 & 1 & 0 \\
0 & 1 & 1 & 0 & 0 \\
1 & 0 & 0 & 0 & 0
\end{array} \right)
\left( \begin{array}{c}
a_2 a_4 a_6 \\
a_4 a_6 \\
a_6 a_2 \\
a_2 a_4 \\
a_2 + a_4 + a_6
\end{array} \right)
$$

\end{figure}

\begin{figure}

\centering

\psfrag{a}{\LARGE{$a_2$}}
\psfrag{b}{\LARGE{$a_1$}}
\psfrag{c}{\LARGE{$a_4$}}
\psfrag{d}{\LARGE{$a_3$}}
\psfrag{e}{\LARGE{$a_6$}}
\psfrag{f}{\LARGE{$a_5$}}

\scalebox{0.5}{\includegraphics{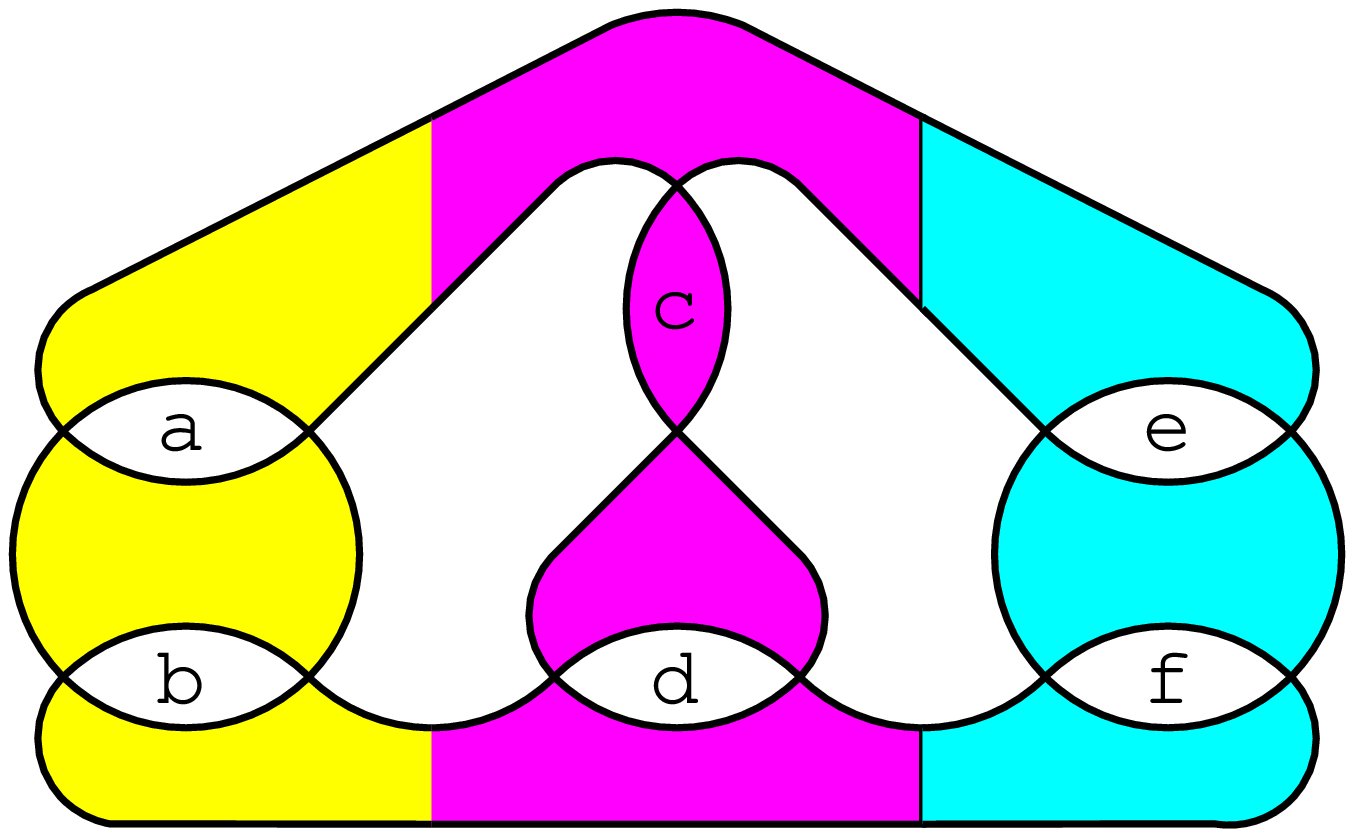}}

\caption{ \ }
$$
(a_1 + a_2, a_1 a_2) \left( \begin{array}{cc}
0 & 1 \\
1 & 0
\end{array} \right)  \left( \begin{array}{cc}
0 & a_3 a_4 + 1 \\
a_3 a_4 + 1 & a_3
\end{array} \right) \left( \begin{array}{cc}
0 & 1 \\
1 & 0
\end{array} \right) \left( \begin{array}{cc}
a_5 + a_6 \\
a_5 a_6
\end{array} \right) =
$$
$$
(1, a_1) M \left( \begin{array}{cc}
1 & a_2 \\
a_2 & 0
\end{array} \right) M  \left( \begin{array}{cc}
0 & a_3 a_4 + 1 \\
a_3 a_4 + 1 & a_3
\end{array} \right) M \left( \begin{array}{cc}
1 & a_5 \\
a_5 & 0
\end{array} \right) M \left( \begin{array}{cc}
1 \\
a_6
\end{array} \right)
$$
\end{figure}

\begin{figure}

\centering

\psfrag{a}{\LARGE{$a_2$}}
\psfrag{b}{\LARGE{$a_1$}}
\psfrag{c}{\LARGE{$a_4$}}
\psfrag{d}{\LARGE{$a_3$}}
\psfrag{e}{\LARGE{$a_6$}}
\psfrag{f}{\LARGE{$a_5$}}

\scalebox{0.5}{\includegraphics{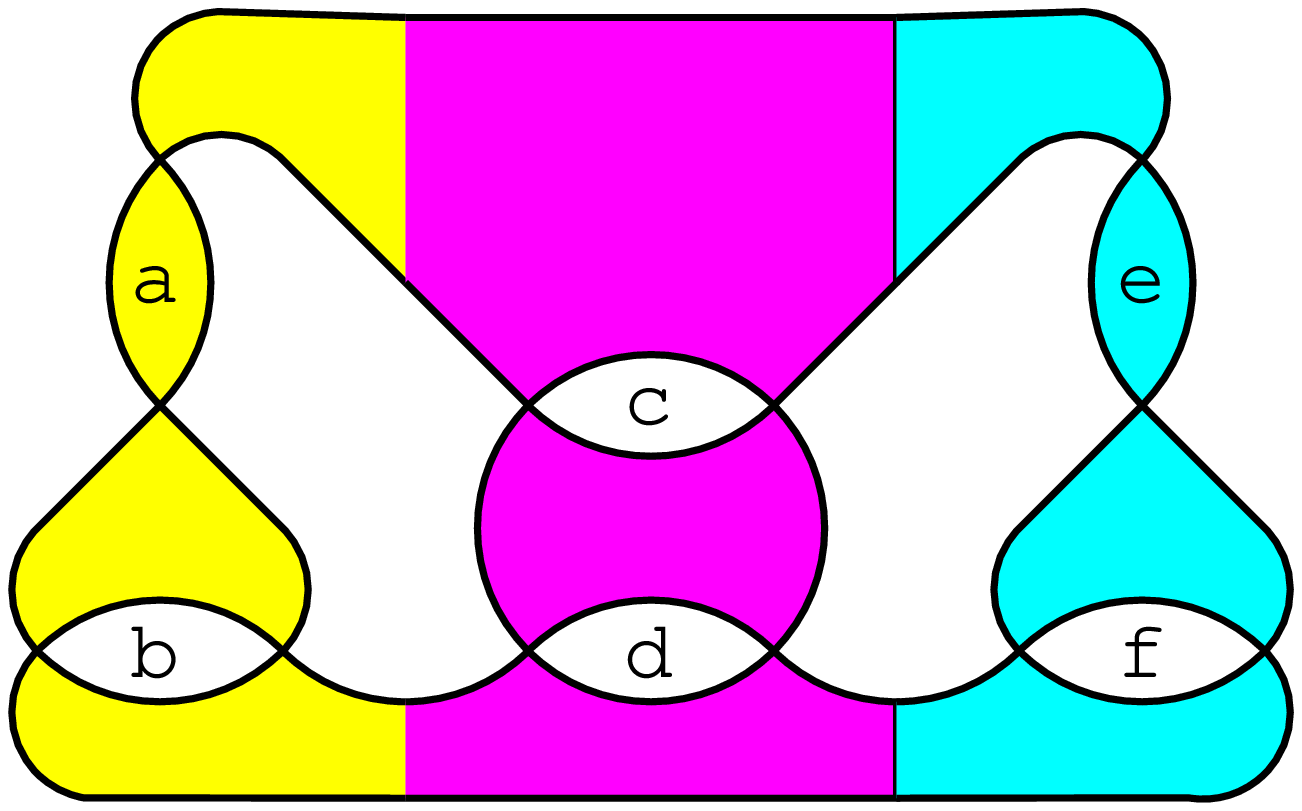}}

\caption{ \ }
$$
(a_1 a_2 + 1, a_1) \left( \begin{array}{cc}
0 & 1 \\
1 & 0
\end{array} \right)  \left( \begin{array}{cc}
0 & a_3 + a_4 \\
a_3 + a_4 & a_3 a_4
\end{array} \right) \left( \begin{array}{cc}
0 & 1 \\
1 & 0
\end{array} \right) \left( \begin{array}{cc}
a_5 a_6 + 1 \\
a_5
\end{array} \right) =
$$
$$
(1, a_1) M \left( \begin{array}{cc}
a_2 & 1 \\
1 & 0
\end{array} \right) M \left( \begin{array}{cc}
0 & a_3 + a_4 \\
a_3 + a_4 & a_3 a_4
\end{array} \right) M \left( \begin{array}{cc}
a_6 & 1 \\
1 & 0
\end{array} \right) M \left( \begin{array}{cc}
1 \\
a_5
\end{array} \right)
$$
\end{figure}

\begin{figure}

\centering

\psfrag{a}{\LARGE{$a_2$}}
\psfrag{b}{\LARGE{$a_1$}}
\psfrag{c}{\LARGE{$a_4$}}
\psfrag{d}{\LARGE{$a_3$}}
\psfrag{e}{\LARGE{$a_6$}}
\psfrag{f}{\LARGE{$a_5$}}

\scalebox{0.5}{\includegraphics{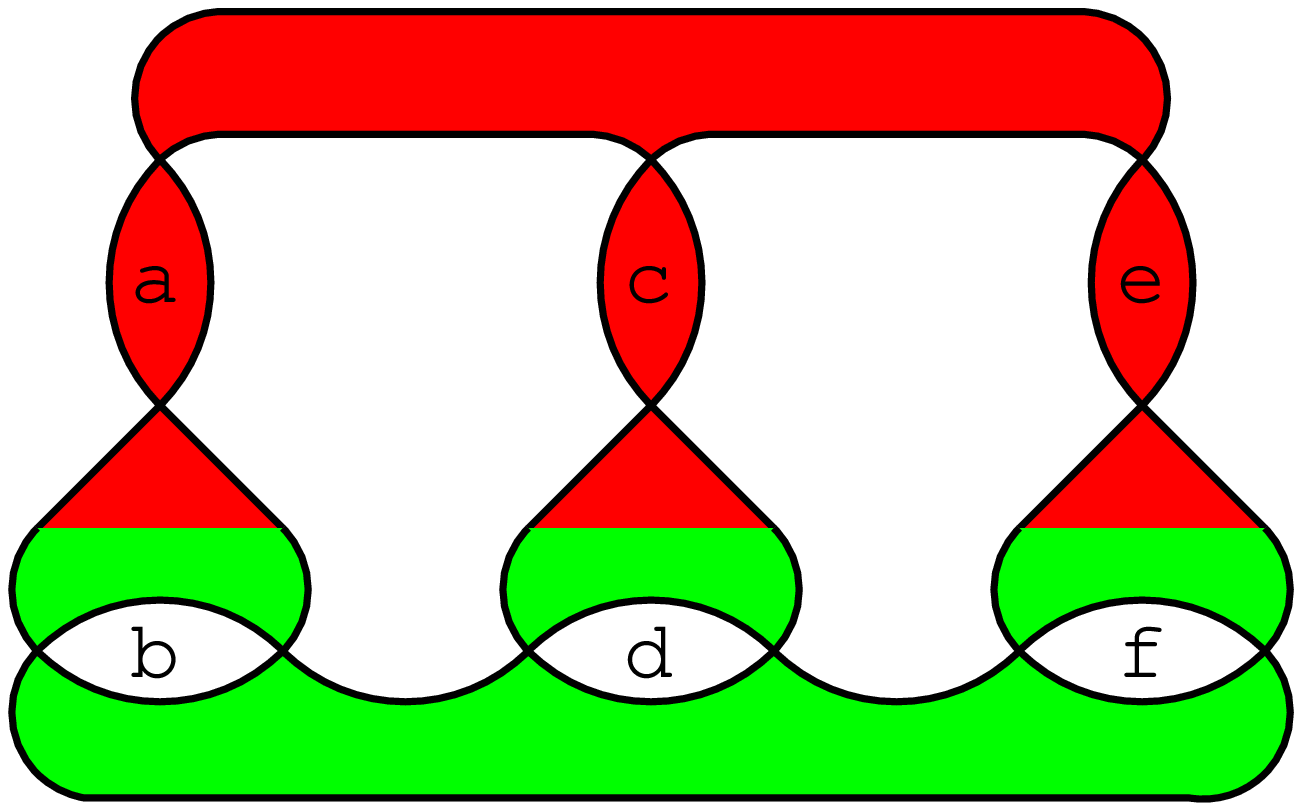}}

\caption{ \ }
$$
\left( \begin{array}{c}
a_1 a_3 a_5 \\
a_3 a_5 \\
a_5 a_1 \\
a_1 a_3 \\
a_1 + a_3 + a_5
\end{array} \right)^{\rm T}
\left( \begin{array}{ccccc}
0 & 0 & 0 & 0 & 1 \\
0 & 0 & 1 & 1 & 0 \\
0 & 1 & 0 & 1 & 0 \\
0 & 1 & 1 & 0 & 0 \\
1 & 0 & 0 & 0 & 0
\end{array} \right)
\left( \begin{array}{c}
 1\\
a_2 \\
a_4 \\
a_6 \\
a_2 a_4 + a_4 a_6 + a_6 a_2
\end{array} \right)
$$

\end{figure}


\begin{figure}
\setcounter{subsection}{7}

\subsection{Families of knots with seed the Borromean link $C_2^3$ with sixteen terms (seven cases)}

\centering

\psfrag{a}{\LARGE{$a_5$}}
\psfrag{b}{\LARGE{$a_2$}}
\psfrag{c}{\LARGE{$a_3$}}
\psfrag{d}{\LARGE{$a_4$}}
\psfrag{e}{\LARGE{$a_6$}}
\psfrag{f}{\LARGE{$a_1$}}

\scalebox{0.4}{\includegraphics{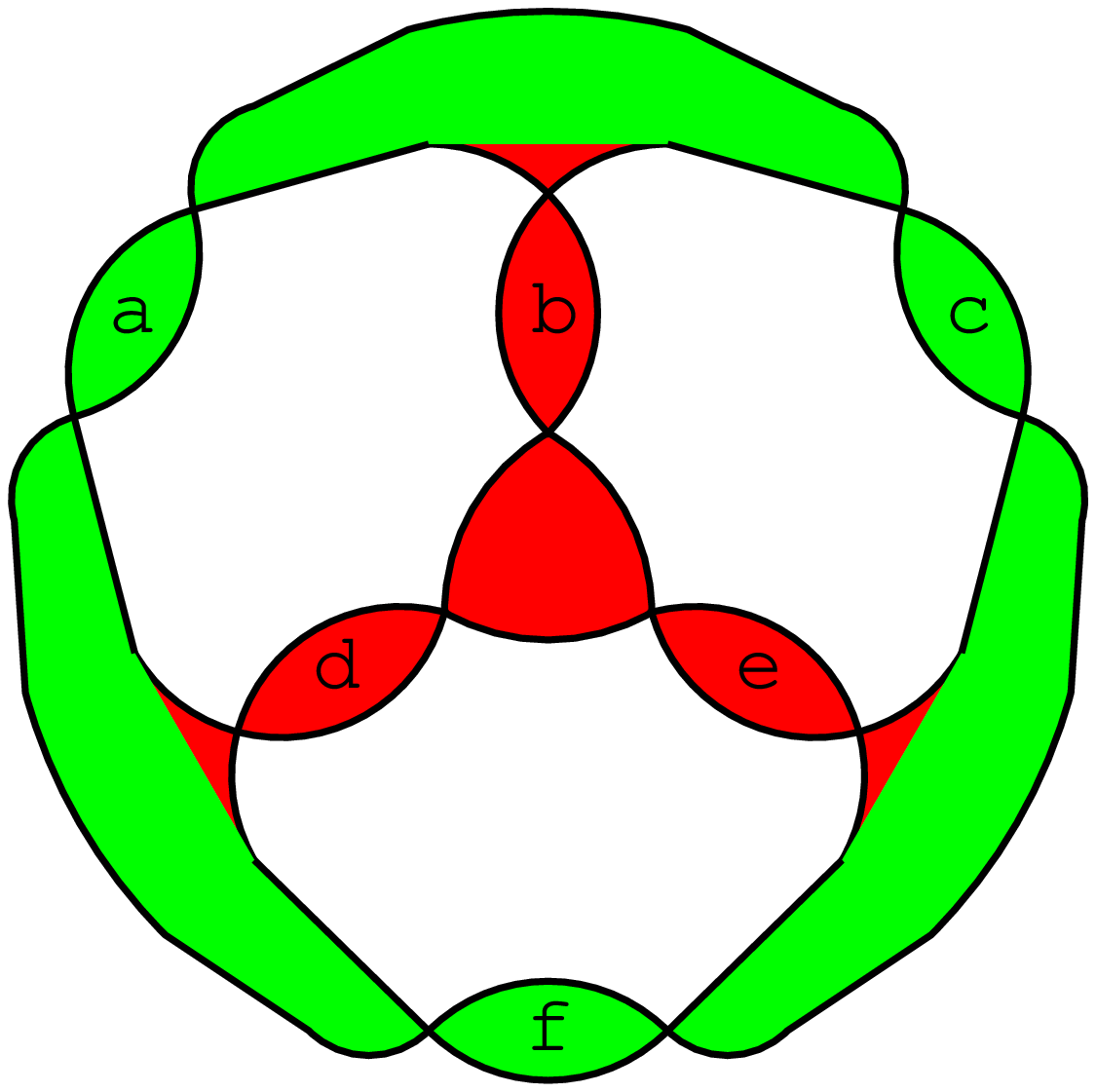}}
\caption{ \ }
$$
\left( \begin{array}{c}
a_1 + a_3 + a_5 \\
a_3 a_5 \\
a_5 a_1 \\
a_1 a_3 \\
a_1 a_3 a_5
\end{array} \right)^{\rm T}
\left( \begin{array}{ccccc}
0 & 0 & 0 & 0 & 1 \\
0 & 0 & 1 & 1 & 0 \\
0 & 1 & 0 & 1 & 0 \\
0 & 1 & 1 & 0 & 0 \\
1 & 0 & 0 & 0 & 0
\end{array} \right)
\left( \begin{array}{c}
1 \\
a_2 \\
a_4 \\
a_6 \\
a_2 a_4 + a_4 a_6 + a_6 a_2
\end{array} \right)
$$

\end{figure}

\begin{figure}

\centering

\psfrag{a}{\LARGE{$a_5$}}
\psfrag{b}{\LARGE{$a_2$}}
\psfrag{c}{\LARGE{$a_3$}}
\psfrag{d}{\LARGE{$a_4$}}
\psfrag{e}{\LARGE{$a_6$}}
\psfrag{f}{\LARGE{$a_1$}}

\scalebox{0.4}{\includegraphics{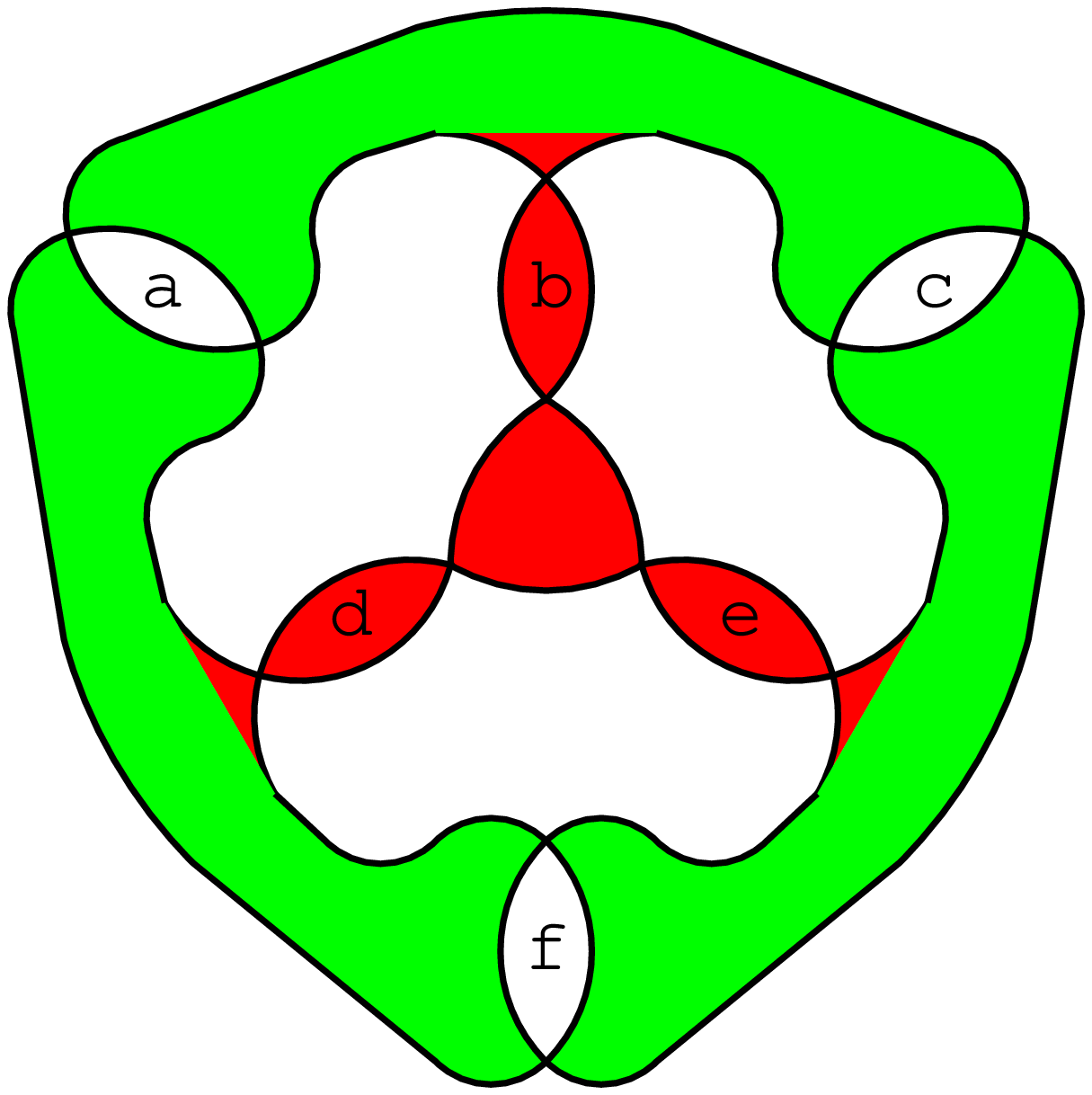}}
\caption{ \ }
$$
\left( \begin{array}{c}
a_1 a_3 + a_3 a_5 + a_5 a_1 \\
a_1 \\
a_3 \\
a_5 \\
1
\end{array} \right)^{\rm T}
\left( \begin{array}{ccccc}
0 & 0 & 0 & 0 & 1 \\
0 & 0 & 1 & 1 & 0 \\
0 & 1 & 0 & 1 & 0 \\
0 & 1 & 1 & 0 & 0 \\
1 & 0 & 0 & 0 & 0
\end{array} \right)
\left( \begin{array}{c}
1 \\
a_2 \\
a_4 \\
a_6 \\
a_2 a_4 + a_4 a_6 + a_6 a_2
\end{array} \right)
$$

\end{figure}

\begin{figure}

\centering

\psfrag{a}{\LARGE{$a_5$}}
\psfrag{b}{\LARGE{$a_2$}}
\psfrag{c}{\LARGE{$a_3$}}
\psfrag{d}{\LARGE{$a_4$}}
\psfrag{e}{\LARGE{$a_6$}}
\psfrag{f}{\LARGE{$a_1$}}

\scalebox{0.4}{\includegraphics{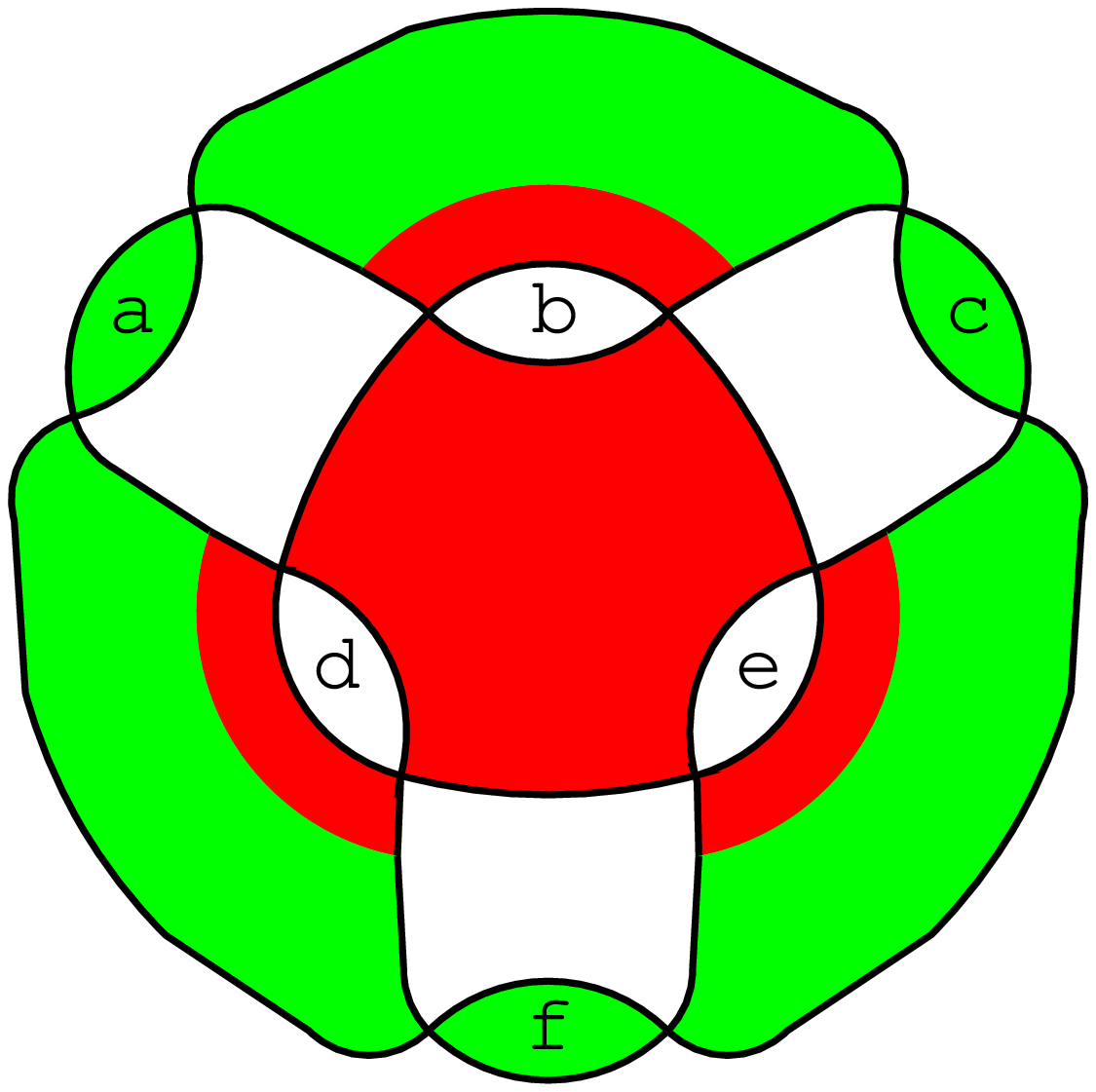}}
\caption{ \ }
$$
\left( \begin{array}{c}
a_1 + a_3 + a_5 \\
a_3 a_5 \\
a_5 a_1 \\
a_1 a_3 \\
a_1 a_3 a_5
\end{array} \right)^{\rm T}
\left( \begin{array}{ccccc}
0 & 0 & 0 & 0 & 1 \\
0 & 0 & 1 & 1 & 0 \\
0 & 1 & 0 & 1 & 0 \\
0 & 1 & 1 & 0 & 0 \\
1 & 0 & 0 & 0 & 0
\end{array} \right)
\left( \begin{array}{c}
a_2 a_4 a_6 \\
a_4 a_6 \\
a_6 a_2 \\
a_2 a_4 \\
a_2 + a_4 + a_6
\end{array} \right)
$$

\end{figure}

\begin{figure}

\centering

\psfrag{a}{\LARGE{$a_5$}}
\psfrag{b}{\LARGE{$a_2$}}
\psfrag{c}{\LARGE{$a_3$}}
\psfrag{d}{\LARGE{$a_4$}}
\psfrag{e}{\LARGE{$a_6$}}
\psfrag{f}{\LARGE{$a_1$}}

\scalebox{0.4}{\includegraphics{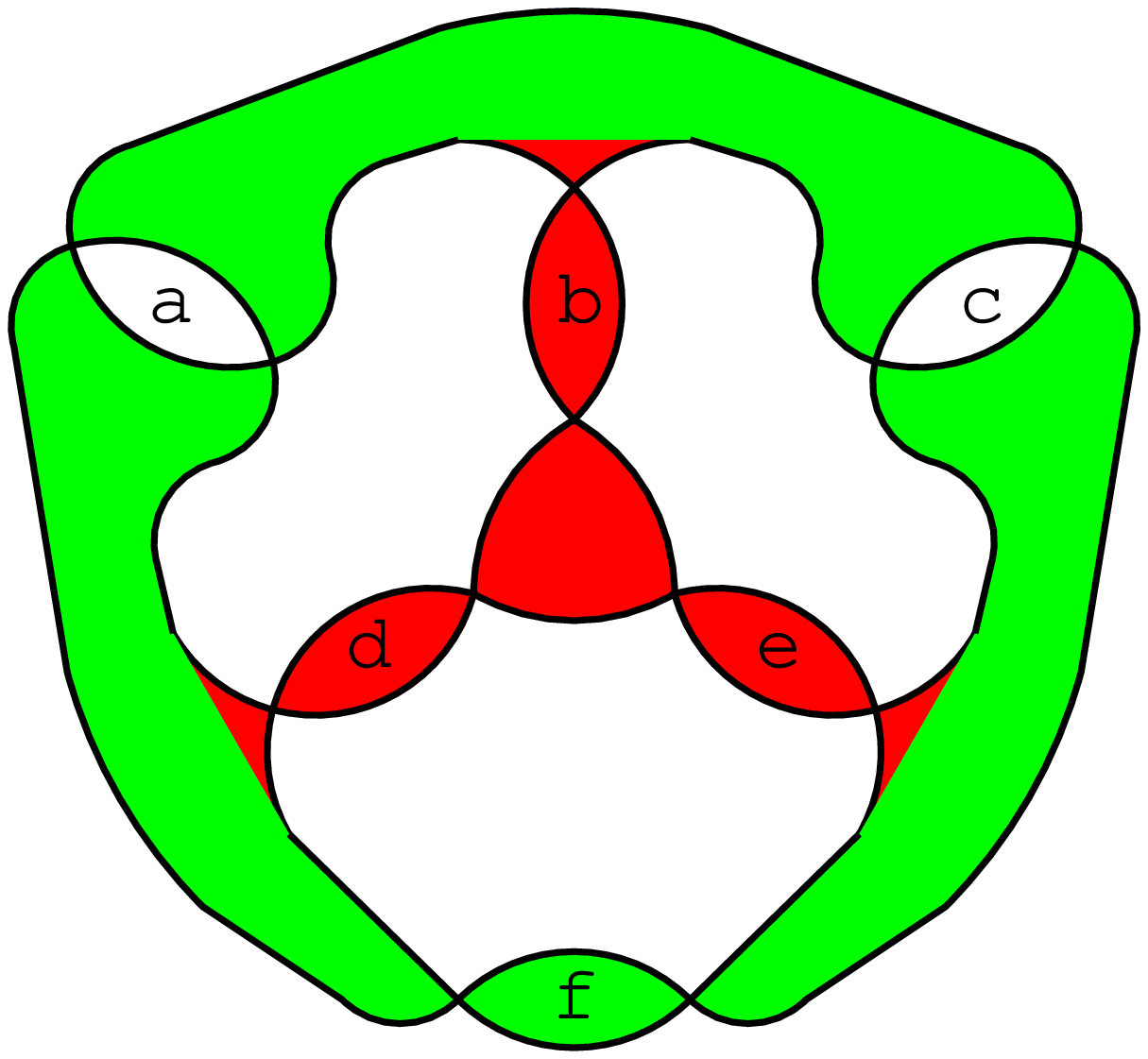}}
\caption{ \ }
$$
\left( \begin{array}{c}
a_1 a_3 a_5 + a_3 + a_5 \\
1 \\
a_1 a_3 \\
a_1 a_5 \\
a_1
\end{array} \right)^{\rm T}
\left( \begin{array}{ccccc}
0 & 0 & 0 & 0 & 1 \\
0 & 0 & 1 & 1 & 0 \\
0 & 1 & 0 & 1 & 0 \\
0 & 1 & 1 & 0 & 0 \\
1 & 0 & 0 & 0 & 0
\end{array} \right)
\left( \begin{array}{c}
1 \\
a_2 \\
a_4 \\
a_6 \\
a_2 a_4 + a_4 a_6 + a_6 a_2
\end{array} \right) =
$$
$$
\left( \begin{array}{c}
1 + a_2 a_3 + a_2 a_5\\
a_2 a_3 a_5 \\
a_5 \\
a_3 \\
a_3 a_5
\end{array} \right)^{\rm T}
\left( \begin{array}{ccccc}
0 & 0 & 0 & 0 & 1 \\
0 & 0 & 1 & 1 & 0 \\
0 & 1 & 0 & 1 & 0 \\
0 & 1 & 1 & 0 & 0 \\
1 & 0 & 0 & 0 & 0
\end{array} \right)
\left( \begin{array}{c}
a_1 a_4 a_6 \\
a_4 a_6 \\
a_6 a_1 \\
a_1 a_4 \\
a_1 + a_4 + a_6
\end{array} \right)
$$
\end{figure}

\begin{figure}

\centering

\psfrag{a}{\LARGE{$a_5$}}
\psfrag{b}{\LARGE{$a_2$}}
\psfrag{c}{\LARGE{$a_3$}}
\psfrag{d}{\LARGE{$a_4$}}
\psfrag{e}{\LARGE{$a_6$}}
\psfrag{f}{\LARGE{$a_1$}}

\scalebox{0.4}{\includegraphics{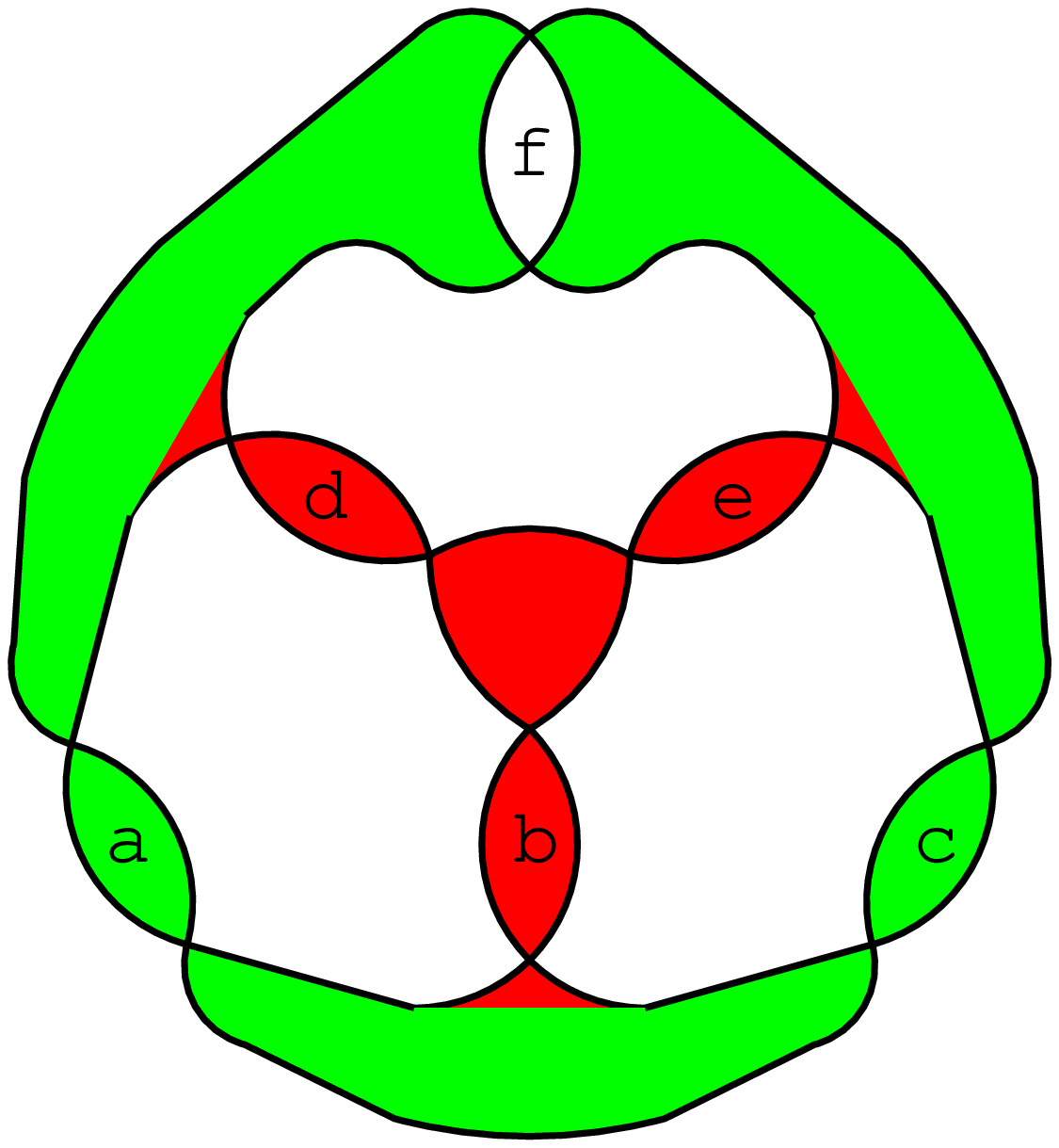}}
\caption{ \ }
$$
\left( \begin{array}{c}
a_1 a_3 + a_1 a_5 + 1 \\
a_1 a_3 a_5 \\
a_5 \\
a_3 \\
a_3 a_5
\end{array} \right)^{\rm T}
\left( \begin{array}{ccccc}
0 & 0 & 0 & 0 & 1 \\
0 & 0 & 1 & 1 & 0 \\
0 & 1 & 0 & 1 & 0 \\
0 & 1 & 1 & 0 & 0 \\
1 & 0 & 0 & 0 & 0
\end{array} \right)
\left( \begin{array}{c}
1 \\
a_2 \\
a_4 \\
a_6 \\
a_2 a_4 + a_4 a_6 + a_6 a_2
\end{array} \right) =
$$
$$
\left( \begin{array}{c}
a_1 a_3 a_6 + a_3 + a_6\\
1 \\
a_1 a_3 \\
a_1 a_6 \\
a_1
\end{array} \right)^{\rm T}
\left( \begin{array}{ccccc}
0 & 0 & 0 & 0 & 1 \\
0 & 0 & 1 & 1 & 0 \\
0 & 1 & 0 & 1 & 0 \\
0 & 1 & 1 & 0 & 0 \\
1 & 0 & 0 & 0 & 0
\end{array} \right)
\left( \begin{array}{c}
a_2 a_4 a_5 \\
a_4 a_5 \\
a_5 a_2 \\
a_2 a_4 \\
a_2 + a_4 + a_5
\end{array} \right)
$$
\end{figure}

\begin{figure}

\centering

\psfrag{a}{\LARGE{$a_6$}}
\psfrag{b}{\LARGE{$a_1$}}
\psfrag{c}{\LARGE{$a_5$}}
\psfrag{d}{\LARGE{$a_2$}}
\psfrag{e}{\LARGE{$a_3$}}
\psfrag{f}{\LARGE{$a_4$}}

\scalebox{0.4}{\includegraphics{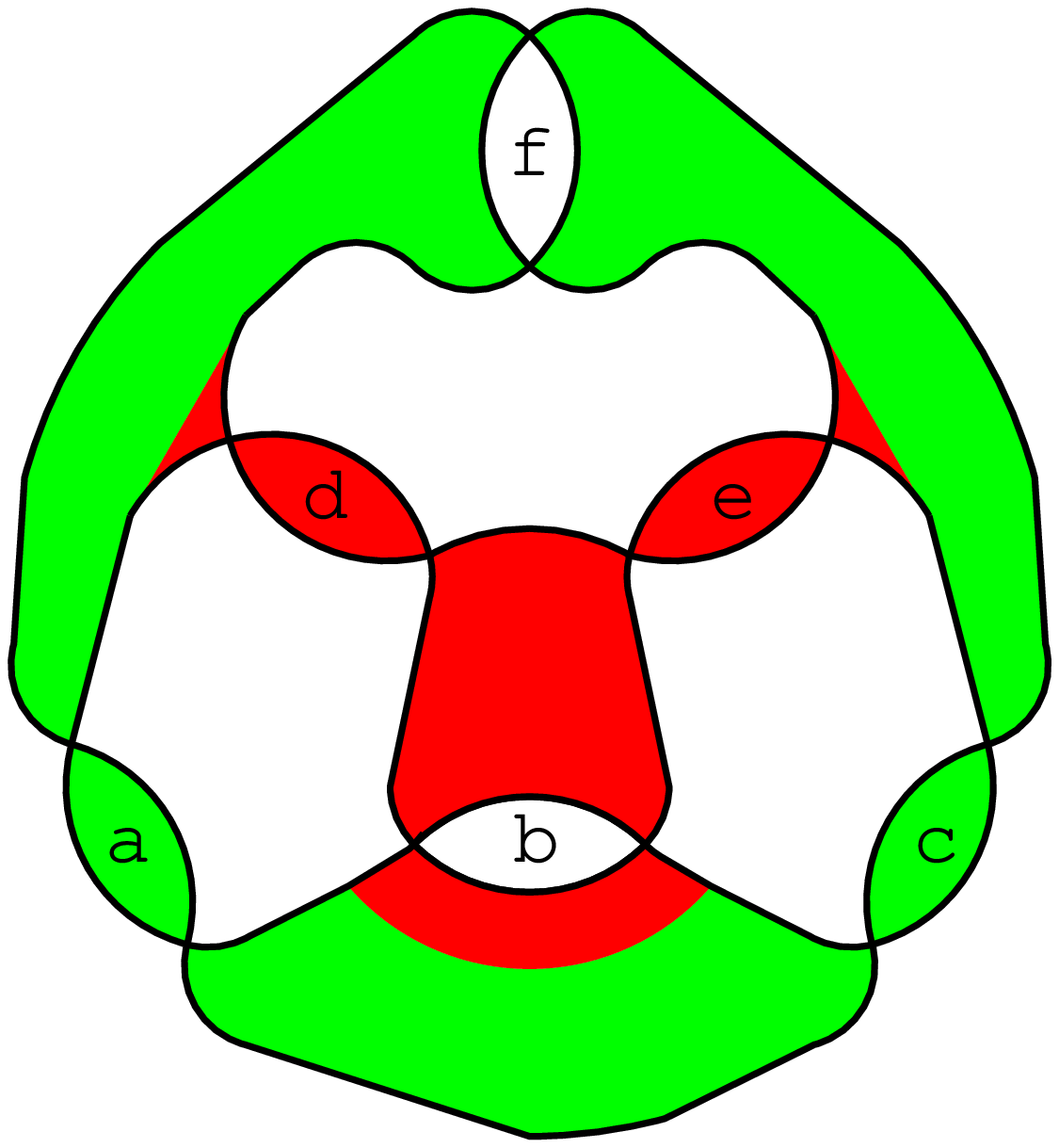}}
\caption{ \ }
$$
\left( \begin{array}{c}
a_4 a_5 + a_4 a_6 + 1 \\
a_5 \\
a_4 a_5 a_6 \\
a_6 \\
a_5 a_6
\end{array} \right)^{\rm T}
\left( \begin{array}{ccccc}
0 & 0 & 0 & 0 & 1 \\
0 & 0 & 1 & 1 & 0 \\
0 & 1 & 0 & 1 & 0 \\
0 & 1 & 1 & 0 & 0 \\
1 & 0 & 0 & 0 & 0
\end{array} \right) \left( \begin{array}{c}
a_1 \\
a_1 a_3 \\
1 \\
a_1 a_2 \\
a_1 a_2 a_3 + a_2 + a_3
\end{array} \right)
$$
\end{figure}

\begin{figure}

\centering

\psfrag{a}{\LARGE{$a_1$}}
\psfrag{b}{\LARGE{$a_2$}}
\psfrag{c}{\LARGE{$a_3$}}
\psfrag{d}{\LARGE{$a_4$}}
\psfrag{e}{\LARGE{$a_5$}}
\psfrag{f}{\LARGE{$a_6$}}

\scalebox{0.6}{\includegraphics{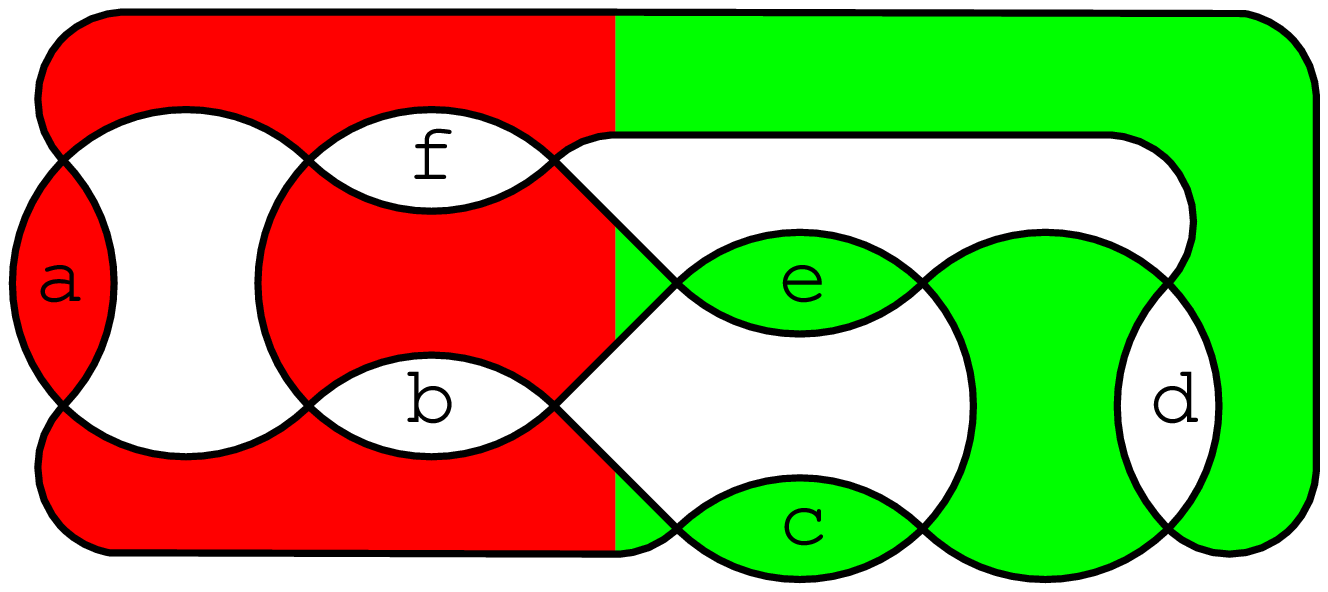}}
\caption{ \ }
$$
\left( \begin{array}{c}
a_1 a_2 a_6 + a_2 + a_6\\
a_1 a_2 \\
a_1 a_6 \\
1 \\
a_1
\end{array} \right)^{\rm T}
\left( \begin{array}{ccccc}
0 & 0 & 0 & 0 & 1 \\
0 & 0 & 1 & 1 & 0 \\
0 & 1 & 0 & 1 & 0 \\
0 & 1 & 1 & 0 & 0 \\
1 & 0 & 0 & 0 & 0
\end{array} \right)
\left( \begin{array}{c}
a_4 \\
1 \\
a_3 a_4 \\
a_4 a_5 \\
a_3 a_4 a_5 + a_3 + a_5
\end{array} \right) =
$$
$$
\left( \begin{array}{c}
a_2 a_5 + a_2 a_3 + 1 \\
a_2 a_3 a_5 \\
a_3 \\
a_5 \\
a_3 a_5
\end{array} \right)^{\rm T}
\left( \begin{array}{ccccc}
0 & 0 & 0 & 0 & 1 \\
0 & 0 & 1 & 1 & 0 \\
0 & 1 & 0 & 1 & 0 \\
0 & 1 & 1 & 0 & 0 \\
1 & 0 & 0 & 0 & 0
\end{array} \right)
\left( \begin{array}{c}
a_4 a_6 \\
a_6 \\
a_1 a_4 a_6 \\
a_4 \\
a_1 a_4 + a_1 a_6 + 1
\end{array} \right)
$$

\end{figure}

\end{document}